\documentclass[12pt]{amsproc}
\usepackage{tikz}
\usepackage{amsmath}
\usepackage{amssymb}
\usepackage{amsthm}
\usepackage{stmaryrd}
\usepackage{cases}
\usepackage{graphicx}
\usepackage{upgreek}
\usepackage{imakeidx}
\usepackage{mathrsfs}

\makeindex

\usepackage{xcolor}
\definecolor{MyGrey}{rgb}{.804,.804,.756}

\newtheorem{theorem}{Theorem}
\newtheorem{lemma}[theorem]{Lemma}
\newtheorem{corollary}[theorem]{Corollary}
\theoremstyle{definition}

\theoremstyle{remark}

\numberwithin{equation}{section}

\newcommand{\C}{\mathbb{C}}
\newcommand{\Z}{\mathbb{Z}}
\newcommand{\R}{\mathbb{R}}

\newcommand{\X}{{\mbox{\sf X}}}

\begin{document}
\title{Braiding and Folding Branched Covers}

\author{J. Scott Carter}
\address{University of South Alabama\\ Department of Mathematics and Statistics\\Mobile, AL}
\curraddr{}
\email{carter@southalabama.edu}
\thanks{}

\author{Seonmi Choi}
\address{Seowon University\\Department of Mathematics Education\\
Cheongju, 28674, Republic of Korea}
\curraddr{}
\email{smchoi@seowon.ac.kr}
\thanks{}

\author{Byeorhi Kim}
\address{Pohang University of Science and Technology\\Center for Research in Topology\\ Pohang, Republic of Korea}
\curraddr{}
\email{byeorhikim@postech.ac.kr}
\thanks{}

\begin{abstract}
A {\it braiding} of a branched cover of the $3$-sphere that is branched over a knot is a continuous map of the cover into the product of the sphere with a $2$-dimensional disk that has the property that the projection onto the sphere factor induces the covering. A {\it folding} maps the cover in to the product of a sphere with an interval. Moreover, away from the branch set, the map is a general position immersion. Cyclic branched covers can be braided so that the map is an embedding when the disk factor is $2$-dimensional. Dihedral branched covers can  be folded. 

In as much as possible, the braidings and foldings that are presented are quite detailed. In particular, the paper focuses upon a folding of the dihedral cover of the $3$-sphere that is branched along a torus knot of type $(2,5)$. The cover also is homeomorphic to the $3$-sphere. 
\end{abstract}

\maketitle

\begin{center}{\it To  Seiichi Kamada on the occasion of his 60th birthday. And to  both Lou Kauffman  and Cameron Gordon on the occasion of their 80th birthdays.} 
\end{center}

\section{Introduction}

Consider  the complex valued function $g(z) = z^2$ with its domain restricted to the unit disk $ D^2= \{ r e^{ {\boldsymbol {i}} \theta}: 0 \le r \le 1, \ \theta \in [0, 2\pi) \}$.  The range of the function is also $D^2$, and each point in the range, with the exception of $0$,  has two points in its pre-image. The squaring function gives the prototypical example of a simple branched cover. 
More generally, restrict the $n$\/th power function $z \mapsto z^n$ to the unit disk $D^2$. Each point in the disk $D^2$ except $0$ has $n$ points in its pre-image. And this function gives the prototypical example of an $n$-fold cyclic branched cover.  

The $n$\/th power function, for $2\le n$, extends to the Riemann sphere $S^2= \C \cup \{ \infty \}$. There are two cyclic $n$-fold branch points: one is at $0$ the other is at $\infty$. In case $n=2$, this is called a {\it simple branched cover of the $2$-sphere with two simple branch points}.  The general case is called the {\it $n$-fold cyclic branched cover of $S^2$ with two branch points.} In general, cyclic branched covers will be built by extrapolating these examples. 

Two extrapolations come to mind immediately. The first is the case of cyclic branched covers of the $2$-sphere $S^2$ in which there are an arbitrary even number of branch points. The second is the case of cyclic branched covers of the $3$-sphere $S^3$ in which the branch set is an arbitrary smooth knot or link. Both cases are described in a topological fashion as follows.

The branch set is  an oriented codimension-two  smooth submanifold, denoted $K$, that has a tubular neighborhood. Take $n$ copies of the sphere in question. Remove an open tubular neighborhood of the branch set in each; 
cut each along a  codimension-one submanifold ({\it Seifert submanifold}) that is bounded by the branch set. 
Sew or glue the negative side of the Seifert submanifold in the $j$\/th sphere to the positive side of the Seifert submanifold in the $(j+1)$\/st sphere cyclically. Finally into the neighborhoods that were removed, glue the branch set $K$ times  the graph of $z \mapsto z^n$ restricted to the unit disk.  The last product mentioned is a cartesian product. The graph in question is a subset of the $4$-dimensional space $D^2 \times D^2$. 

There are, of course, questions about uniqueness of the construction that was outlined above. In particular, the cyclic branched covers seem to depend upon the Seifert manifold that was used for cutting and pasting. We refer the reader to the papers  \cite{EF},\cite{HildenPac78}, \cite{HLMTrans83}, and \cite{Kolay} for more details about branched coverings.

Another construction of branched coverings is given by means of representations into a finite group. Let us elaborate the construction above with some more rigor.

Consider a  smoothly embedded, compact, codimension-two submanifold $K$  of the $m$-sphere. In case $m=3$, the manifold is a knot or link, and in case $m=2$, it is a finite set of points. In either case, there is a tubular neighborhood, $N(K)$, that is  also embedded. We consider the compact $m$-manifold with boundary   $E = E(K) = S^m \setminus {\mbox{\rm int}}(N(K))$. If $X = S^m \setminus K$, then $X$ deformation retracts onto $E$, and, in  case $m=3$, we colloquially call either set {\it the knot or link exterior}.   Depending upon  the value of $m$, the boundary of $E$ is a collection of  circles or tori.

Suppose that there is a non-trivial  group homomorphism $G \stackrel{\phi}{\longleftarrow} \pi_1(E)$ 
from the fundamental group of the exterior into a finite group $G.$\footnote{Functions are indicated by left-pointing arrows since the result of a composition $Z \stackrel{g}{\leftarrow} Y \stackrel{f}{\leftarrow} X$,  is written as $g(f(x))$. Our convention differs from the standard one, but it has its conveniences.}  We think of $G$ as a subgroup of the symmetric group. For the purposes of this paper,  the group $G$ will be specialized. If $G$ is a cyclic group, $\Z /n$, then the generator is thought of as the $n$-cycle $(1,2,\ldots ,n) \in \Sigma_n$ or its inverse; the homomorphism passes through the abelianization homomorphism.  In case $G = D_n$ --- the dihedral group of order $2n$, then $D_n$ is a subgroup of the  permutation group $\Sigma_n$ on $n$ elements since  the dihedral group is  the group of symmetries of a regular polygon that has $n$ vertices.  
The techniques that we discuss here also extend to the alternating groups $A_3$ and $A_5$, but we won't cover these examples herein.

 There is an associated $n$-fold irregular cover of  the exterior $E$ that is determined by $\phi$; see e.g.  \cite{HildenPac78}. Assume that the permutation representation of $G$ lies in the permutation group $\Sigma_n$. Then take $n$ copies of the exterior $E(K)$, cut each along a collection of  proper submanifolds that are dual to the fundamental groups elements, and re-glue these cut manifolds by following the permutation representation. In both the $2$-dimensional and $3$-dimensional cases, there is quite a bit of ambiguity that requires further detailed explanations. Much of the content herein contain such amplifications. 
  
 A branched cover $M$  of the sphere branched over $K$ is obtained by extending this cover into the tubular neighborhood in a standard manner. In the case of a cyclic group, that extension was described via the $n$\/th power function. In the dihedral  case, the extension to the tubular neighborhood is  more subtle. We ask for the reader's patience. Details will be presented when the text is ripe.

\subsection{Main Results}
\begin{sloppypar}
Let $M=M^m$ denote a cyclic  or dihedral branched cover of the $m$-sphere, $S^m$,  that is branched over a codimension-two submanifold  $K$ of $S^m$.   {\it A braiding of $M$}  (resp. {\it a folding of $M$})  is a continuous map $S^m \times D \stackrel{f}{\longleftarrow} M$ that satisfies conditions $1$ and $2$ when $D= [0, n+1] \times [-1,1]$ (resp. when $D=[0,n+1]$). Here 
the dihedral group $D_n$ is of order $2n$, or in the case of the cyclic branched cover it is the length of the cycle. 
\end{sloppypar}

\begin{enumerate}
\item
 Away from a neighborhood of the branch set, the map $f$ is a general position immersion, and 
 \item
  The composition $S^m \stackrel{p_1}{\longleftarrow}S^m \times D \stackrel{f}{\longleftarrow} M$ of $f$ with the projection onto the  first factor induces the branched covering. 
  \end{enumerate}
  We also will say that the branched cover $M$ is {\it braided} or {\it folded}.  In practice, $m=1, 2,$ or $3$. In the course of our study, we will indicate that a folded branched cover always lifts to an immersed braided branched cover. Whether or not it can be embedded, is unknown. See~\cite{CKFold} or \cite{Kolay} for details and examples of immersed braidings of branch covers that do not embed to embedded braidings. 
  
Terminological differences are emphasized. 
  
  \begin{itemize}
  \item
  A folding maps the branched cover $M$ of $S^m$ into $S^m \times [0,n+1]$. Away from the branch set the image is {\it a general position immersion}. In case $m=1$, there are double points that resemble an $\X$: $\{ (x,y) \in \R^2: xy=0\}$. In case $m=2$, the double points are of the form of the intersection of the two coordinate planes in $3$-space: $\{ (x,y,z) \in \R^3: xy=0\}$ . The triple points resemble the intersection of  the $3$ coordinate planes at the origin: $\{ (x,y,z) \in \R^3: xyz=0\}$. In case $m=3$, the double points are of the form: $\{ (x,y,z,w) \in \R^4: xy=0\}$. The triple points are of the form: $\{ (x,y,z,w) \in \R^4: xyz=0\}$. And there are quadruple points of the form: $\{ (x,y,z,w) \in \R^3: xyzw=0\}$. That is every multiple point  has a neighborhood in which the image is homeomorphic to one of these coordinate sets. 
  
  \item
 An immersed braiding is a lifting of  a folding of $M$ into $S^m \times[0, n+1] \times [-1,1]$. In case, $m=2$ the branched cover $M$ of $S^2$ is immersed into the $4$-dimensional space $S^2 \times[0, n+1] \times [-1,1]$. As such, it has isolated general position double points. In case $m=3$, the general position double point set of the immersed braiding will form a closed $1$-dimensional manifold. 
 
 \item
 An embedded braiding is an immersed braiding that has no double points. Embedded braidings are usually called {\it braidings}.
 
 \end{itemize}

In the case $m=1$, the branch set is empty. An embedded braiding is merely a braiding.   
A folding occurs when a closed braid diagram projects onto the annulus upon which the diagram would be drawn. In this case, the braid, as an element in the braid group, projects to a permutation that is presented as intertwining arcs in a disk region of the annulus. These arc close to loops as in the braid case.

The dihedral group $D_n$ of order $2n$ is a subgroup of the permutation group $\Sigma_n$. The annular region into which $M$ is folded is of the form $S^m \times [0, n+1].$ The permutations appear as strings permuting the integral points in the interval factor.

{\it Thus the idea of braiding or folding a branched cover extrapolates the idea of a braid. It also extrapolates Kamada's notion of a surface braid. } 
See \cite{CKFold} and \cite{Kolay} for further information.

We begin by considering branched covers of $S^2$ that are branched over a finite set of points. Let $b$ denote a positive integer, $1<b$. Let $K=\{p_1, p_2, \ldots, p_b\}$ denote a finite set of points that are contained in $S^2$. Let $N(K)$ denote a union of sufficiently small disks that form an embedded tubular neighborhood of $K$. 
The exterior $E(K)= S^2 \setminus {\mbox{\rm int}} (N(K))$ is a planar surface whose fundamental group $\pi_1(E(K))$ is a free group $F_\ell$ of rank $\ell=b-1$. 

\begin{theorem}
\label{Dihedral2}
In the notation above,  let  $ D_n \stackrel{\phi}{\longleftarrow} \pi_1(E(K)) = F_{\ell}$ denote a group homomorphism into the dihedral group of order $2n$. Let $M$ denote the associated branched cover of $S^2$ that is branched along $K$. 
There is an immersed  braiding 
\[ S^2 \times \left([0, n+1] \times [-1,1] \right) \stackrel{f \ \  \ }{\longleftarrow M}.\]
In particular,  the  restriction of  continuous, general position map $f$  to the exterior $E(K)$ is an immersion, and the composition 
\[ S^2 \stackrel{p_1}{\longleftarrow} S^2 \times \left([0, n+1] \times [-1,1]\right) \stackrel{f \ \ \ \ }{\longleftarrow M} \] is the branched covering map.
\end{theorem}

\begin{theorem}
\label{CyclicThm}
Let $M$ denote an $n$-fold cyclic branched cover of the $3$-sphere $S^3$ that is branched over a  knot or link $K \subset S^3$. There is a  braiding 
\[ S^3 \times \left([0, n+1] \times [-1,1] \right) \stackrel{f \ \  \ }{\longleftarrow M} .\]
In particular,  the  restriction of  continuous, generic map $f$  to the exterior $E(K)$ is an embedding, and the composition 
\[ S^3 \stackrel{p_1}{\longleftarrow} S^3 \times \left([0, n+1] \times [-1,1]\right) \stackrel{f \ \ \ \ }{\longleftarrow M} \] is the branched covering map.
\end{theorem}

These results also appear in~ \cite{HildenPac78}. But our context is more diagrammatic, and it uses more modern technology. We explicitly demonstrate the foldings and braidings by means of  permutation and braid charts. In the $3$-dimensional case, the folding is described via a sequence of permutation charts.

\begin{theorem}
\label{DihedralThm}
Let $M$ denote a dihedral branched cover of $S^3$ that is branched over a knot or link  $K \subset S^3$.   The cover is induced by a homomorphism 
$D_n \stackrel{\phi}{\longleftarrow} \pi_1(E(K))$ from the  exterior of $K$ into the dihedral group. Suppose that $1<n=2k+1$ is an odd positive integer. So $1\le k$.

There is a  folding 
\[ S^3 \times [0, n+1]  \stackrel{f \ \  \ }{\longleftarrow M} .\]
So the composition 
\[ S^3 \stackrel{p_1}{\longleftarrow} S^3 \times [0, n+1] \stackrel{f \ \ \ \ }{\longleftarrow M} \] is the branched covering map.
\end{theorem}

\begin{sloppypar}
Because every $3$-manifold can be expressed as a $3$-fold simple branched cover of the $3$-sphere that is branched over a knot or link~\cite{HildenBull}, the focus of the work~\cite{HLMTrans83} is upon  $3$-fold dihedral covers. Indeed, the work or the first named author  with Kamada in~\cite{CKFold} and~\cite{CarterDA} concentrated on such dihedral covers.  
The proofs and techniques developed here expand upon the initial work in ~\cite{CKFold}. 
\end{sloppypar}

\subsection{Presentation of the material}

The main theorems have been stated and put into context in the material above. The theorems indicate that   branched covers of spheres can be mapped continuously into a sphere times a disk in such a way that the projection onto the sphere induces the covering. 

For example, Section~\ref{charts} starts by recalling that classical braids are  examples of such braidings as are surface braids that are obtained from Kamada's braid charts. Closely related to braid charts are permutation charts. These are the main tools for our study. Both  are labeled graphs embedded into a $2$-dimensional disk. Braid charts of degree $n$  give braided surfaces in $S^2\times [0,n+1] \times [-1,1]$. Permutation charts give 
folded surfaces in $S^2\times [0,n+1]$ --- surfaces in $3$-dimensions that are in general position away from the branch set.

Section~\ref{diag4dih} contains a significant portion of the techniques that will be used to establish the results of the paper. Some of this section might be considered to be elementary, but we would rather bore an expert than intimidate a student. We hope that we have done neither. The focus is upon the $10$ element dihedral group $D_5.$ The techniques generalize  to arbitrary general dihedral groups. An analogue of braid charts is developed in this section that will be useful in providing blueprints for the dihedral branched covers of both the $2$-dimensional and $3$-dimensional spheres. In the $2$-dimensional case, dihedral charts are labeled oriented graphs that also can be thought of in terms of maps into classifying spaces. We also give a glossary between the chart for a fixed reflection or rotation in $D_5$ and the corresponding  representation via a permutation chart. Surfaces in $3$-space that are in general position in $3$-space away from the branch sets are described by permutation charts as indicated in Fig.~\ref{relations}.

Section~\ref{thesteps} contains the principal and motivating example. The $D_5$ cover of $S^3$ that is branched over the torus knot $T(2,5)$ will be explicitly folded into $S^3 \times [0,6]$ by translating the sequence of dihedral charts that are illustrated in this section into permutation charts. 

A key aspect of this translation is to indicate that certain operations to dihedral charts can be realized in their permutation counterparts. Section~\ref{digress} quantifies the changes between successive charts as invertibility of the relators in the dihedral group when the relators are thought of as morphisms in a category constructed from a group presentation. The differences between successive dihedral charts are consequences of some of the permutation chart moves that are depicted in Fig.~\ref{almostKmoves}. Section~\ref{thelemmas} contains the proofs that those moves can be realized within the permutation charts.  
The proofs of the main results are contained in Sections~\ref{mequals2},
~\ref{ProofofDihedral}, and ~\ref{cyclicproof}. Section~\ref{Concl} is a short recap, and it points out that the techniques herein are more widely applicable. We conclude with our aspirations about our techniques.

\section{Braid and Permutation Charts}
\label{charts}

Seiichi Kamada \cite{KamJKTR} introduced the notion of a braid chart in order to represent orientable surfaces as simple surface braids in $4$-dimensional space. A surface braid is an embedding  of an orientable surface into the cartesian product of $S^2 \times D^2$ such that the projection onto the sphere factor induces a simple branched covering. We recall that {\it simple} branch points are modeled upon the complex quadratic function at the origin.  Kamada's theory of surface braids and their charts is very rich. It provides a method  of studying knotted orientable surfaces from a diagrammatic point of view. See, for example, the many works of Nagase and Shima \cite{NS1,NS2} etc. 

\begin{figure}
\[
\scalebox{.5}{\includegraphics{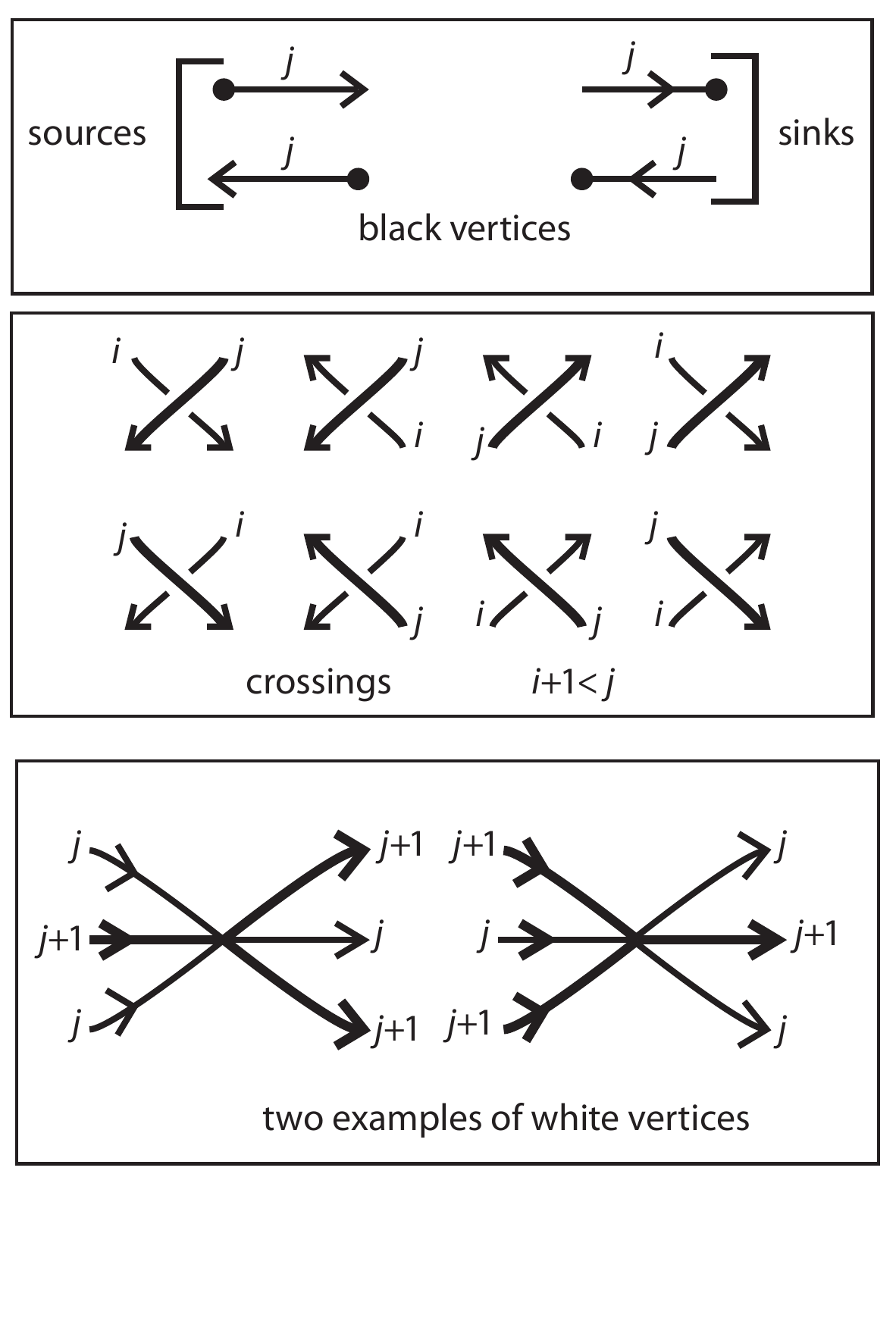}}
\]
\caption{Types of vertices in a braid chart}
\label{BraidChart}
\end{figure}

{\it A braid chart  of degree $n$} is an oriented   graph in which the edges are labeled by integers in the set $\{1,2,\ldots, n-1 \}$ and which is  embedded in the $2$-dimensional disk $D^2= [0,p+1]\times [0, q+1]$.  There are three types of vertices:
\begin{enumerate}
\item
{\it A black vertex} is incident to a single edge that is incoming or outgoing. If the incident edge is labeled with the integer $j \in \{ 1,2, \ldots, n-1 \}$, then the black vertex is said to be of {\it index $j$}.
\item {\it A crossing} is incident to four edges that are labeled cyclically $j,k,j,k$ where $1<|j-k|$. The edges with the same labels are consistently oriented. So one edge that is labeled $j$ points into the vertex and the other that is labeled by $j$ points outward. Similarly for $k$. One imagines that the antipodal edges are flowing through the crossing.
\item {\it A white vertex} has six edges incident. The incident edges are cyclicly labeled $(j, j \pm 1, j, j \pm1, j, j \pm 1)$. Three consecutive edges are incoming and the other three are outgoing. 
\end{enumerate}

Notice that the crossings and the white vertices correspond to the relations in the group presentation of the $n$-string braid group, 
\begin{small}{
\[ B_n = \left\langle \sigma_1, \sigma_2, \ldots, \sigma_{n-1} : \begin{array}{cl} \sigma_j \sigma_k = \sigma_k \sigma_j & {\mbox{\rm for}}\ \  1<|j-k|, \\
\sigma_j \sigma_{j+1} \sigma_j = \sigma_{j+1} \sigma_j \sigma_{j+1} & {\mbox{\rm for}}  \ \ j=1, \ldots, n-2  \end{array} \right\rangle . \]}\end{small}

\begin{figure}
\[
\scalebox{.65}{\includegraphics{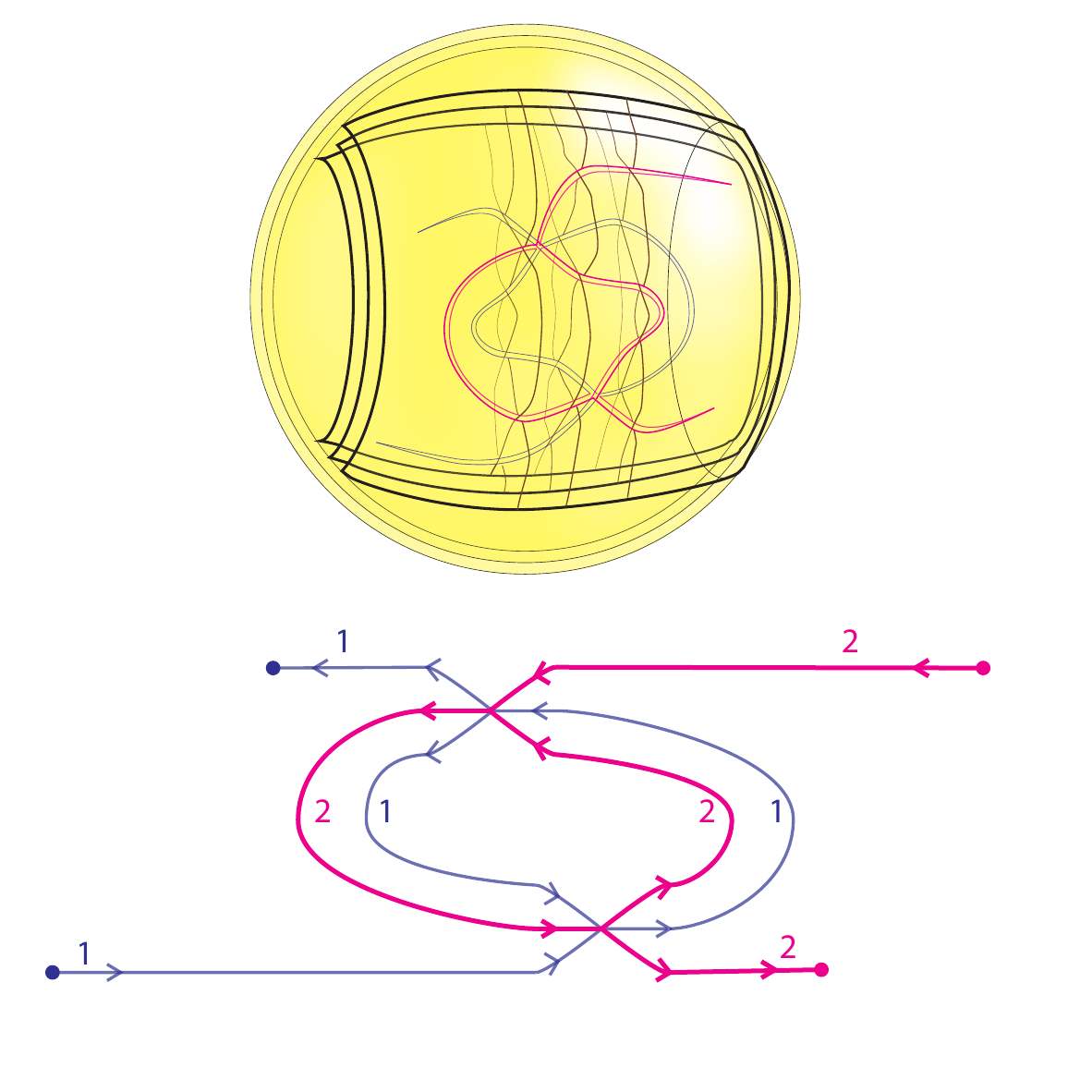}}
\]
\caption{A chart and the associated braided (or folded) surface}
\label{FromCKnotes}
\end{figure}

Figure~\ref{BraidChart} indicates the types of vertices. Our depictions of white vertices  and crossings differs from that given by Kamada, but the content is the same.  Figure~\ref{FromCKnotes} indicates how a chart represents a surface braid that is drawn as a broken surface diagram.

The disk $[0,p+1]\times [0, q+1]$ is arranged as a rectangle with the first coordinate along the horizontal and the second along the vertical.  The integer $p$ indicates the total number of vertices and critical points along the edges. The integer $q$ indicates the maximal number of times an arc perpendicular to the edges intersects the edge set. 
In Fig.~\ref{FromCKnotes}, $p=10$ and $q=6$. That is to say that the chart may be situated in that disk so that its edges are usually pointing in the horizontal direction. The projection onto the interval $[0,p+1]$ acts as a Morse function for the chart. If an edge points right, then the corresponding braid generator is $\sigma_j$. If an edge points left, it corresponds to the inverse $\sigma_j^{-1}. $ A generic critical point corresponds to a type II move, or equivalently the group relation $ \sigma_j^{\pm 1} \sigma_j^{\mp 1} = 1$.  For concision of terminology, the vertices (black, white, and crossings) of the graph are also considered to be critical points. Moreover, often the critical points will be arranged to lie at distinct critical levels.

The disk also is considered to be a sub-disk of the $2$-dimensional sphere $S^2$. So the braid chart describes a braid movie --- a sequence of embedded braids that successively differ by branch points and braid relations.  For each still in the movie, one can form the braid closure. The movie starts and ends with an identity braid that closes to a collection of $n$ nested circles. These circles bound a sequence of disks of increasing radii outside of the realm in which the chart appears. In this way, the braid chart implicitly describes a braiding of an orientable surface in $S^2 \times D^2$. This second disk factor can be parametrized as $D^2 = [0,n+1] \times [-1,1]$ where $n$ is the braid index of the chart, and the crossings of the braid generators occur in the interval $[-1,1]$.  More details can be found in Kamada's book \cite{KamBook}. 

\begin{figure}
\[
\scalebox{.25}{\includegraphics{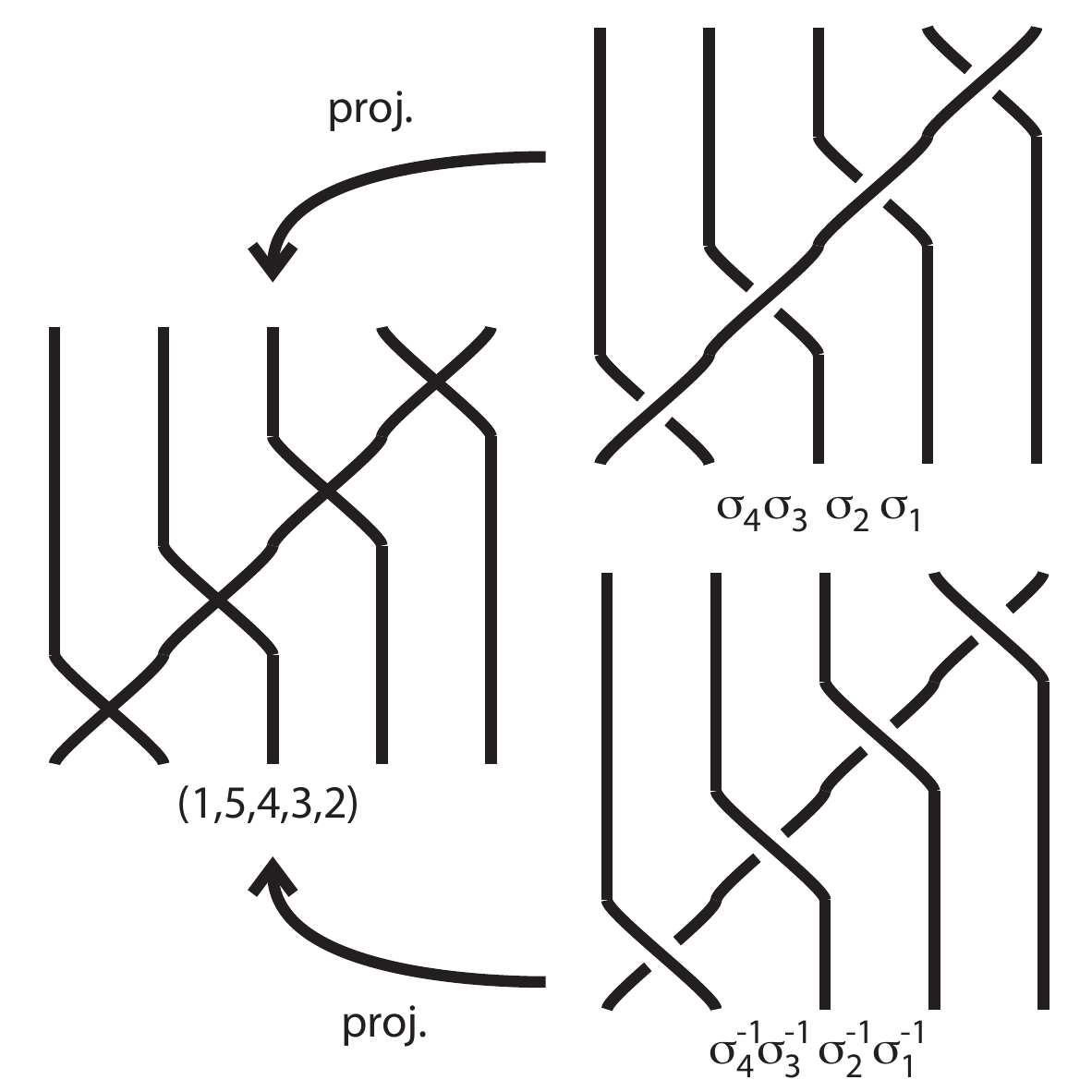}}
\]
\caption{Projections of braids to permuations}
\label{proj}
\end{figure}

The projection $\Sigma_n \stackrel{p}{\longleftarrow} B_n$ from the braid group to the symmetric group  
\begin{small}{
\[ \Sigma_n = \left\langle t_1, t_2, \ldots, t_{n-1} : \begin{array}{cl} t_j t_k = t_k t_j & {\mbox{\rm for}}\ \  1<|j-k|, \\
t_j t_{j+1} t_j = t_{j+1} t_j t_{j+1} & {\mbox{\rm for}}  \ \ j=1,2, \ldots, n-2,  \\
t_j^2 =1 & {\mbox{\rm for}} \ \ j=1,2,\ldots, n-1\end{array} \right\rangle,  \]}\end{small}
that is given by $p(\sigma_i)=t_i$,
can be described by means of projecting the standard braid generators to planar transverse and generic double points as indicated in Fig.~\ref{proj}. 
Similarly, a braid chart projects to a {\it permutation chart} by means of removing the orientations upon the edges of the graph.

\section{Diagrammatics of the dihedral group}
\label{diag4dih}

This section will concentrate upon the $10$ element dihedral group $D_5$ that is given via the presentation:
\[ D_5 = \langle r, x : r^2 = x^5 = (xr)^2 = 1 \rangle. \]
The element $r$ is called {\it a reflection} and the element $x$ is called {\it a rotation}. 
Using the relation $xrxr=1$, we obtain that $xr=r^{-1}x^{-1}$. And since $r^2=x^5=1$, $xr=rx^{-1}=rx^4$. The elements in the group can be classified as rotations $X= \{ 1,x,x^2, x^3, x^4 \}$, and reflections $R=\{ r,rx, rx^2, rx^3, rx^4\}$. Thus the standard normal forms are presented as 
\[ D_5 = \{ 1,x,x^2,x^3,x^4, r, rx, rx^2, rx^3, rx^4 \}.\]

The dihedral group of order $2n$, that is denoted as $D_n$ has a similar presentation with $n$ replacing $5$ above. And almost everything that we say about $D_5$ can be extended to the more general case. We are only interested in the case in which $n=2k+1$ is an odd integer.

\begin{figure}
\[
\scalebox{.4}{\includegraphics{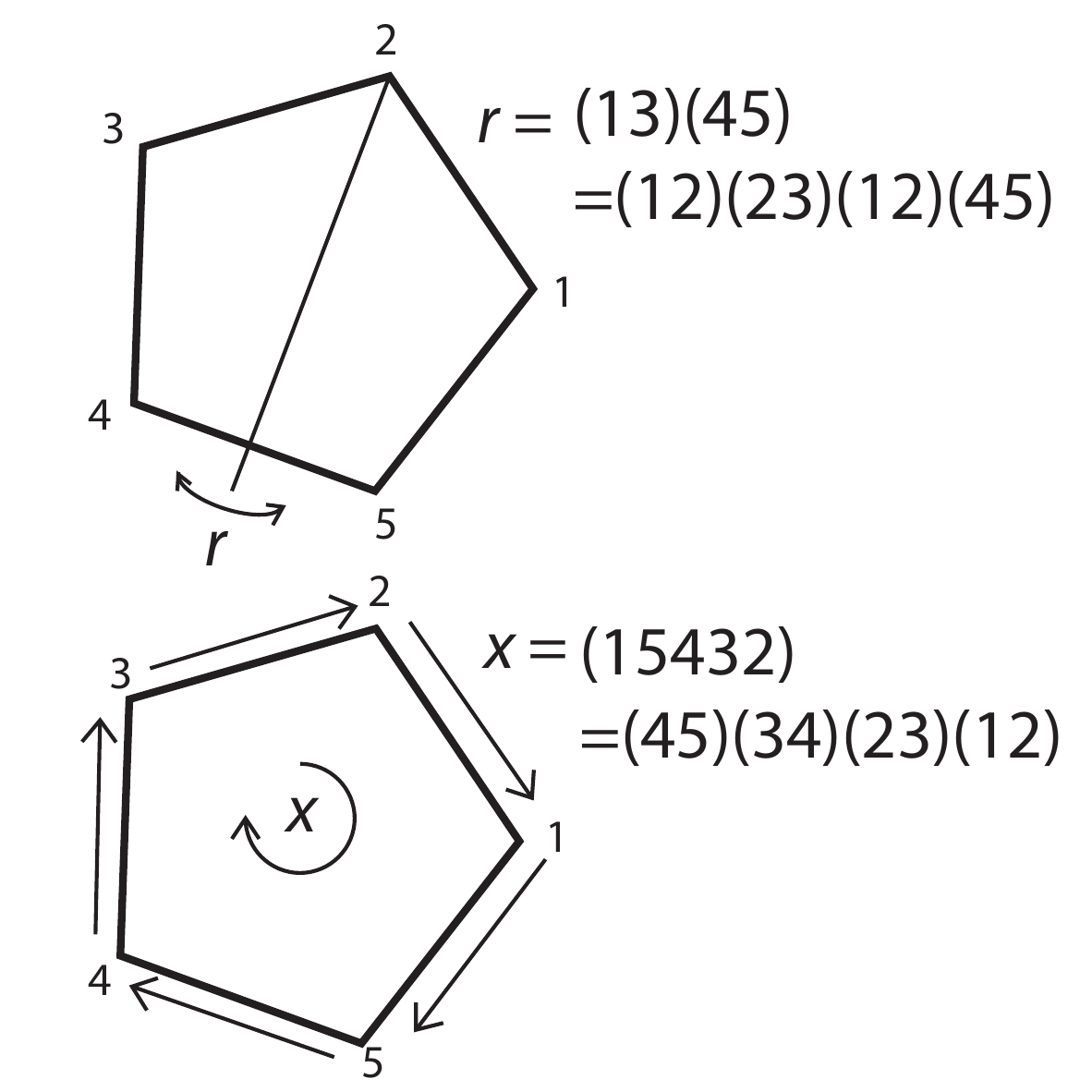}}
\]
\caption{A choice of generators for the dihedral group}
\label{pentup}
\end{figure}

The dihedral group $D_5$ is the collection of symmetries of a regular pentagon. The reflections are reflections in an axis that passes through a vertex. For reasons that are not immediately apparent, we will let $r$ denote the reflection through the axis that passes through the vertex labeled $2$   and $x$ will denote a clockwise rotation. See Fig.~\ref{pentup}. In the figure, these elements are represented as permutations in $\Sigma_5$ that have been rewritten in terms of the generators, $t_j$,  that are adjacent transpositions. For example, $t_1=(12)$, $t_2=(23)$, etc.

\begin{figure}
\[
\scalebox{.35}{\includegraphics{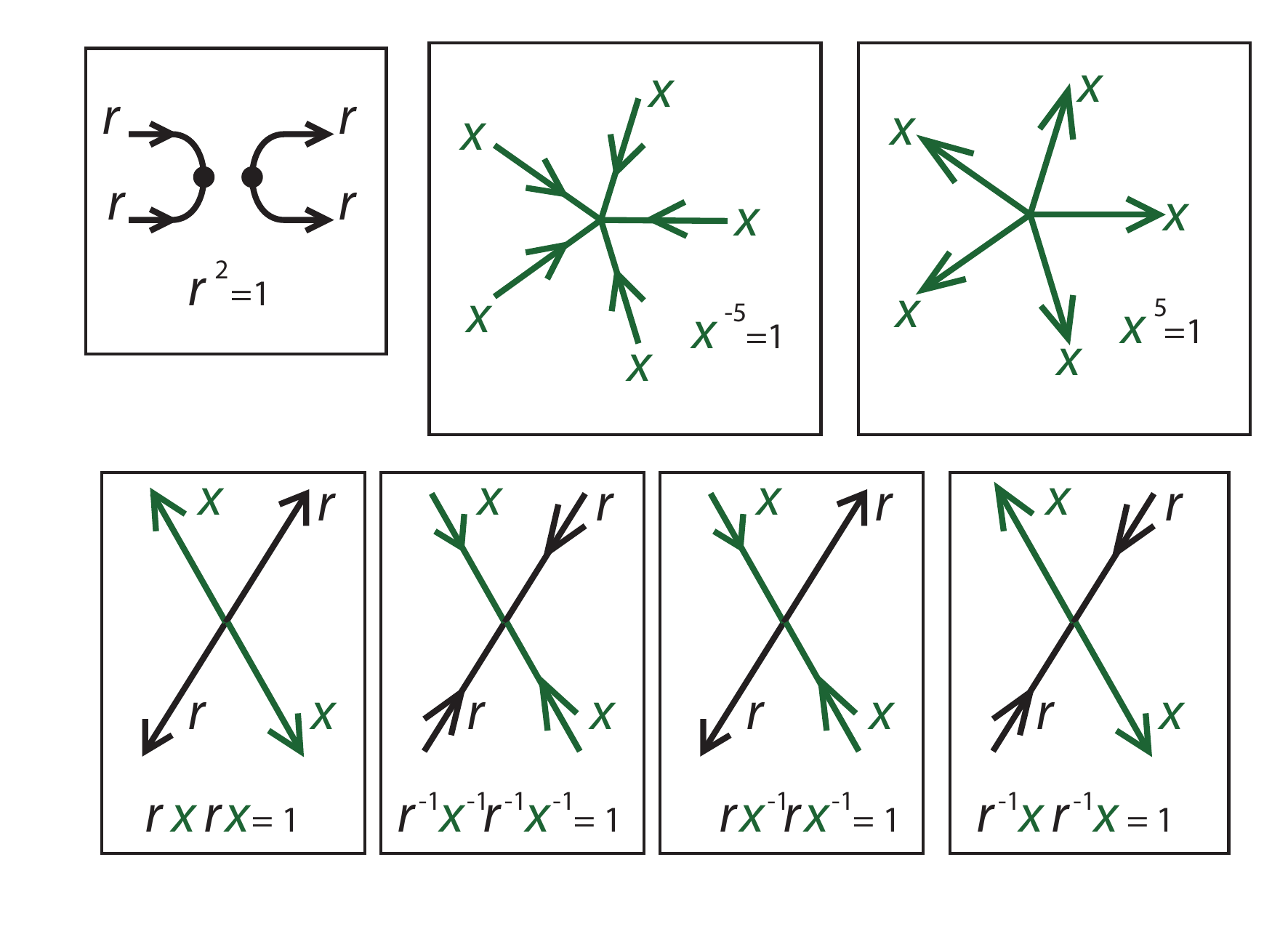}}
\]
\caption{Schematics of the relations in the dihedral group $D_5$}
\label{relations}
\end{figure}

Figure~\ref{relations} indicates the relations in the dihedral group.  A black arc is used to represent the  reflection $r$.  A green arc represents the rotation $x$. A {\it dihedral chart} is an analogue of a braid chart. It is an embedded labeled graph in a disk with four types of vertices. 
\begin{enumerate}
\item A {\it black vertex} is monovalent with an arc that is labeled either by $r$ or by $x$ that is incident to it. Usually, we will only use black vertices that are incident to arcs labeled by $r$.  However, cyclic branched covers (of degree $5$) can also be thought of via representations into the cyclic subgroup of the dihedral group. 
\item A bivalent vertex that  is incident to a pair of inward or outward pointing edges that  are labeled $r$. 
\item A vertex of valence five in which the five incident edges are all pointing inward or outward and each is labeled with an $x$. 
\item A valence four vertex has alternating edges cyclicly labeled $r^{\pm1}$, $x^{\pm1}$, $r^{\pm1}$, $x^{\pm1}$. 
\end{enumerate}
The last three types of vertices are exemplified in Fig.~\ref{relations}. Since $r=r^{-1}$, the valence four vertices  may represent the various relations
$rxrx=1$, $r^{-1}x^{-1}r^{-1}x^{-1}=1$,  $rx^{-1}rx^{-1}=1$, or $r^{-1}xr^{-1}x=1,$ as the figure indicates.


\begin{figure}
\[
\scalebox{.1}{\includegraphics{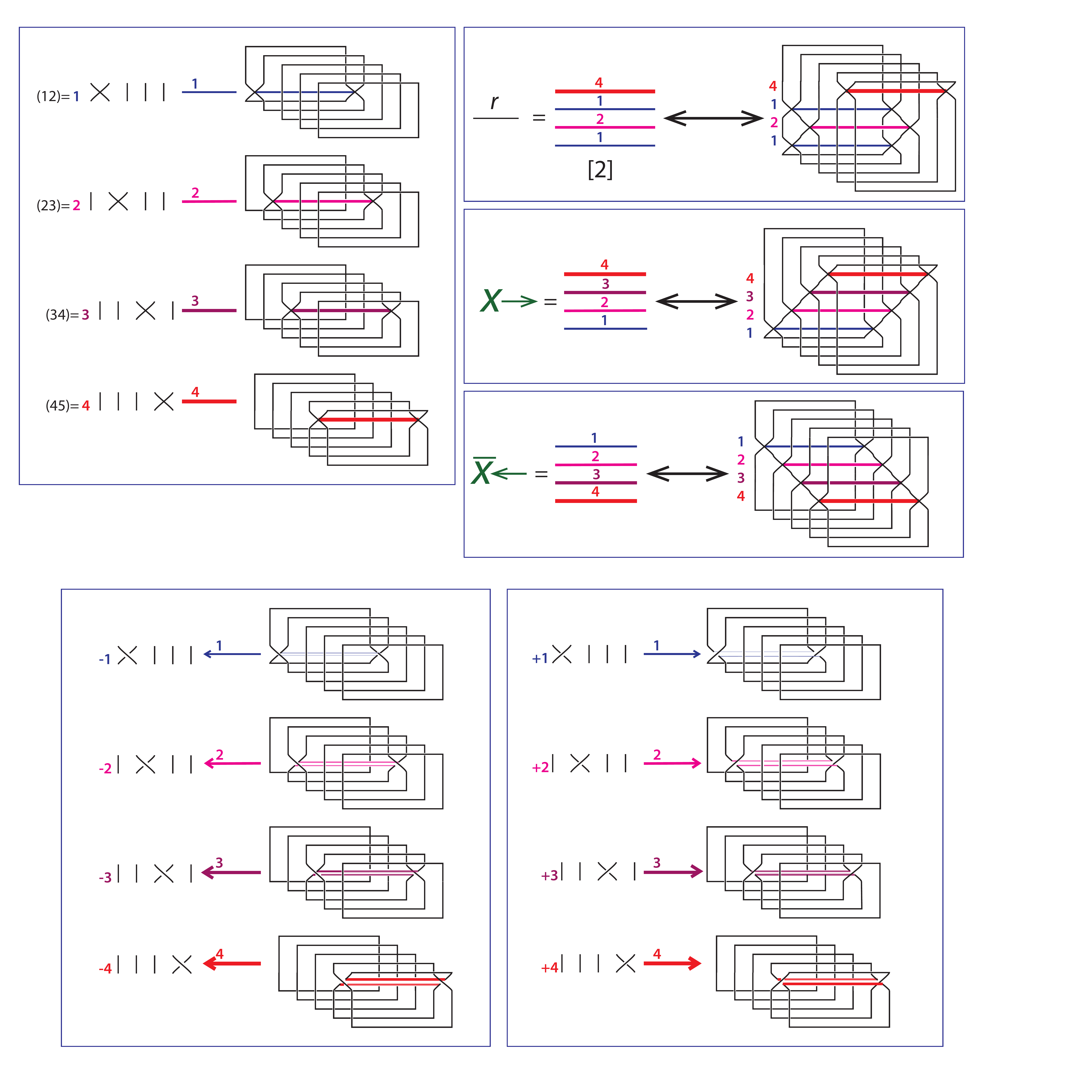}}
\]
\caption{A glossary of diagrams}
\label{legend}
\end{figure}

\begin{center} 
\rule{3in}{0.005in}
\end{center}

\noindent 
{\bf Remark.} The vertices that are not black, in  the braid, permutation,  and  dihedral cases correspond to $2$-cells which are found among the structures in the classifying space for the corresponding group. Suppose that a finite presentation of a group $G=\langle X: R \rangle$ is given. Then the classifying space for the group $BG$, which is a $K(G,1)$-space, can be constructed with a single vertex, an edge for each generator in $X$, a $2$-cell for each relation in $R$, and a host of higher dimensional cells that are attached to annihilate any higher homotopy groups. Consider the planar duals to non-black vertices in the braid, permutation, or dihedral charts. These are polygons that correspond to the $2$-cells in the corresponding $K(G,1)$. 

The classifying space $BG$ classifies irregular covers whose group of deck transformations is $G$. In the construction of branched covers, we remove the open tubular neighborhoods of the codimension-two  branch set, and start from the irregular cover of the exterior, $E(K) = S^m\setminus {\mbox{\rm int}}(N(K))$. Such a cover is determined by the homotopy class of a map $BG \longleftarrow E(K)$. In the case $m=2$, we are considering an irregular cover of a planar surface with boundary. And the vertices in the corresponding chart represent the pull-backs of the $2$-cells in  $BG$. In Section~\ref{ProofofDihedral}, these methods will be used to construct dihedral charts in general.

In the $m=3$ case, (branched covers of $S^3$ that are branched over a knot or link $K$), we may also need to pull back some $3$-cells in the classifying space. In \cite{CKFold} or \cite{KamBook}, the braid chart moves are listed. See also Fig.~\ref{almostKmoves}. Some of these correspond to $3$-cells that are in the classifying space of the braid group. One such move (the quadruple point move) occurs in the construction of the $D_5$ cover of $S^3$ that is branched along the torus knot $T(2,5)$. 
The reader is reminded that permutation charts are obtained from braid charts by removing the orientations along the edges. 

In fact, almost all of our examples will be constructed via permutation charts. However, the dihedral charts that were introduced in this section are useful as a schematic for the corresponding underlying permutation charts. And, as we will see, two key lemmas will be expressed by means of dihedral charts, but their  proofs are implemented at the level of permutations. These lemmas can also be interpreted as certain homotopies of maps between the exterior and the classifying space.

\begin{figure}[htb]
\[
\scalebox{.15}{\includegraphics{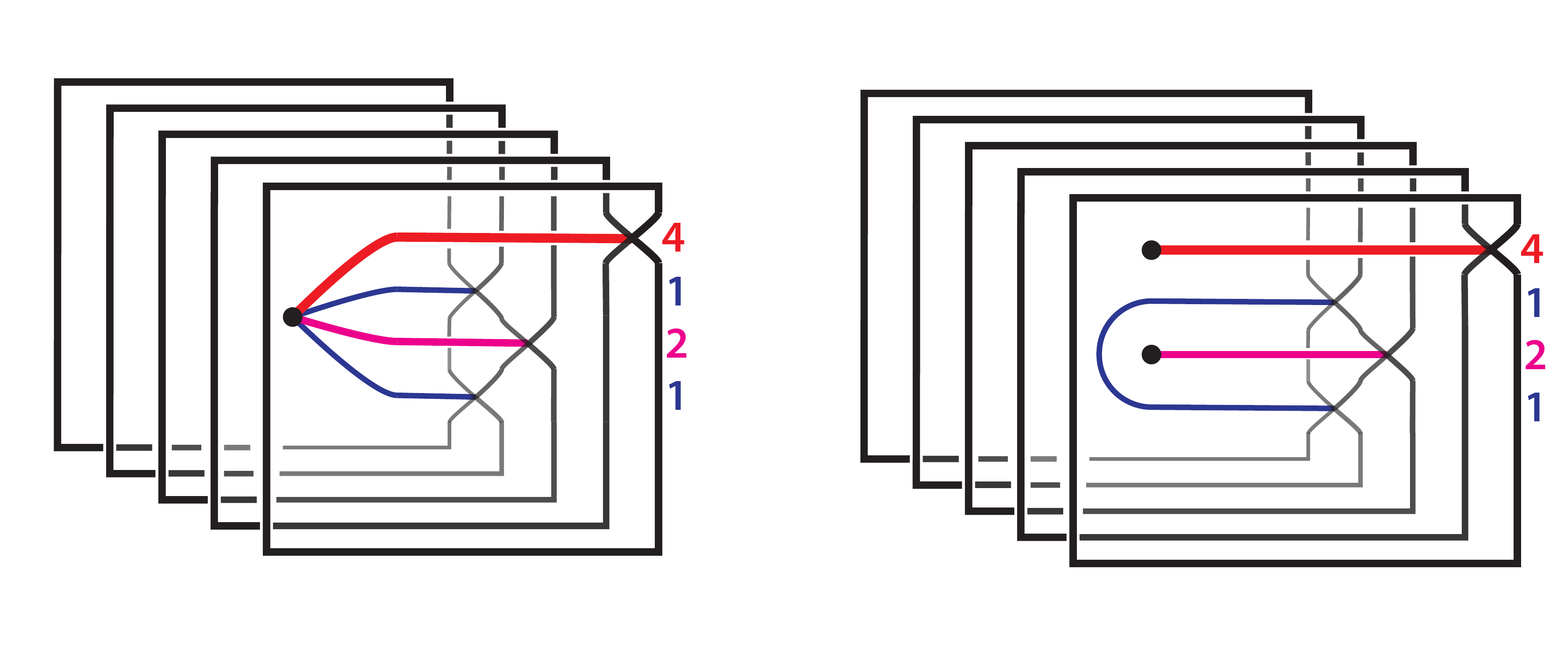}}
\]
\caption{The principal dihedral branch point and its resolution into simple branch points}
\label{resolve}
\end{figure}

\subsection{The reflections in $D_5$ as conjugates}

\begin{sloppypar}
Consider the reflections in $D_5$. Label the vertices of the pentagon $1$ through $5$ as in Fig.~\ref{pentup}. As a permutation, the reflection $r$ --- that fixes the vertex $2$ --- is given as $r=(13)(45)$. We write $[2]=(13)(45)$ to denote this reflection. We compute further: $rx=(13)(45)(15432)=(14)(23)=[5]$ --- the reflection that fixes the vertex $5$; $rx^2=(14)(23)(15432)=(15)(24)=[3]$ --- the reflection that fixes the vertex $3$; $rx^3= (15)(24)(15432)=(25)(34)=[1]$ --- the reflection that fixes the vertex $1$; $rx^4=(25)(34)(15432)= (12)(35)=[4]$ --- the reflection that fixes the vertex $4$. 
\end{sloppypar}

We are interested in representations of knot groups. The Wirtinger relations at a crossing express one meridional element as a conjugate of another. So when representing such groups into the dihedral group, it is more convenient to express all of the reflections as conjugates of $[2]=(13)(45)$. Here is another routine calculation: $x^{-1}rx=rx^2 =[3]$; $x^{-2}rx^2=rx^4=[4];$ $x^{-3}rx^3=rx^6=rx=[5]$; and $x^{-4}r x^4=rx^3=[1].$

\subsection{Dihedral covers of $S^2$ branched over two points}

In this section,  more diagrammatic notation is established. The generators $t_1$ through $t_4$ in the permutation group will be schematized, and each will be assigned a color that is  consistent throughout the illustrations. 

\begin{figure}
\[
\scalebox{.2}{\includegraphics{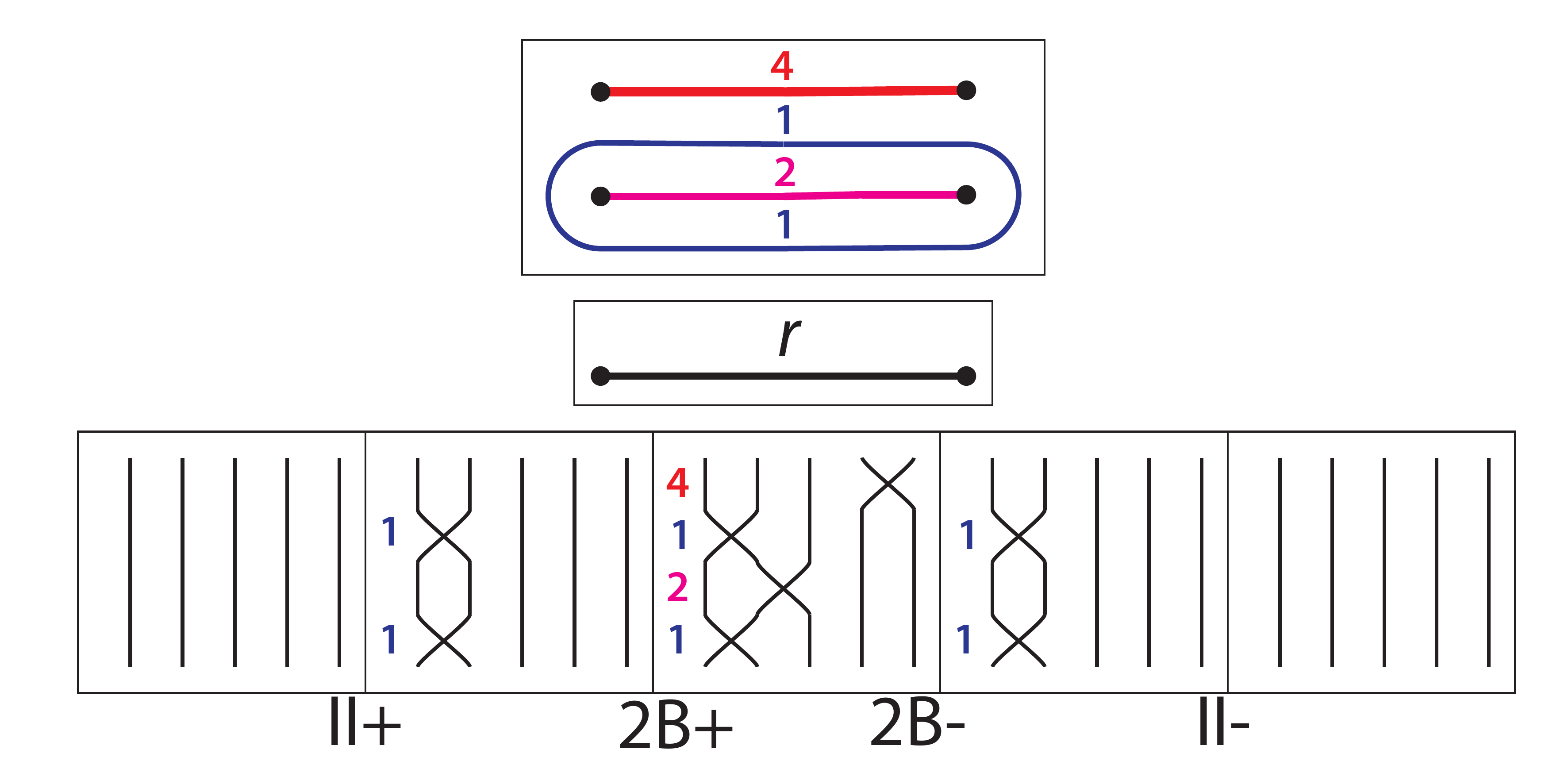}}
\]
\caption{A dihedral cover of $S^2$}
\label{submarine}
\end{figure}

Then we construct foldings of dihedral branched covers of $S^2$ that are branched over $2$ points. Because the reflections are expressed as conjugates of the reflection $r=[2]$, each of the foldings in Fig.~\ref{submarine} and Fig.~\ref{fromnone2five} can be lifted to an embedded braiding in $S^2 \times [0,6] \times [-1,1]$. That is, the permutation charts can be oriented to become braid charts.

Figure~\ref{legend} contains a legend of the graphical notation. The two boxes in the bottom of the figure indicate the two possible local liftings of a permutation folding to an embedded braiding. That is, a permutation folding lands in $S^3 \times [0,6]$ and an embedded braiding lands in $S^3\times [0, 6] \times [-1,1]$. So the permutation generator $t_j$ lifts to either $\sigma_j$ or  its inverse $\overline{\sigma_j}$. In order to avoid subscripts, the permutation generators in $\Sigma_5$ are written as $1$ through $4$, and their associated colors are indicated. The corresponding braid generators are written as $\pm1, \ldots, \pm 4$ where the $(+)$ and $(-)$ signs are  included. The three boxes on the top right indicate the permutation representatives that we have chosen for $r$, $x$, and $x^{-1}$. Note that the orientations of the green arrows indicate where the representative cycle is $(15432)$ for $x$, or $(12345)$ for $x^{-1} = \overline{x}$.

On the left of Fig.~\ref{resolve}, a branch point appears to the left the reflection $r=[2]=(13)(45)$. 
The branch point itself is over-emphasized as a large black vertex. It occurs as a {\it deus ex machina} phenomenon that is intrinsically difficult to depict. On the other hand, the reflection $(13)=(12)(23)(12)$ can be resolved to the simple branch point between the second and third sheet which is conjugated by the transposition $(12)$. Meanwhile, the reflection $(45)$ is independent of this. On the right of the figure, the singularity has been resolved into a pair of simple branch points. We resolved the branch point into simple branch points merely for our own mental conveniences.

\begin{figure}[htb]
\[
\scalebox{.23}{\includegraphics{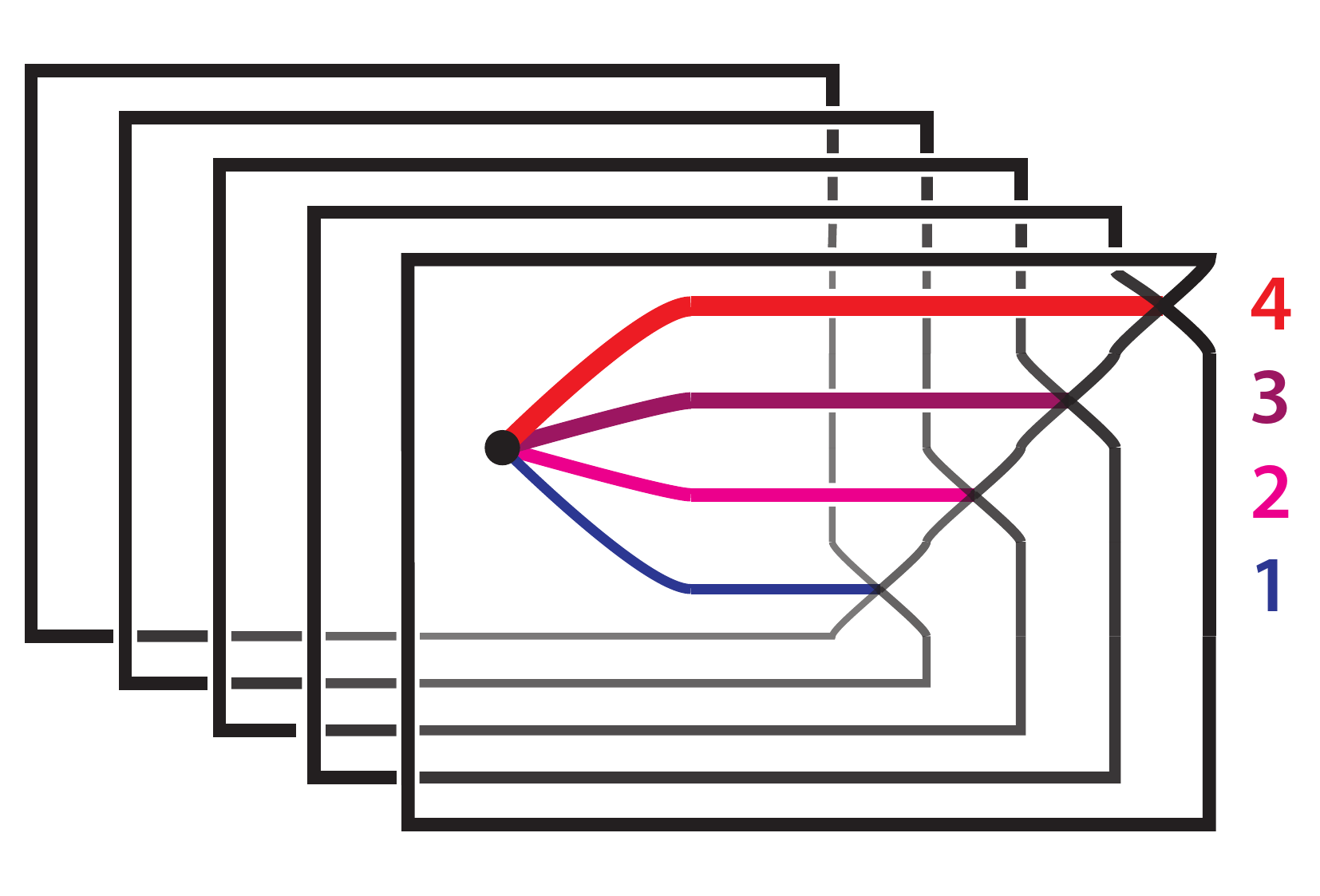}}
\]
\caption{A degree $5$ cyclic branch point modeled upon $z\mapsto z^5$}
\label{ex5}
\end{figure}

\begin{figure}[htb]
\[
\scalebox{.23}{\includegraphics{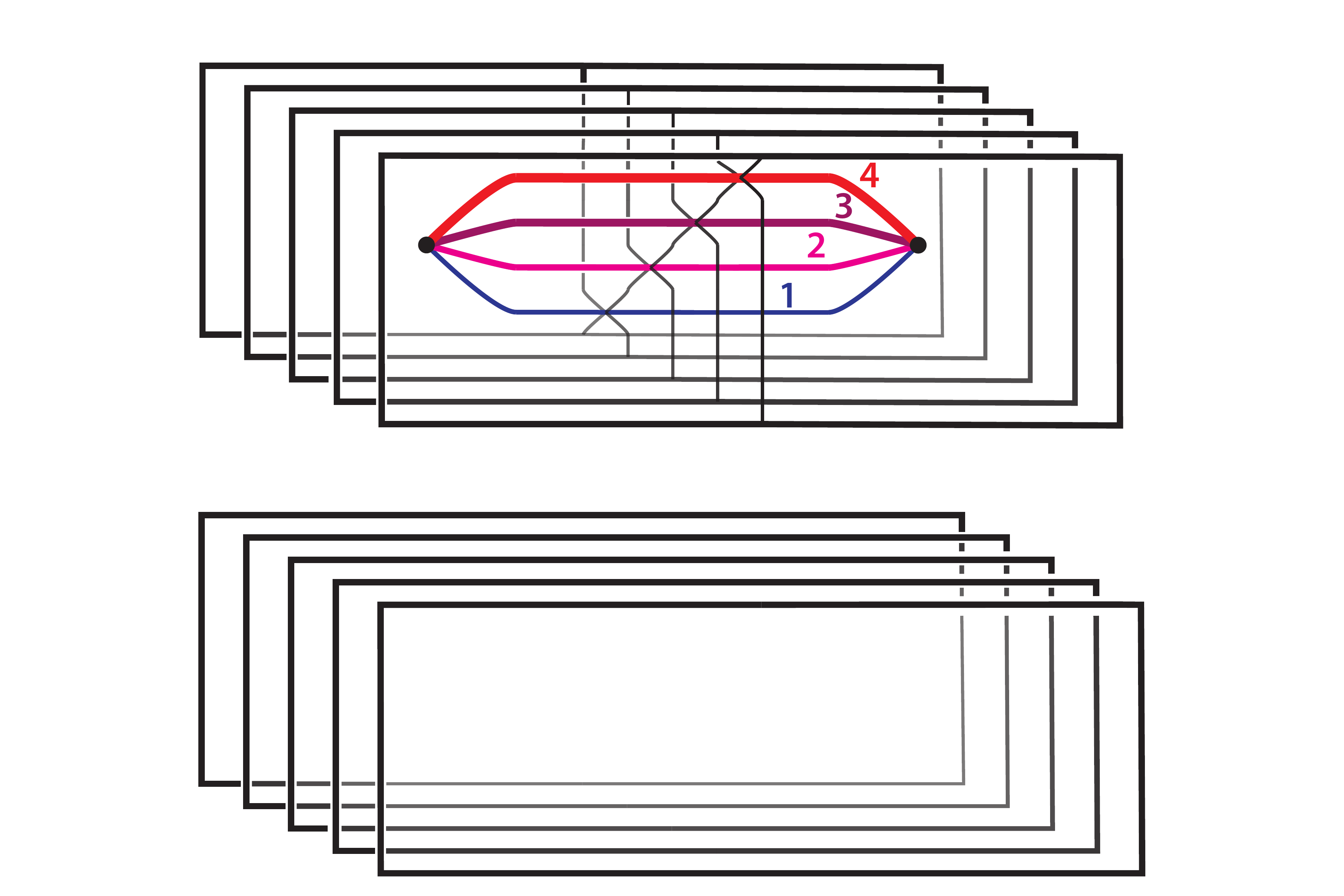}}
\]
\caption{A five-fold cyclic branched cover and its resolution}
\label{fromnone2five}
\end{figure}

\begin{figure}[htb]
\[
\scalebox{.23}{\includegraphics{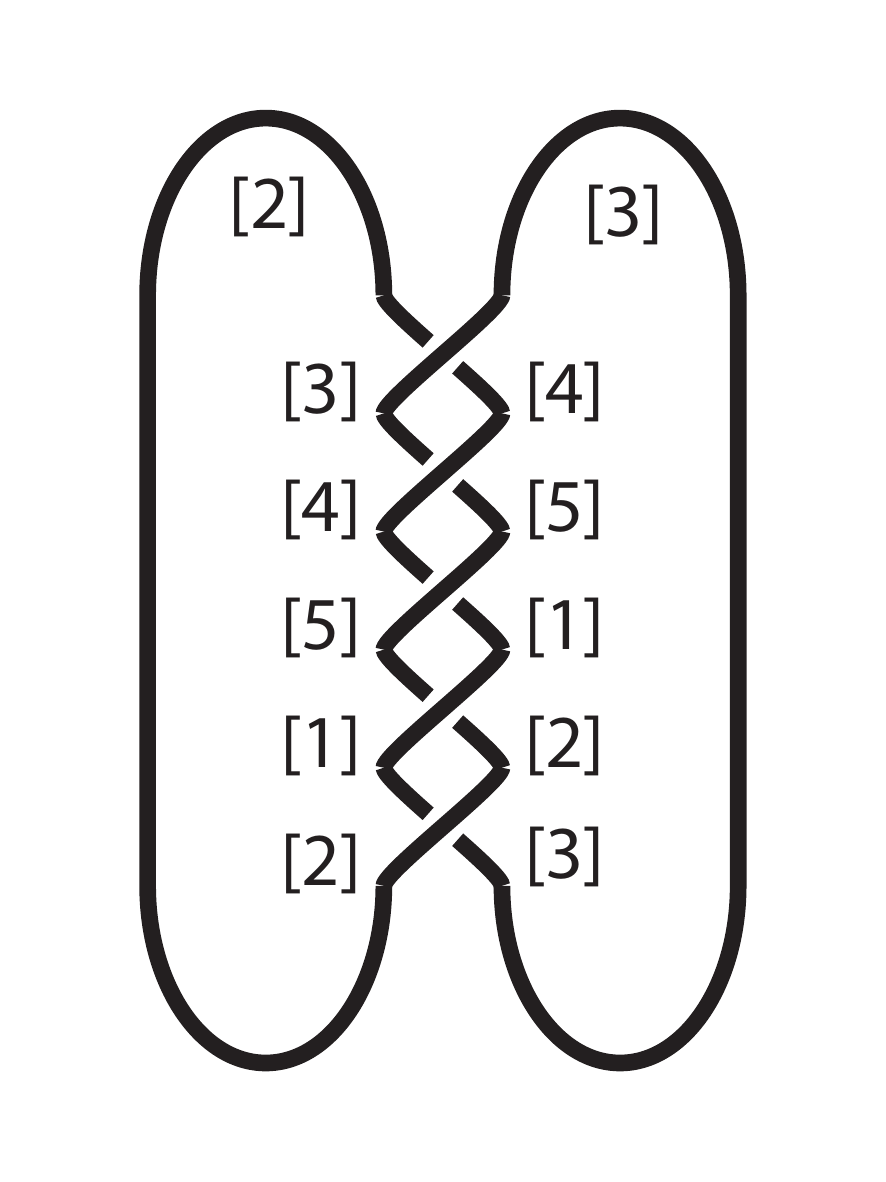}}
\]
\caption{The torus knot $T(2,5)$}
\label{teatwofive}
\end{figure}

Meanwhile, Fig.~\ref{submarine} illustrates a permutation chart 
that has four black vertices. The chart is a subset of the disk $D^2$ which in turn is considered as a subset of the sphere $S^2$. So the chart is indicating a dihedral branched cover of $S^2$ with four singularities. The middle of the illustration indicates the corresponding dihedral chart. From the permutation chart that has these four simple branch points, a braiding can be explicitly envisioned. For example the permutation movie 
\[ {\mbox{\rm id}} \stackrel{{II+}}{\Longrightarrow} (12)(12) 
\stackrel{{2B+}}{\Longrightarrow} (12)(23)(12)(45) \stackrel{{2B-}}{\Longrightarrow} (12)(12) \stackrel{{II-}}{\Longrightarrow}{\mbox{\rm id}} \]
is illustrated at the bottom of the figure. In the singular (merged) situation, the second and fourth stills of the movie are absent. 

The notation  indicates: (A) a type II move that increases the number of crossings ($II+$), (B) the addition of a pair of branch points ($2B+$), (C) 
 the removal of a pair of branch points ($2B-$), and (D) a type II move that decreases the number of crossings ($II-$). 

We can compute directly that the branched covering is the union of three $2$-spheres. Or if you prefer, you may use the Riemann-Hurwitz formula.

The rotation $x\in D_5$ represents an order $5$  generator of  a cyclic subgroup. We form the cone upon the string diagram that represents the cycle $(15432)=(45)(34)(23)(12)$ as illustrated in Fig.~\ref{ex5}. This is a model of the cyclic branch point the corresponds to $z\mapsto z^5$. While we won't use this type of branch point in the construction of dihedral covers, we will use it for the construction of cyclic covers. One can imagine the graph of $g(z)=z^5$, when restricted to $D^2 = \{ r e^{{\boldsymbol i} \theta}: 0\le r \le 1 \ \ \& \ \ 0\le \theta < 2 \pi \}$ as a subset of $D^2 \times D^2$, and project it onto a $3$-dimensional subspace. The illustration in Fig.~\ref{ex5} is a topological depiction of such a projection.

In Fig.~\ref{fromnone2five}, two illustrations are given. The top drawing is a topologist's rendering of the the fifth power function $z\mapsto z^5$ as it would appear  when restricted to a rectangular disk $D$ that contains $0$ and $\infty$ and that is a subset of the Riemann sphere $\C \cup \{ \infty \}$. The top drawing indicates a $5$-fold cyclic branched cover of $S^2$, $M$, that has two branch points and that is folded. The manifold $M$ also is a $2$-sphere.

The bottom drawing means to illustrate the trivial covering of the sphere with five unbranched sheets. Each is a $2$-dimensional sphere. The  covering $M$ can be obtained from  $\cup_{i=1}^5 S^2_i$  by attaching four $1$-handles --- thereby connecting the first sphere to the second, the second to the third, etc. In the opposite direction, four $2$-handles can be used to remove the branch set. In Section~\ref{cyclicproof}, these differences occur at optimal points of the knot diagram. 

Similarly, one can transition from an empty dihedral chart to one that has a simple arc, labeled $r$, that starts and ends at two black vertices. The empty chart represents five embedded spheres. The chart containing the arc labeled $r$ represents a branched cover of $S^2$ as in Fig.~\ref{submarine}. The difference between the charts is the attachment of $1$-, or $2$-handles as the direction of the difference indicates.

\section{A schematic description of a dihedral cover}
\label{thesteps}

Consider the torus knot $T(2,5)$ that is the braid closure of the braid $\sigma_1^5$ as depicted in Fig.~\ref{teatwofive}. A dihedral coloring is also presented. Thus we have a group homomorphism
$D_5 \stackrel{\phi}{\longleftarrow} \pi_1 (E(T(2,5)))$ in which reflections are assigned to each of the meridional generators. This coloring induces a branched covering $S^3 \stackrel{B}{\longleftarrow} M_\phi$. See, e,g,~\cite{HildenPac78}.

 A schematic  set of pictures  is illustrated and labeled Step 1 through Step 21. 
 Each is a dihedral chart as described above. These   will  be used to construct a folding of  the branched cover $M_\phi$.
 
 Specifically, 
the torus knot $T(2,5)$ is the $2$-bridge knot  with the bridges depicted at the top of the diagram in Fig~\ref{teatwofive}. The left-hand bridge is colored by $[2]=r=(13)(45)$, and the right-hand bridge is colored by $[3]=x^{-1}rx$.  
The  illustrations  Step 1 through Step 21 
indicate successive cross-sectional levels of the knot.  Each is  a  $2$-sphere (or rather a $2$-disk contained therein) together with schema to construct the dihedral branched cover of that sphere. The $2$-spheres essentially foliate $S^3$. Empty diagrams which indicate the $5$-fold trivial cover of $S^2$ that would flank the illustrations that we present are not included.  

In a somewhat ironical fashion, such empty diagrams convey some crucial information. The $3$-sphere can be thought of as a tropical $S^2\times I$ together with a north and south pole.  The knot lies in the tropics, and the polar caps are $3$-dimensional balls. The $5$-fold trivial cover of each of these $B^3$\/s is the union of five $3$-dimensional balls. The  two empty diagrams that flank the illustrations represent five concentric $S^2$\/s. These successively bound $B^3$\/s from innermost to outermost. Those ten (five in the north and five in the  south) $B^3$\/s form the cover of $S^3$ outside the region that contains the knot. We will return to this idea in Section~\ref{ProofofDihedral}.

The schematic illustrations  below indicate branched coverings of the $2$-sphere. They are latitudinal slices ($S^2\times \{t\}$) of the tropical $S^2\times I$ .
After Step $21$, there are no longer any branch points present. So the illustration iindicates the union of five  intersecting $2$-spheres in $3$-space. The schematic will be translated into a less symbolic representation of the intersection sets among the spheres. Selected stages of the isotopy of these  will be depicted following this section on schematics.

The difference between the  excluded $0$\/th (empty) diagram and Step 1 is the addition of the arc\ that is labeled by the reflection $r=[2]$ in Step 1. This difference is quantified as the addition of a pair of $1$-handles  between the empty diagram and the first step. Similarly, the $1$\/st and $2$\/nd illustrations differ by means of the addition of an arc, that is labeled by $r$, and that is surrounded by a loop  labeled by the rotation $x$. This loop is a conjugation of $r$, and another pair of handles is attached among the spheres.

Orientations  on the arcs of the illustrations, that are consistent from one to the next, are induced by an orientation of the knot $T(2,5)$ that is indicated in Fig.~\ref{teatwofive}.  

Each illustration is a dihedral chart that, in turn, represents a branched cover of the $2$-sphere.  Those illustrations that have four black dihedral vertices  are related by isotopies of dihedral charts.  At the top of the illustrations, a horizontal line segment that crosses the knot indicates the  rough level in $S^3$ at which the dihedral chart sits. The level is only approximate since, for example, between  
Step 2 and Step 3 that horizontal  line doesn't move. 

Between
 Step 3 and Step 4
the branch points in the center of the figures have moved. This difference codifies the topmost crossing of the knot. The lower center right branch point has moved to the left in front of the upper center left branch point which has moved to the right. Their relative positions are marked by means of the thin vertical lines in the illustrations that remain fixed throughout the sequence of illustrations.

The dihedral charts between Step 2 and Step 19 all have exactly four black dihedral vertices. These charts represent isotopic dihedral branch coverings of $S^2$ with four branch points (sources or sinks with an incident edge labeled $r$). Often the difference between successive steps is merely an isotopy. In particular, the inner black vertices swirl around each other up to Step 11. 

Subsequently, pairs of canceling valence four vertices that represent $rxrx=1$ or some conjugate of this relation are created --- for example between Step 12 and Step 13. When a pair of edges labeled by $r$ run anti-parallel, a $\cup/\cap$ relation can be added. After the illustrations are presented, the text will explain why these changes can occur. Meanwhile, textual commentary will be presented between the steps.

 \newpage
 
 We ask the reader to imagine an empty chart. Between the empty chart and Step 1, a left pointing edge that is labeled by the reflection $r$ is indicated.

\[
\scalebox{.15}{\includegraphics{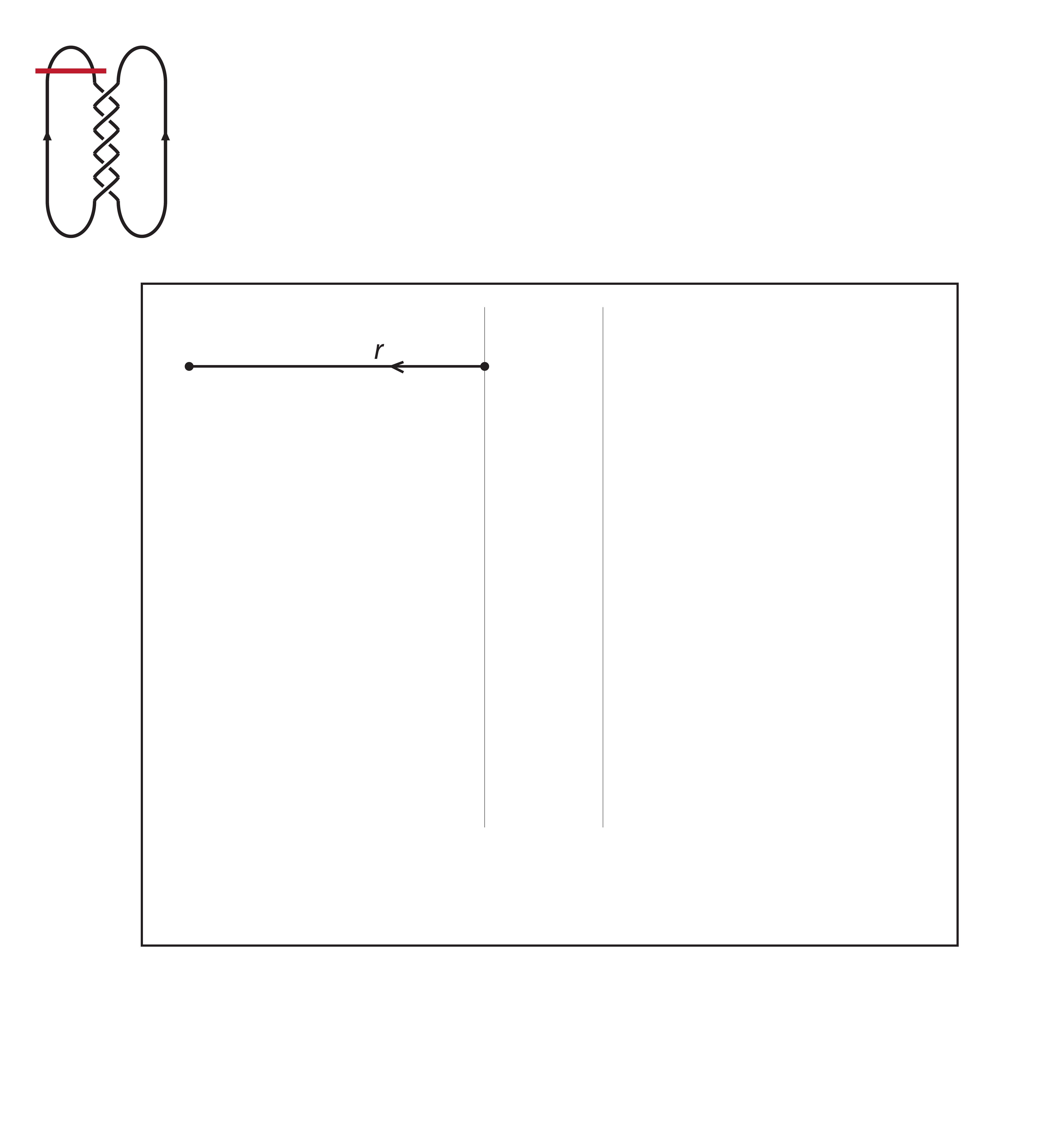}}
\]
\begin{center} {\bf Step 1} \end{center}

\newpage 

Between Steps 1 and 2, a right pointing edge labeled $r$ that is surrounded by a  clockwise loop labeled $x$ has been born. The two maximal bridges in the knot diagram of $T(2,5)$ have been passed. 
\[
\scalebox{.15}{\includegraphics{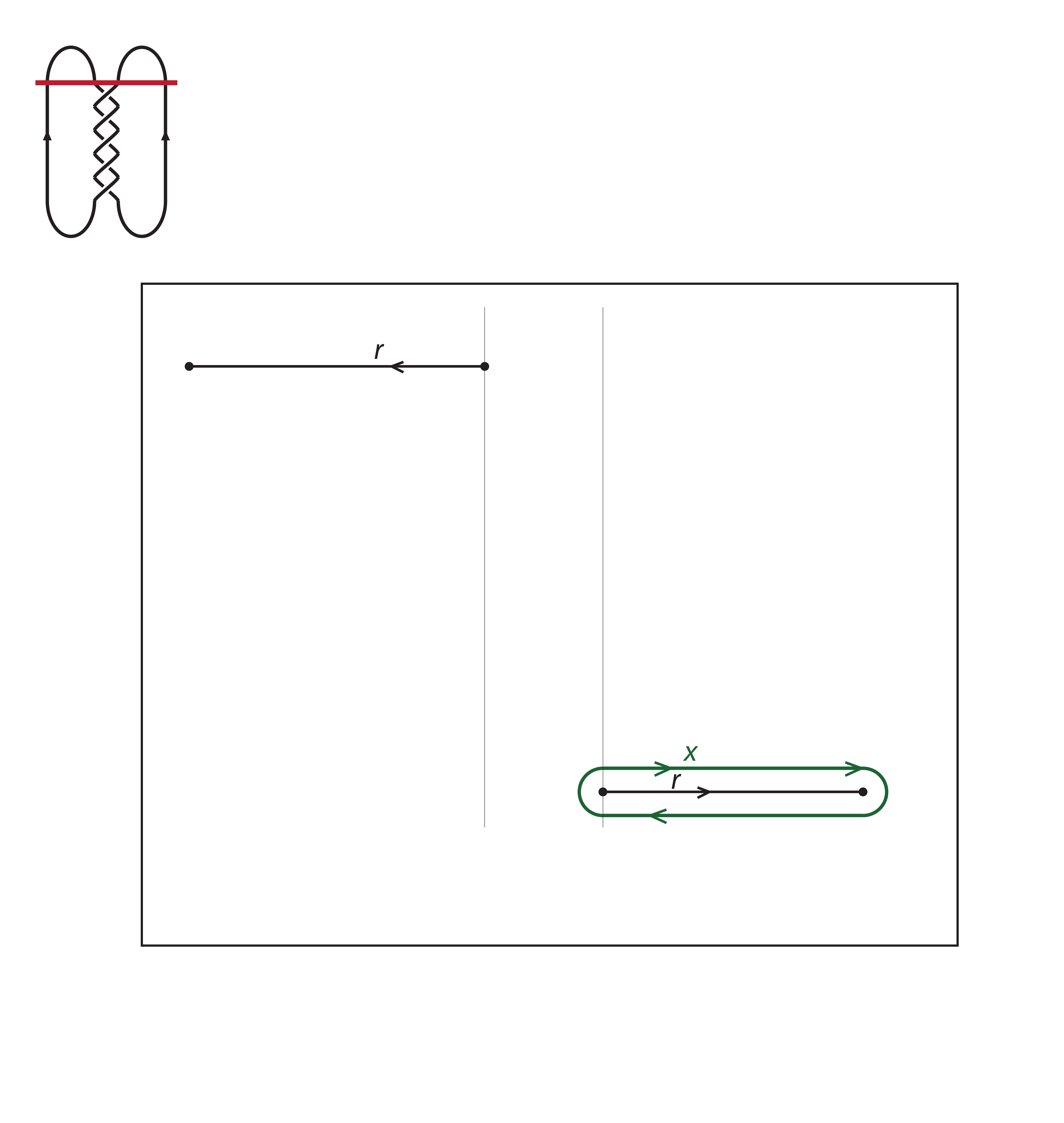}}
\]
\begin{center} {\bf Step 2} \end{center}

\newpage

Between Step 2 and Step 3, the inner black vertices have exchanged their positions in relation to the horizontal axis. The horizontal red arc on the knot diagram is  below the first crossing.

\[
\scalebox{.15}{\includegraphics{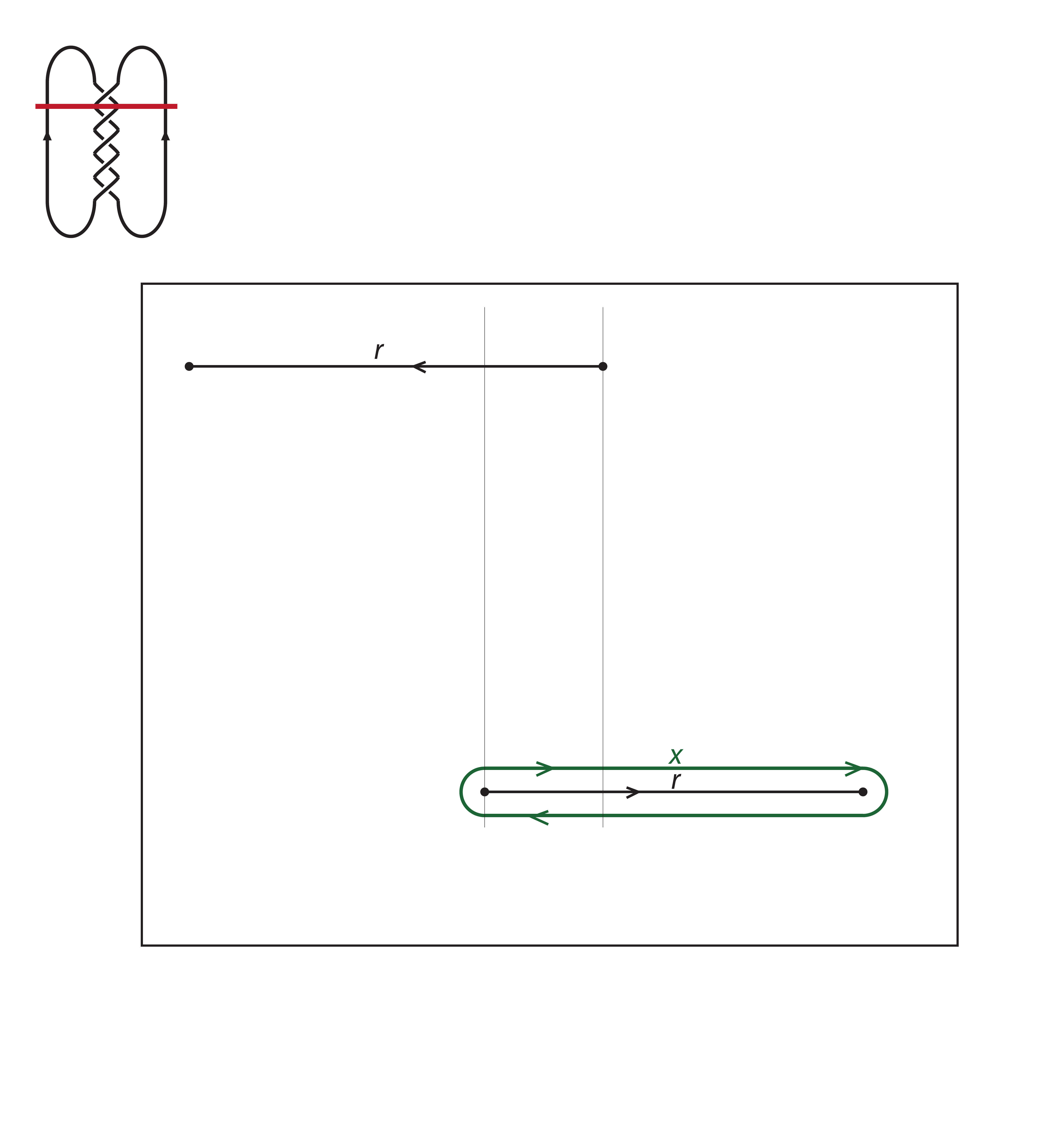}}
\]
\begin{center} {\bf Step 3} \end{center}

\newpage

 Between Step 3 and Step 4, semi-circular arcs on either side of the diagram are used to exchange the vertical positions of the inner black vertices. In this way, the vertex on the right will be able to pass in front of that on the left.
\[
\scalebox{.15}{\includegraphics{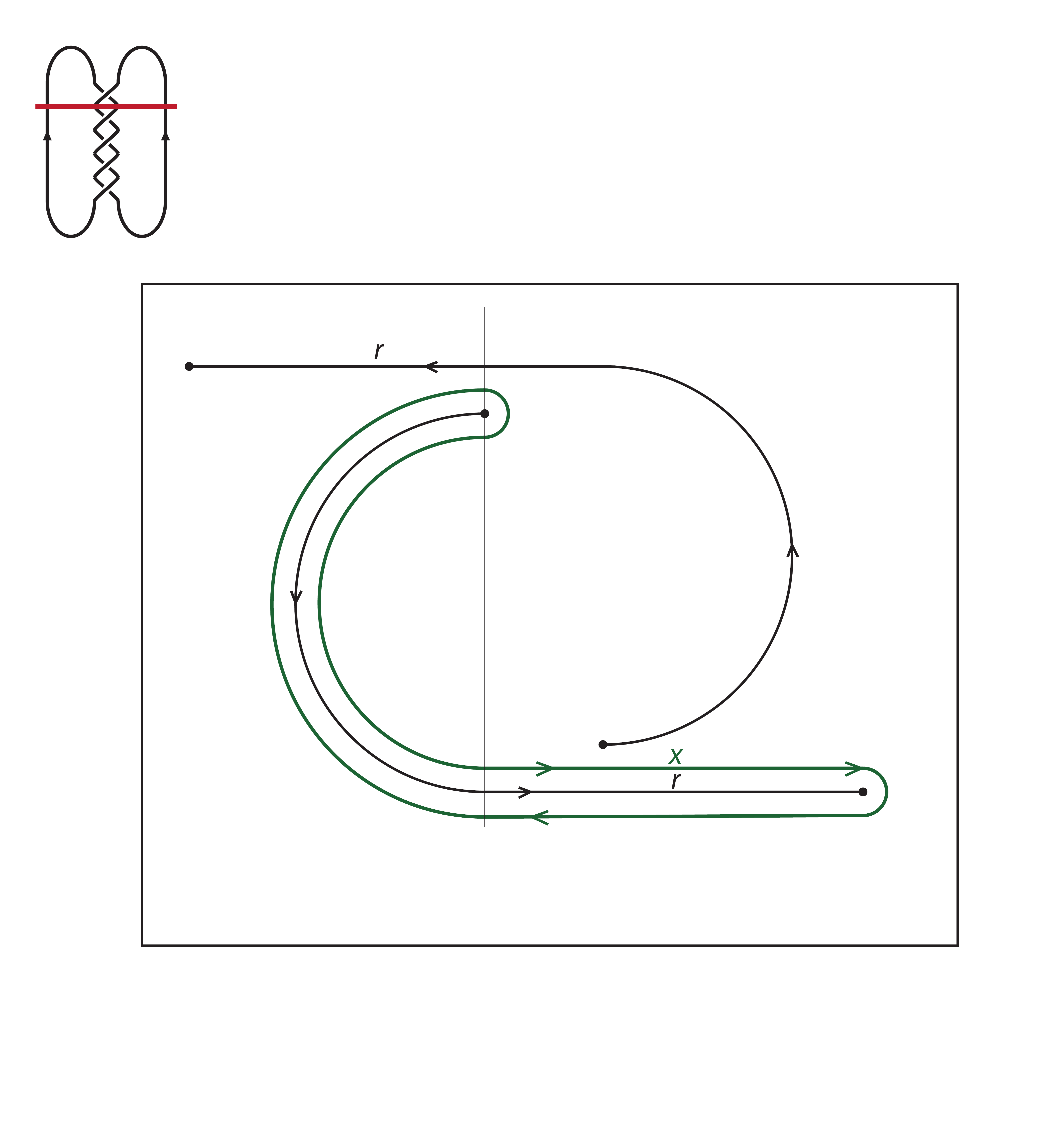}}
\]
\begin{center} {\bf Step 4}  \end{center}

\newpage

Between Step 4 and Step 5, the inner black vertices have exchanged their positions in relation to the horizontal axis. The horizontal red arc on the knot diagram is  below the second crossing.

\[
\scalebox{.15}{\includegraphics{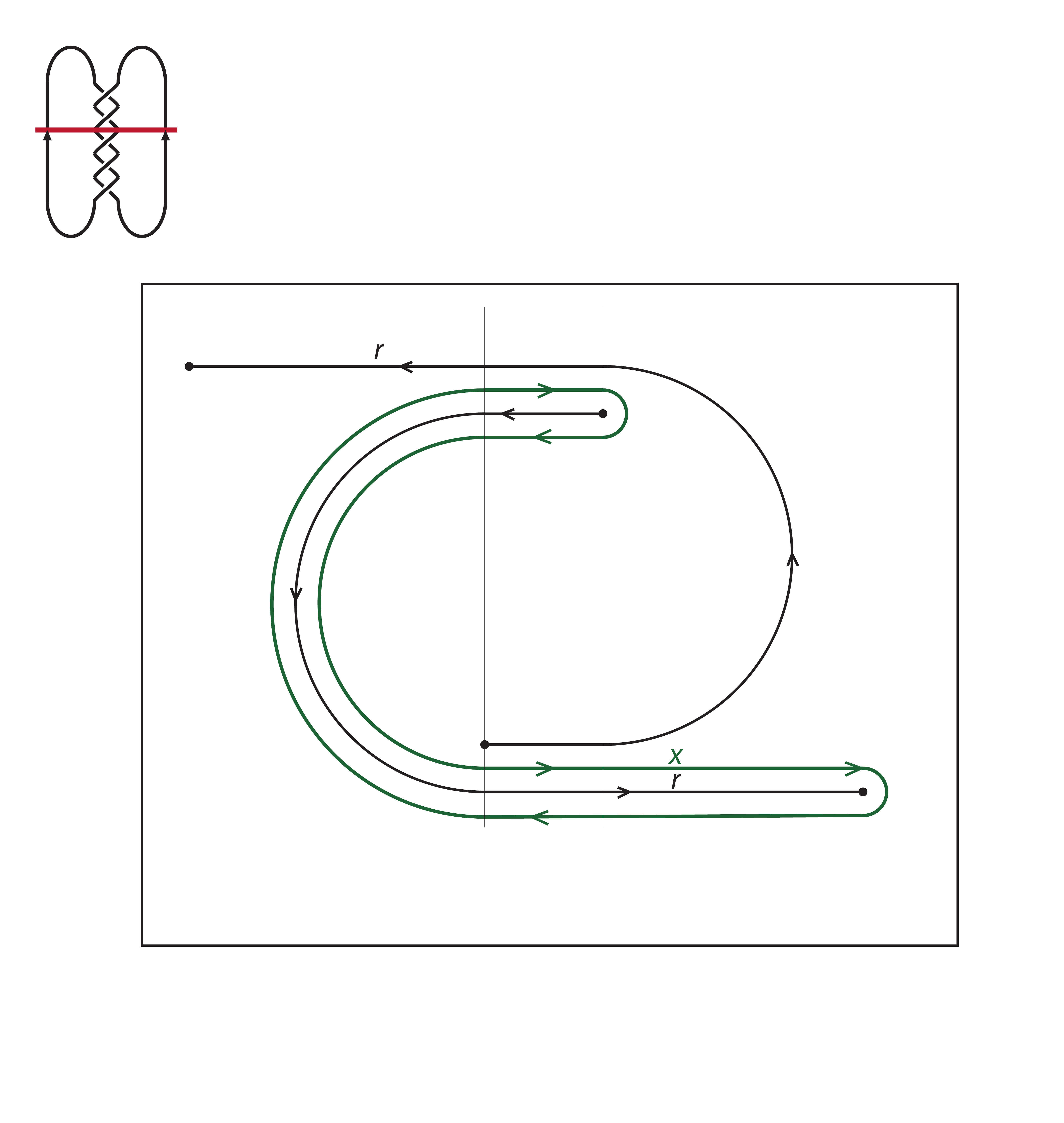}}
\]
\begin{center} {\bf Step 5} \end{center}

\newpage

 Between Step 5 and Step 6, semi-circular arcs on either side of the diagram are used to exchange the vertical positions of the inner black vertices. In this way, the vertex on the right will be able to pass in front of that on the left.

\[
\scalebox{.15}{\includegraphics{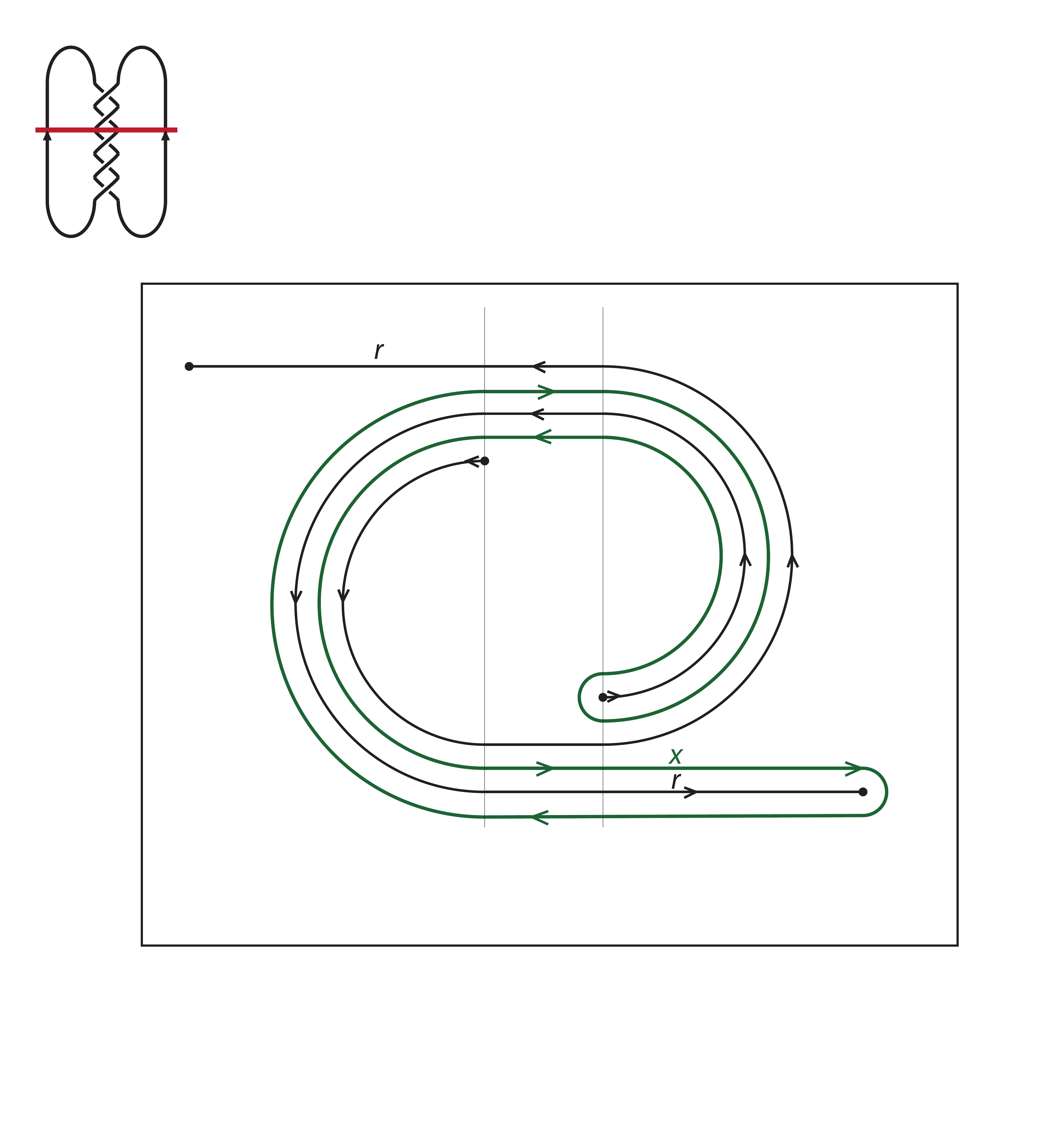}}
\]
\begin{center} {\bf Step 6} \end{center}

\newpage

Between Step 6 and Step 7, the inner black vertices have exchanged their positions in relation to the horizontal axis. The horizontal red arc on the knot diagram is  below the third crossing.
\[
\scalebox{.15}{\includegraphics{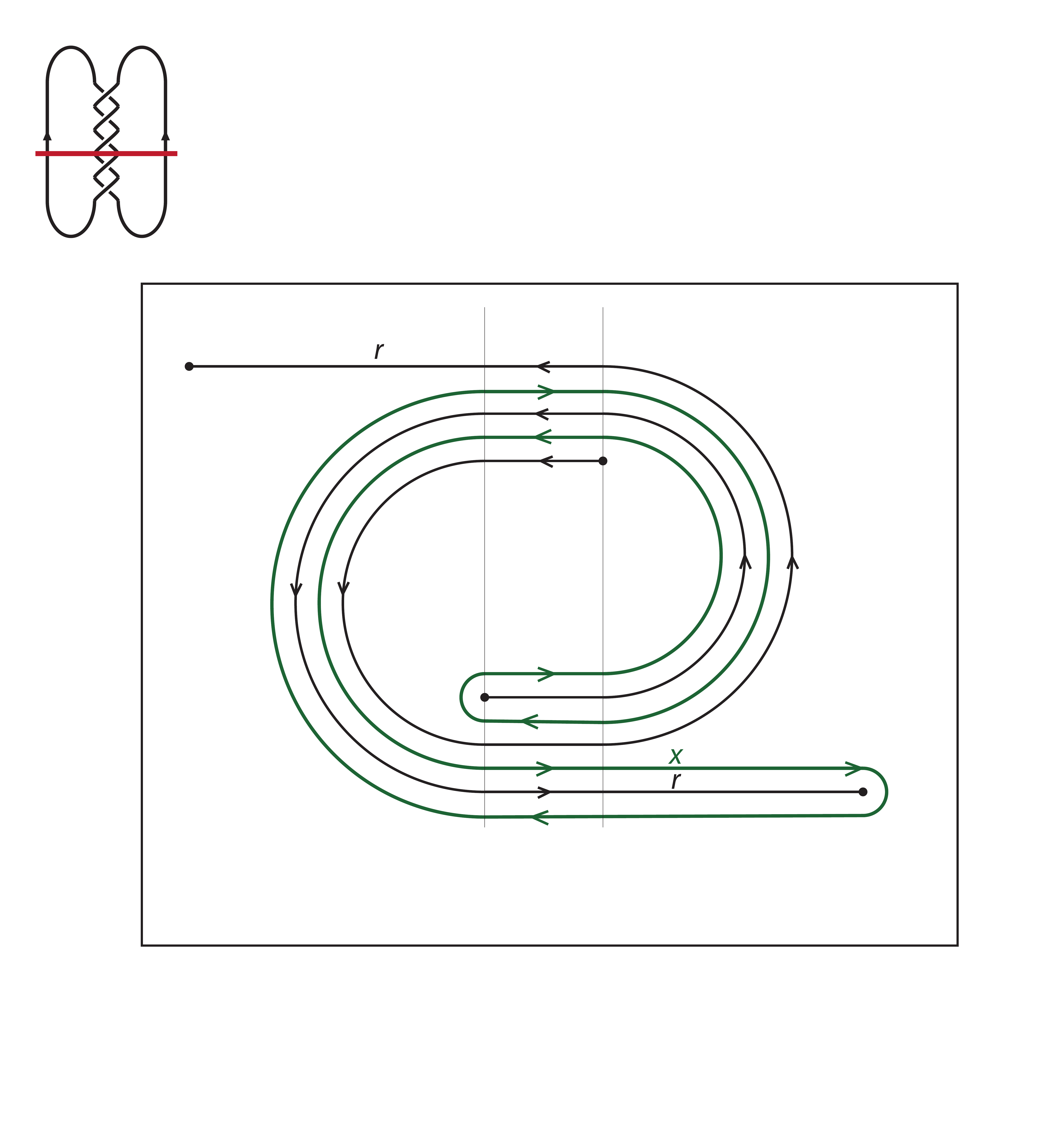}}
\]
\begin{center} {\bf Step 7} \end{center}

\newpage

 Between Step 7 and Step 8, semi-circular arcs on either side of the diagram are used to exchange the vertical positions of the inner black vertices. In this way, the vertex on the right will be able to pass in front of that on the left.

\[
\scalebox{.15}{\includegraphics{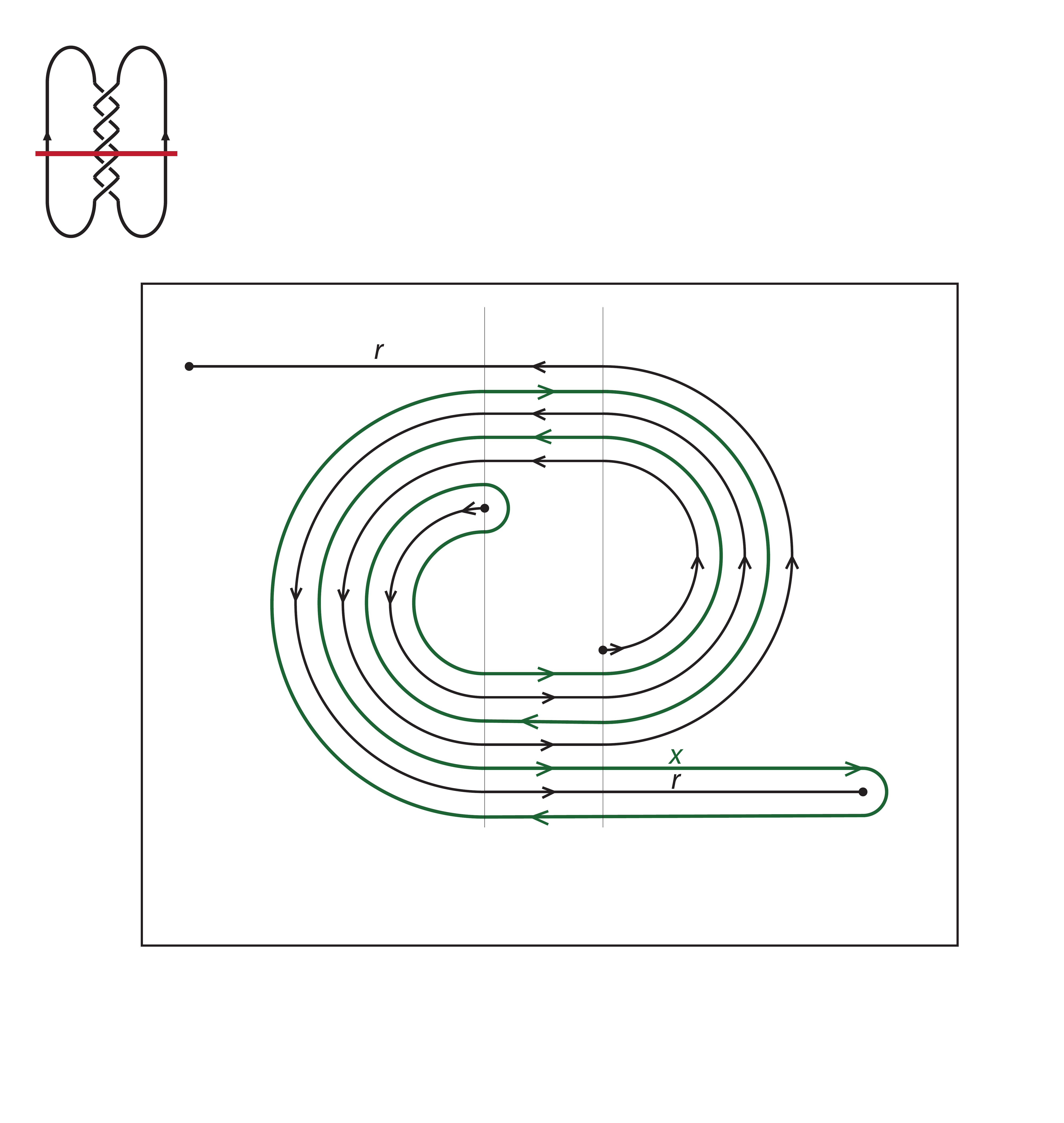}}
\]
\begin{center} {\bf Step 8} \end{center}

\newpage

Between Step 8 and Step 9, the inner black vertices have exchanged their positions in relation to the horizontal axis. The horizontal red arc on the knot diagram is  below the fourth crossing.
\[
\scalebox{.15}{\includegraphics{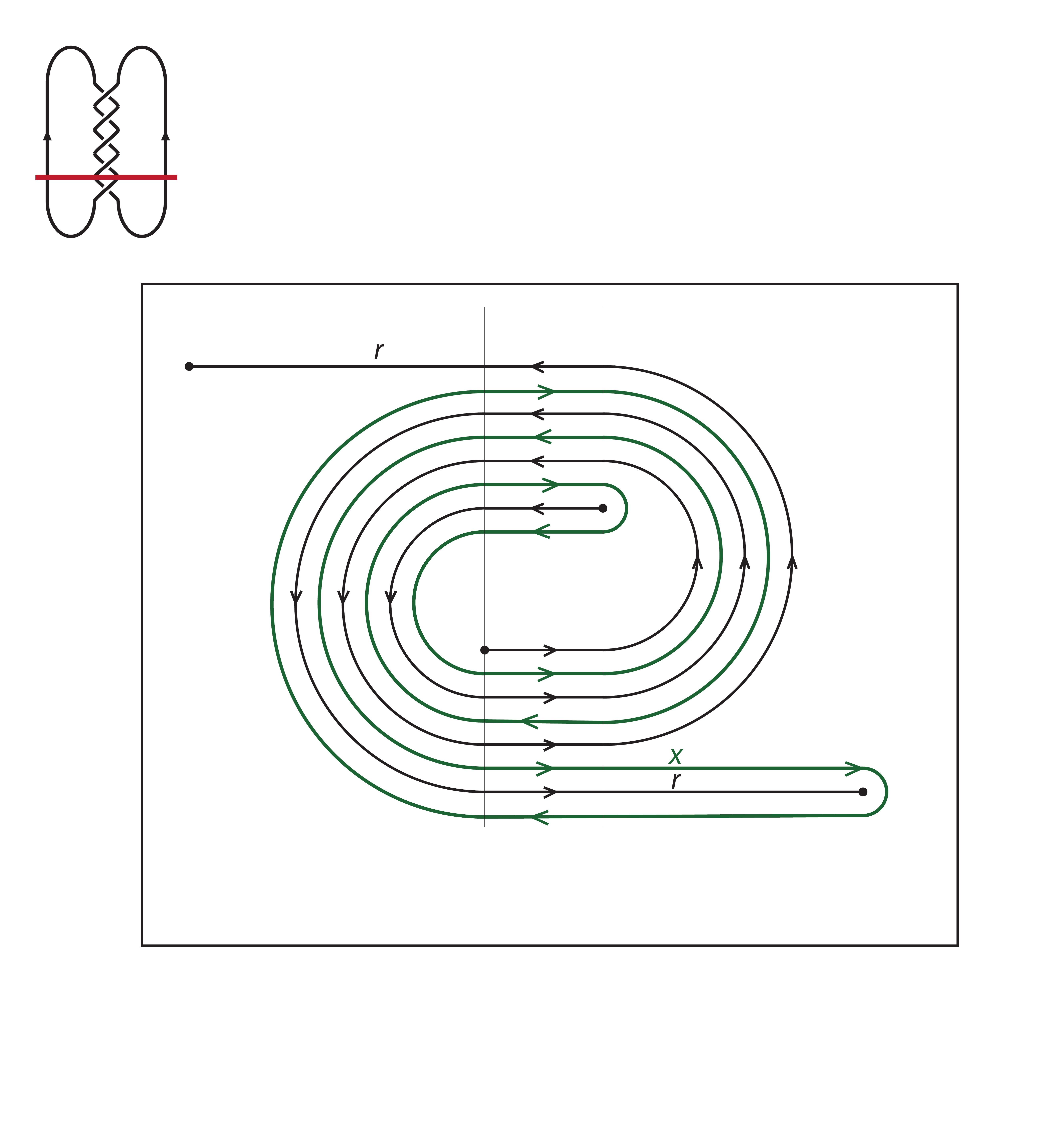}}
\]
\begin{center} {\bf Step 9} \end{center}

\newpage

 Between Step 9 and Step 10, semi-circular arcs on either side of the diagram are used to exchange the vertical positions of the inner black vertices. In this way, the vertex on the right will be able to pass in front of that on the left.
\[
\scalebox{.15}{\includegraphics{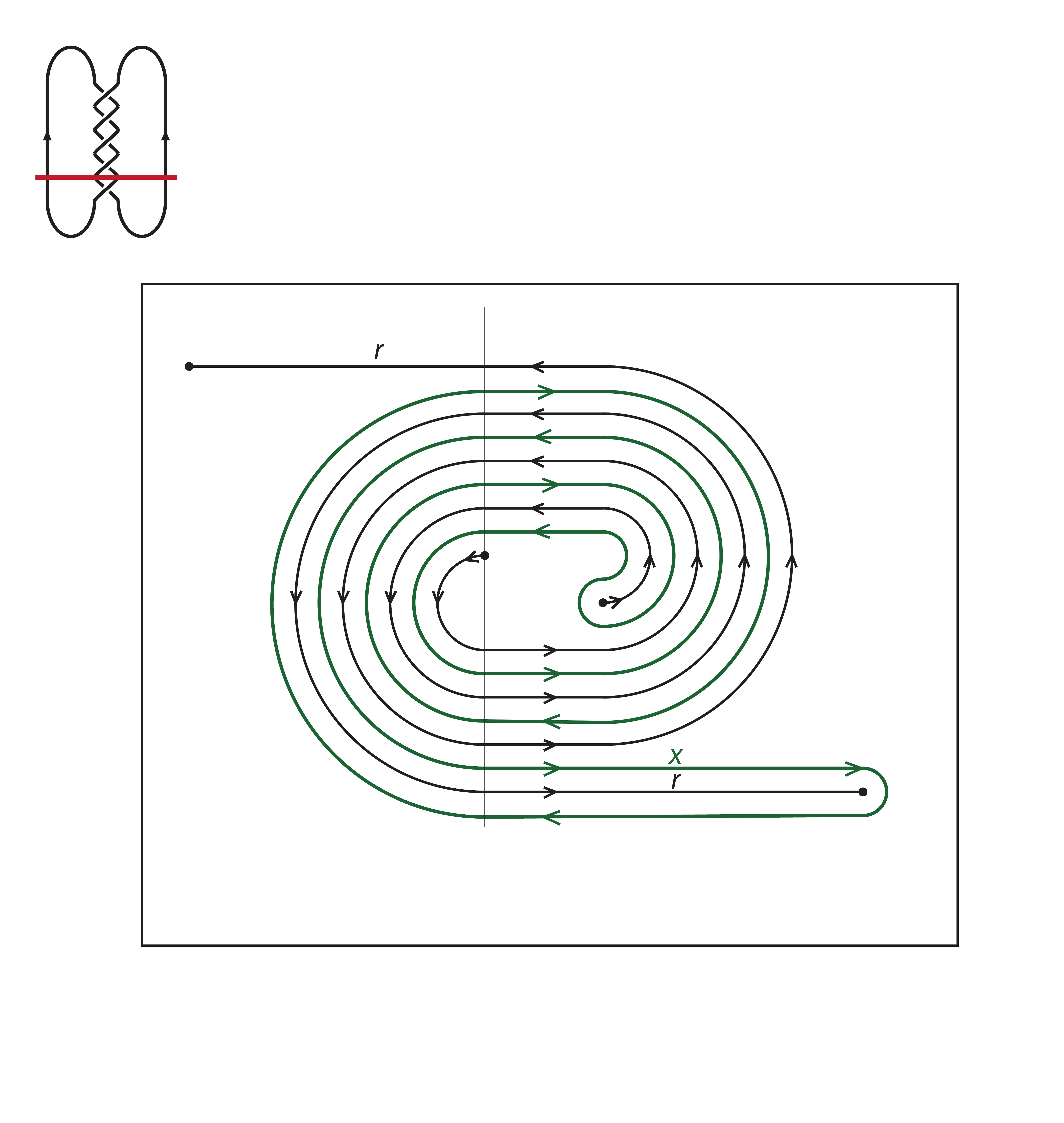}}
\]
\begin{center} {\bf Step 10} \end{center}

\newpage

Between Step 10 and Step 11, the inner black vertices have exchanged their positions in relation to the horizontal axis. The horizontal red arc on the knot diagram is  below the fifth crossing.

\[
\scalebox{.15}{\includegraphics{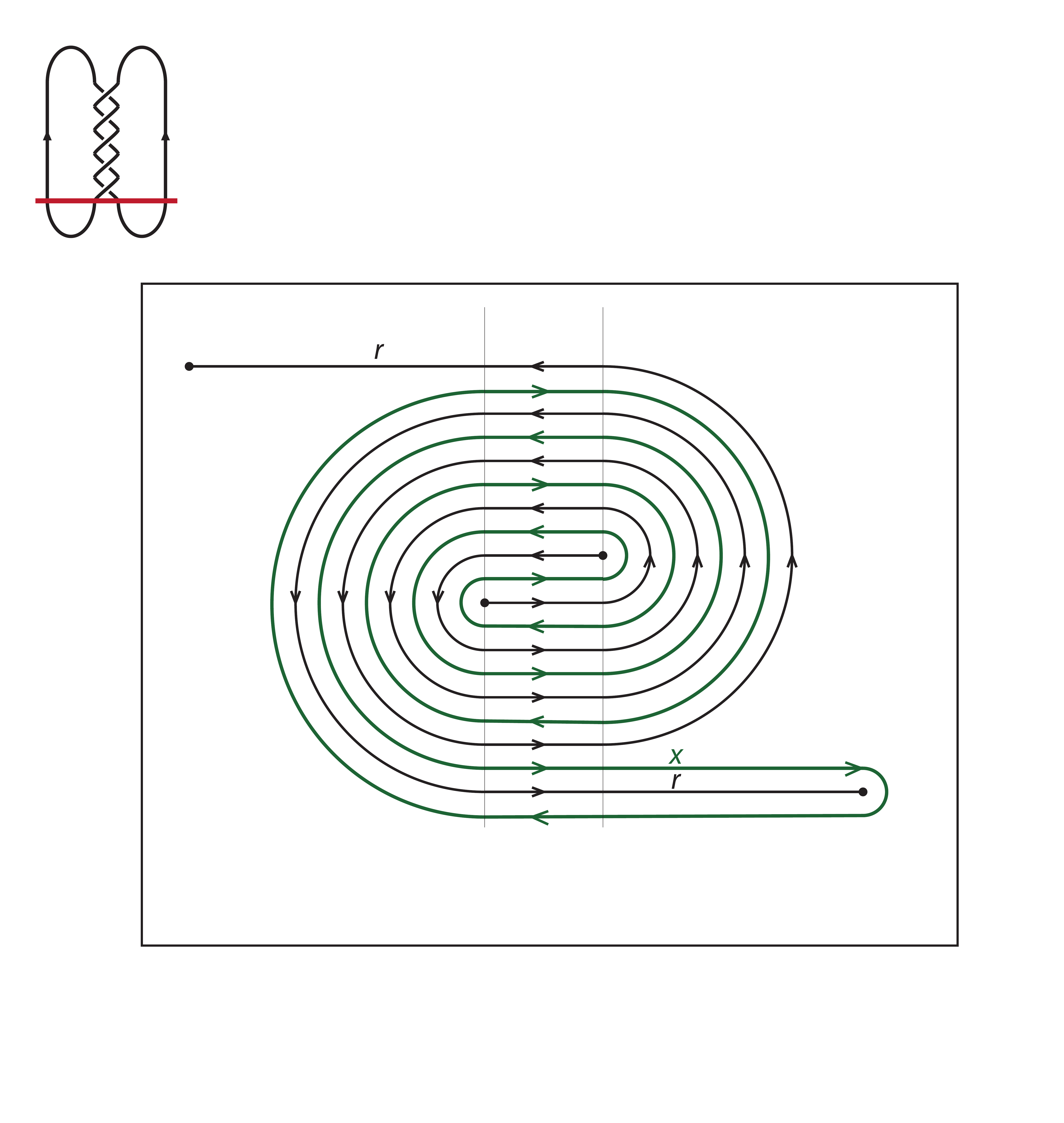}}
\]
\begin{center} {\bf Step 11} \end{center}

\newpage

The chart in Step 12 is isotopic to that in Step 11. The circular arcs have been replaced by vertical arcs. As topologists, we recognize that circles and squares are homeomorphic. 

\[
\scalebox{.15}{\includegraphics{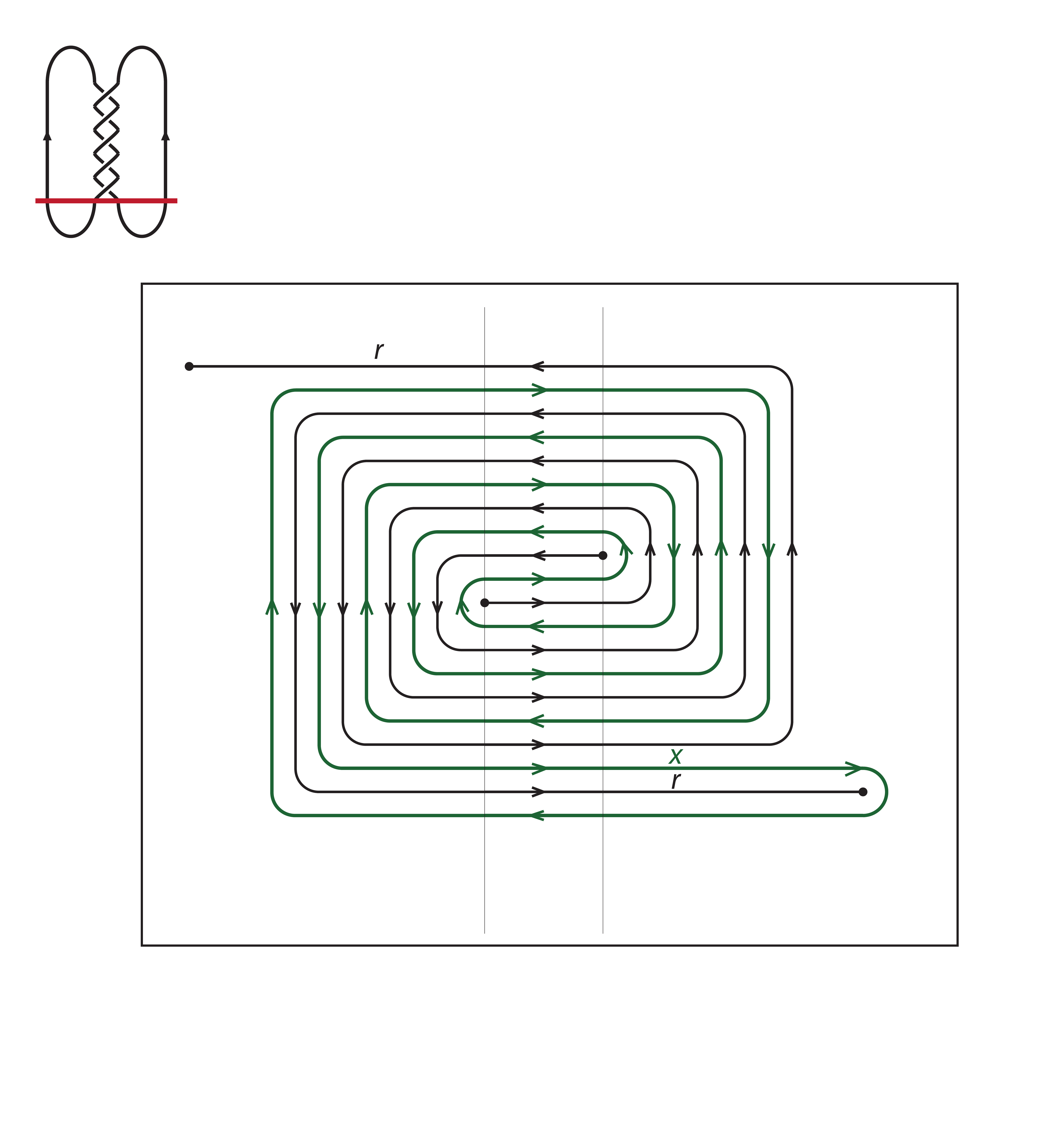}}
\]
\begin{center} {\bf Step 12} \end{center}

\newpage

Between Step 12 and Step 13, a canceling pair of valence four vertices along a pair of arc have been inserted. Such vertices will be added liberally in the illustrations that appear subsequently.

\[
\scalebox{.15}{\includegraphics{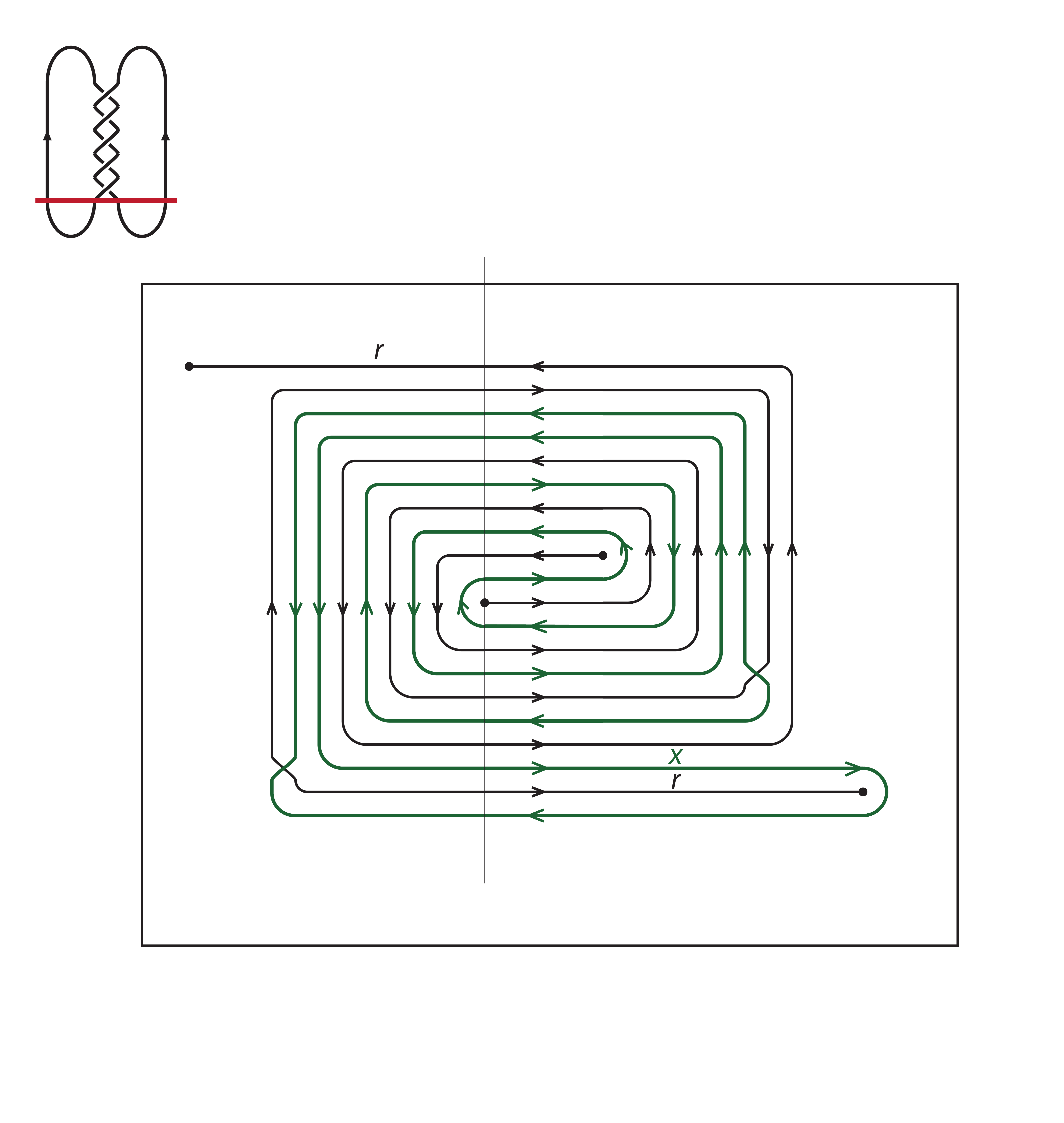}}
\]
\begin{center} {\bf Step 13} \end{center}

\newpage

Between Step 13 and Step 14, an anti-parallel pair of arcs that are labeled by $r$ are replaced with a pair of  ``turn-arounds." This move is analogous to a type II saddle move. Also two pairs of valence four vertices have been added.

\[
\scalebox{.15}{\includegraphics{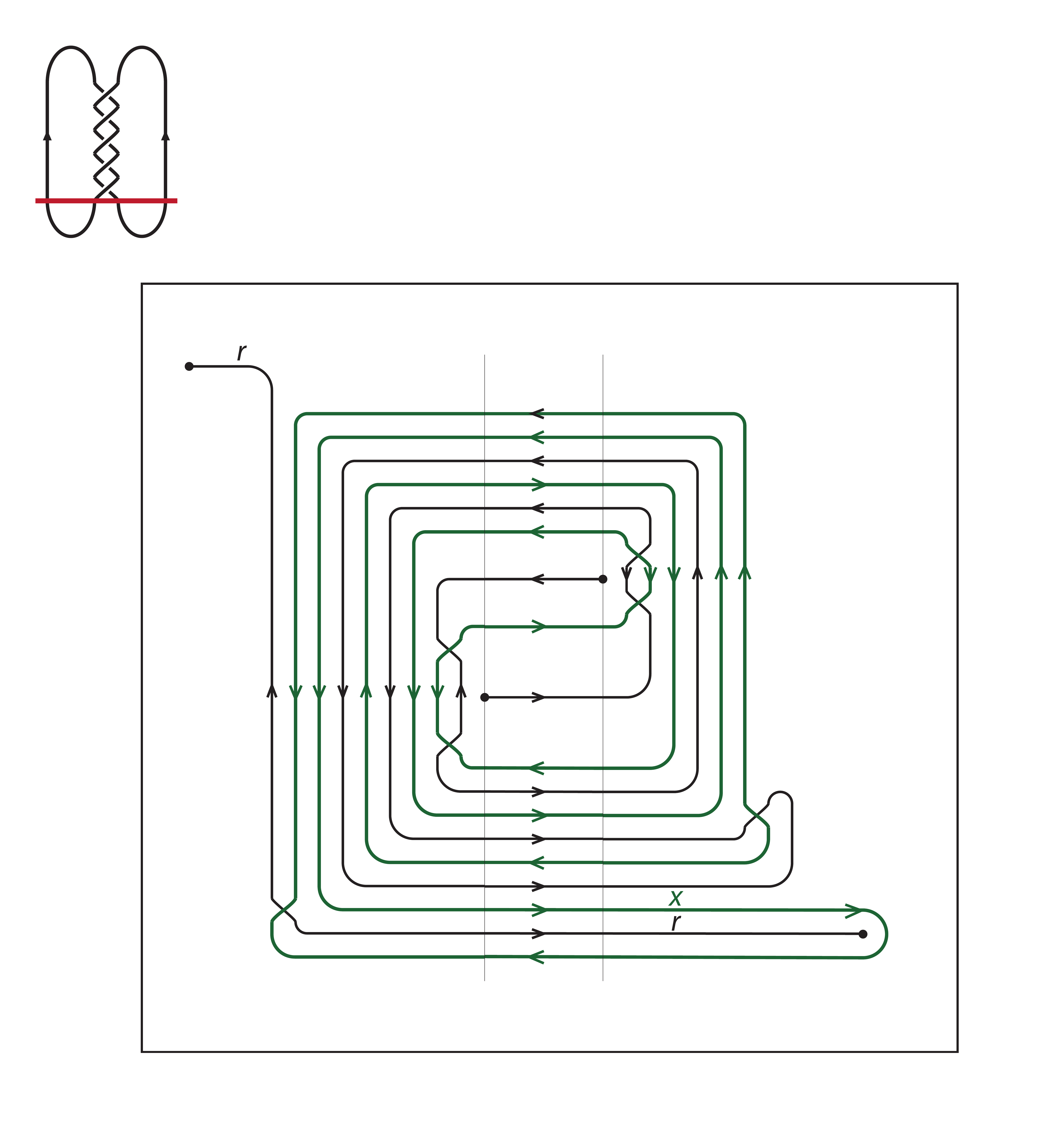}}
\]
\begin{center} {\bf Step 14} \end{center}

\newpage

Between Step 14 and Step 15, an isotopy of the chart has occurred. The bottom-right turn-around and valence four vertex has been moved to the left and upward. The lower, middle, left vertex has moved right and upward.

\[
\scalebox{.15}{\includegraphics{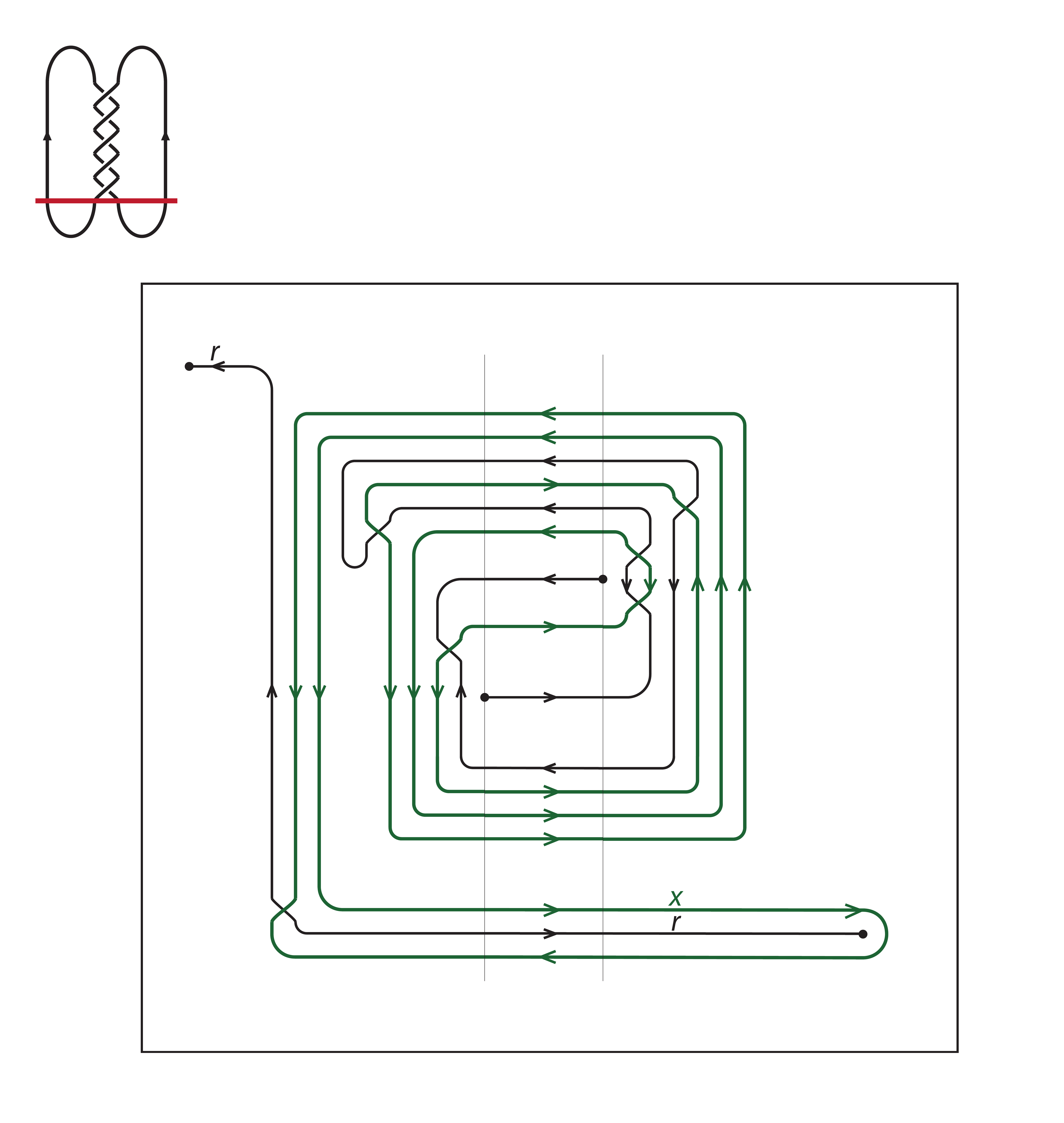}}
\]
\begin{center} {\bf Step 15} \end{center}

\newpage

Between Step 15 and Step 16, another pair of valence four vertices has been created, and the turn-around and vertex that moved previously has moved around the left corner and to the top right.

\[
\scalebox{.15}{\includegraphics{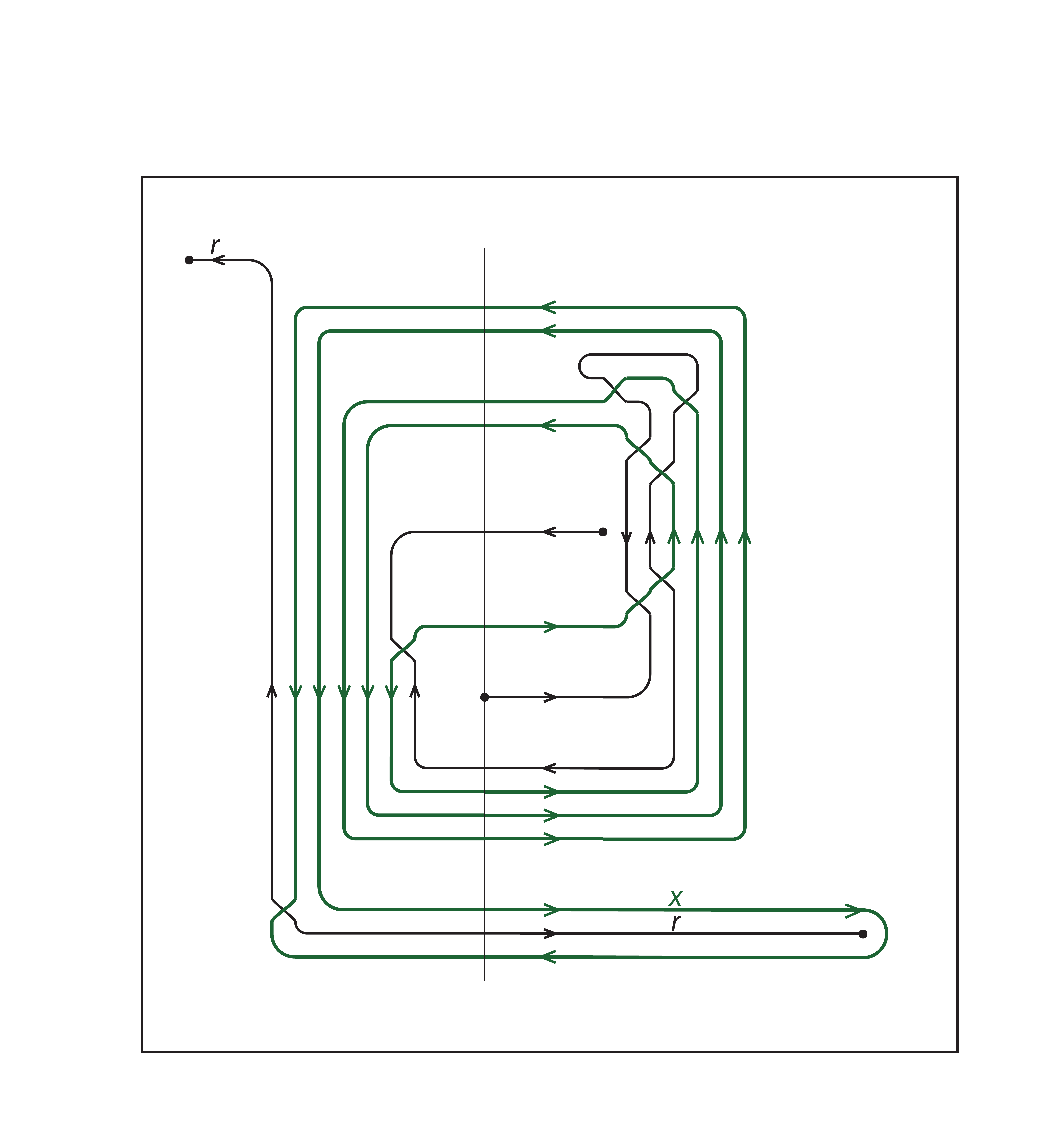}}
\]
\begin{center} {\bf Step 16} \end{center}

\newpage

Between Step 16 and Step 17, the six vertices on the right of the illustration cancel pairwise. The arc that is labeled by $x$ that surrounds the lower right black vertex is bent upward. 

\[
\scalebox{.15}{\includegraphics{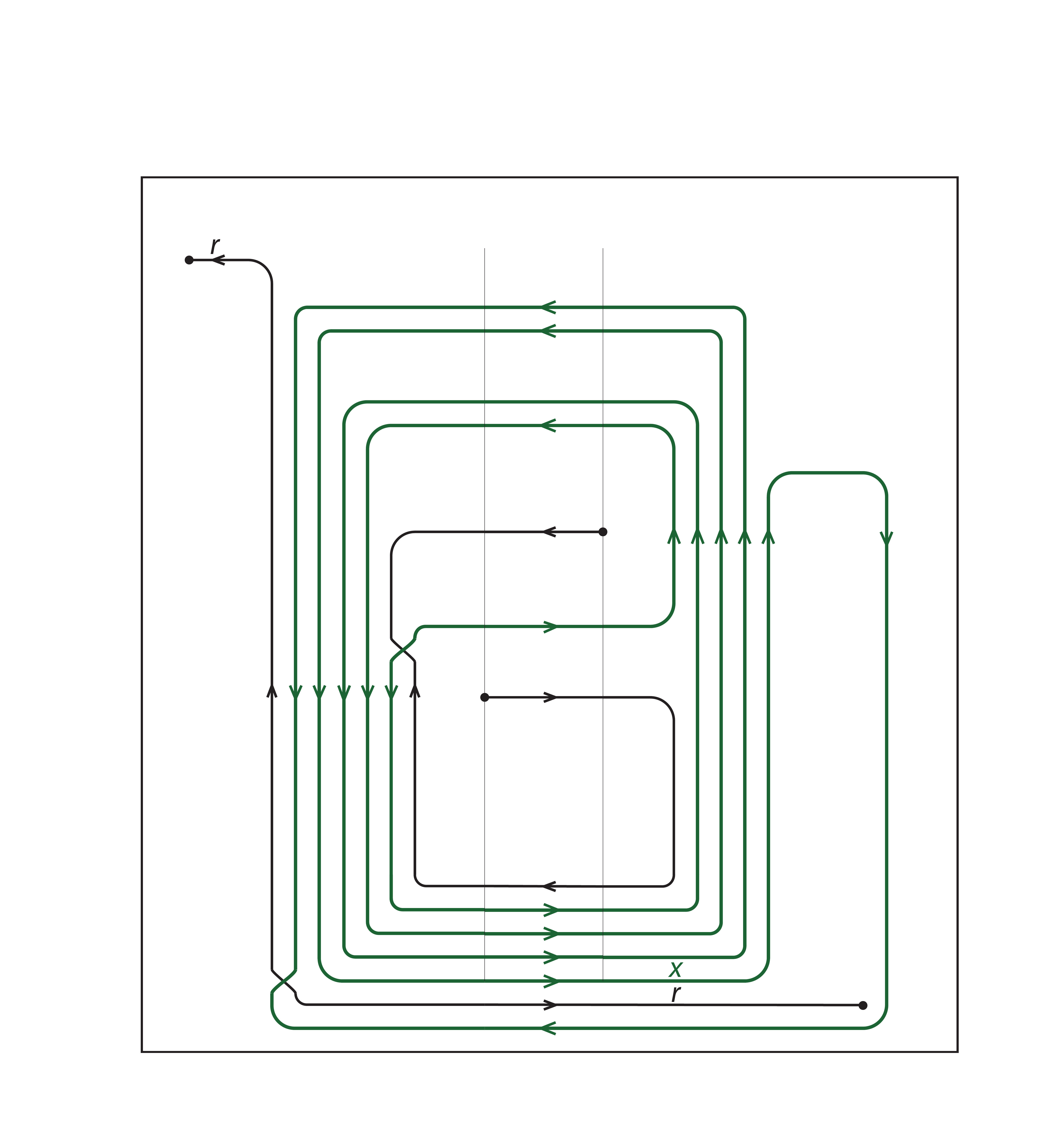}}
\]
\begin{center} {\bf Step 17} \end{center}

\newpage
Between Step 17 and Step 18, 
the  five parallel arcs labeled by $x$  have been split by adding two pairs of cancelling valence five vertices. The top left black vertex has also been repositioned vertically.

\[
\scalebox{.15}{\includegraphics{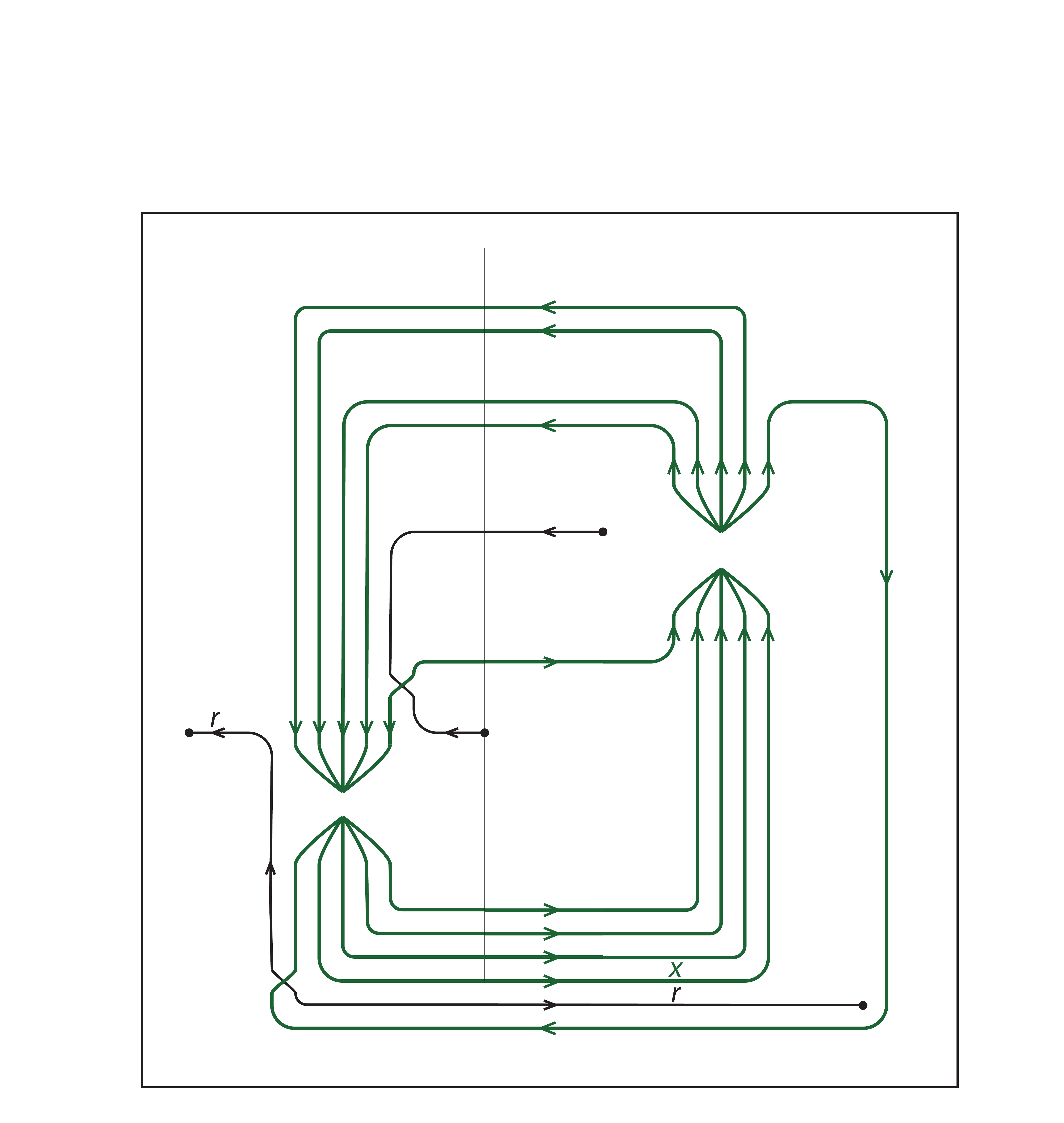}}
\]
\begin{center} {\bf Step 18} \end{center}

\newpage
Between Step 18 and Step 19, the channels that were opened by the addition of the $x^5=1$ relations have been used to reconnect the edges that are labeled $r$ using the ``turn-around" moves --- analogues of the type II saddle moves.

\vspace{-2cm}

\[
\scalebox{.15}{\includegraphics{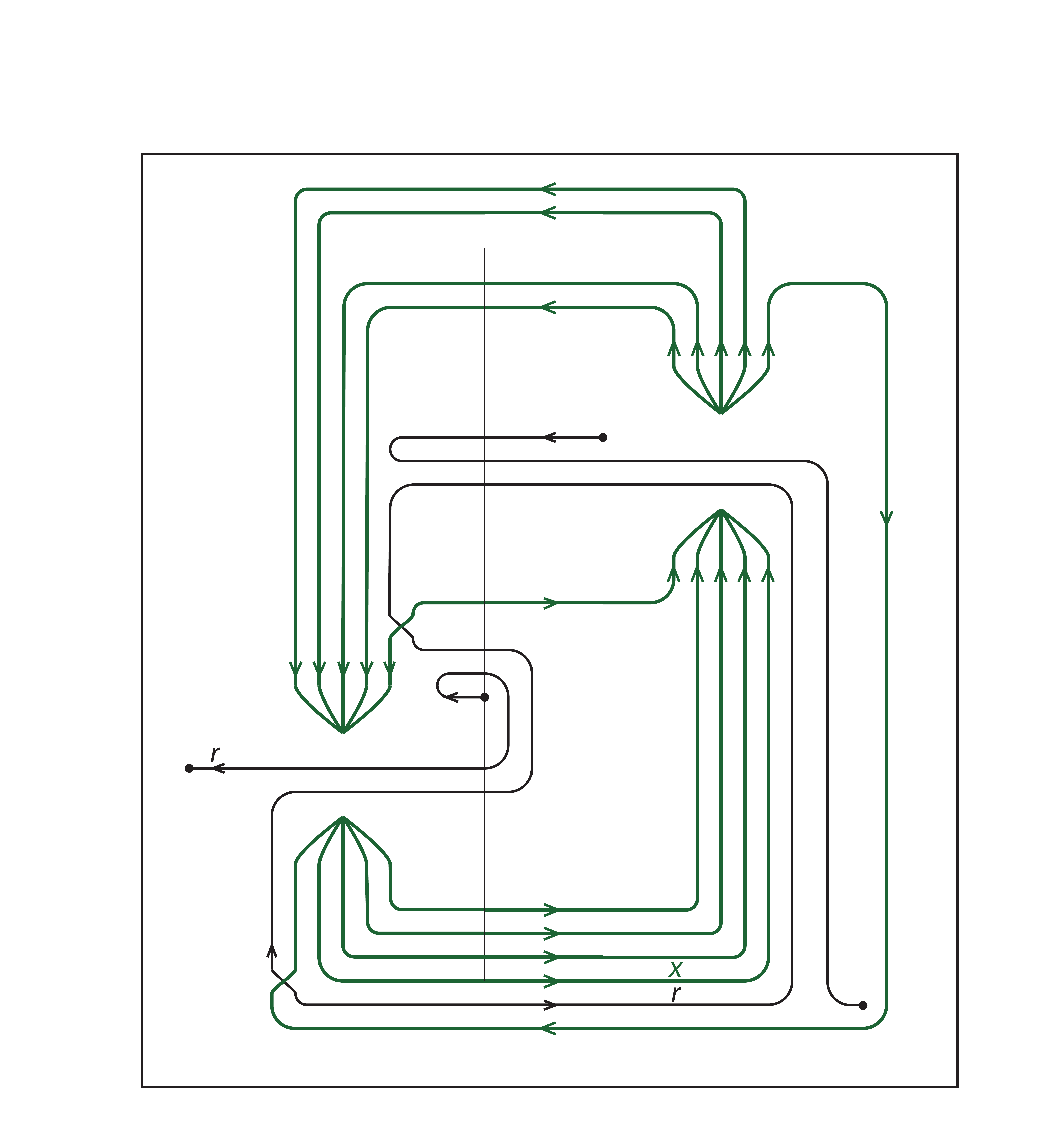}}
\]
\begin{center} {\bf Step 19} \end{center}

\newpage

Between Step 19 and Step 20, the arc labeled $r$ on the left has been removed. The knot diagram has been reintroduced, and the minimum on the left has been passed. 

\[
\scalebox{.15}{\includegraphics{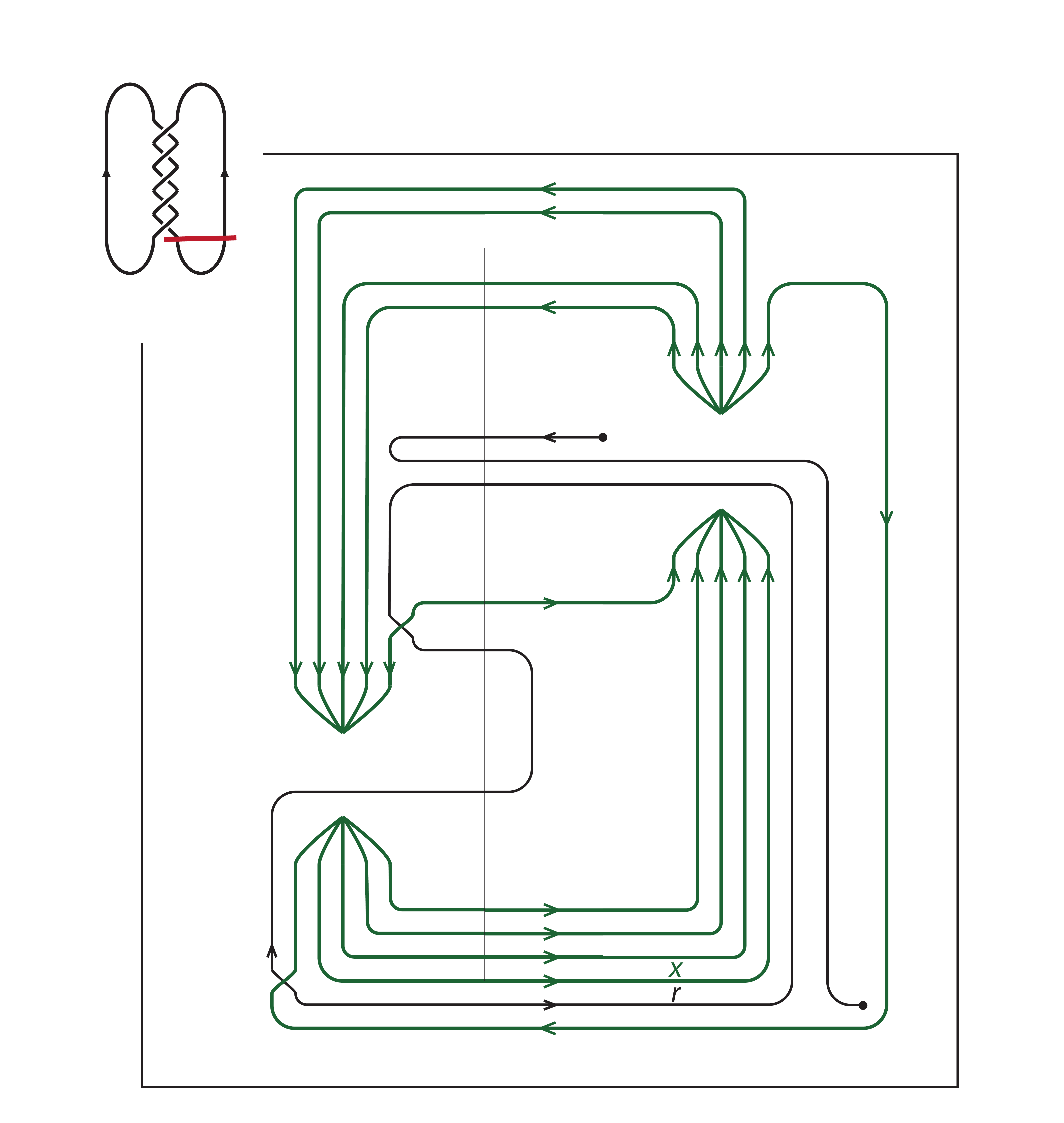}}
\]
\begin{center} {\bf Step 20} \end{center}

\newpage

Between Step 20 and Step 21, the arc labeled $r$ on the right cancels. The knot diagram has been reintroduced, and the minimum on the right has been passed. This dihedral chart represents an immersion of five $2$-spheres in $S^3$. 

\[
\scalebox{.15}{\includegraphics{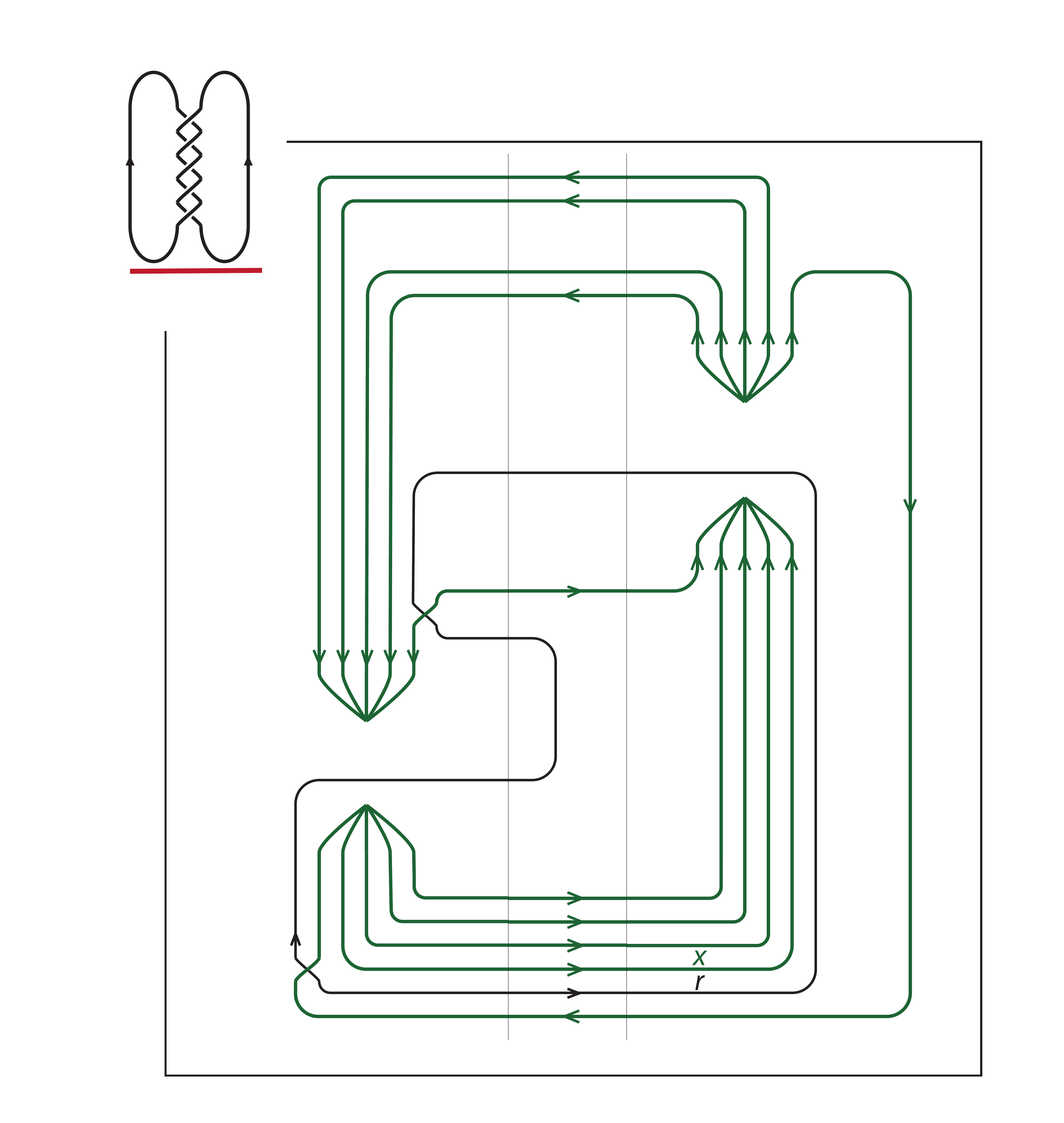}}
\]
\begin{center} {\bf Step 21} \end{center}

 \begin{figure}[htb]
\[
\scalebox{.2}{\includegraphics{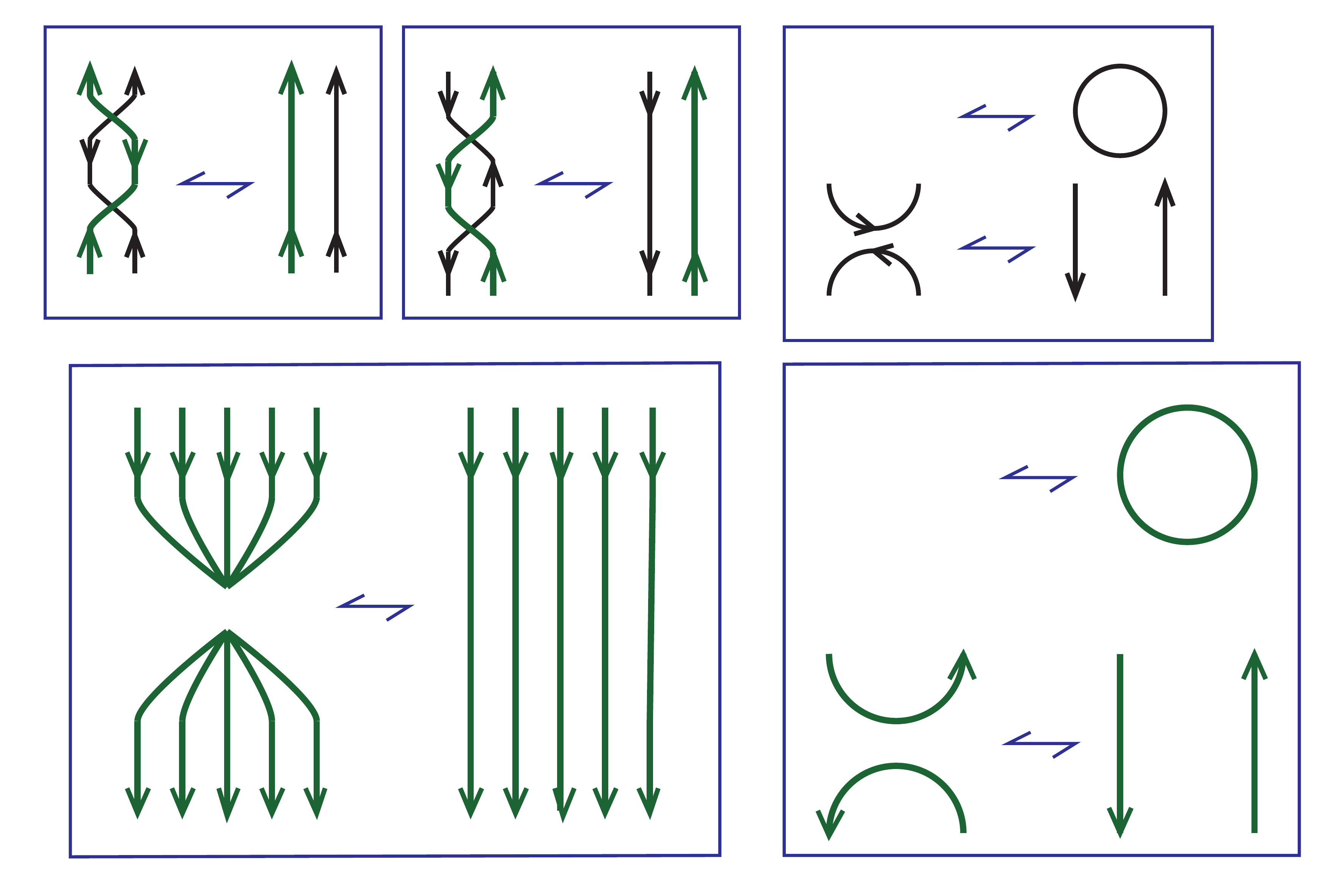}}
\]
\caption{Some moves that relate equivalent dihedral charts}
\label{dichartmoves}
\end{figure}

\section{Digression: Moves to Dihedral Charts / Moves to Permutation Charts}
\label{digress}

In the discussion above, there were essentially five differences between successive steps.

\begin{enumerate}

\item Adding or removing a simple arc labeled $r$ that connects a pair of black dihedral vertices. Addition: Step 0 to Step 1 and Step 1 to Step 2. Removal: Step 19 to Step 20 and Step  20 to Step 21.
\item Isotopies of the  graph's embedding.
\item Adding or removing a pair of  ``canceling" vertices.
\item Replacing a pair of  anti-parallel arcs that are labeled $r$ with a pair of ``turn-arounds."
\item Adding or removing a loop upon which no vertices appear. This occurred without fanfare between Steps 1 and 2. In our conception of the process (that has thus far been unstated) the loop labeled $x$ was created before the edge, labeled $r$,  that it surrounds was added. Perhaps we omitted a Step $1.5$. 
\end{enumerate}
The non-isotopy moves to charts that were used are depicted in Fig.~\ref{dichartmoves}.

 \begin{figure}[htb]
\[
\scalebox{.15}{\includegraphics{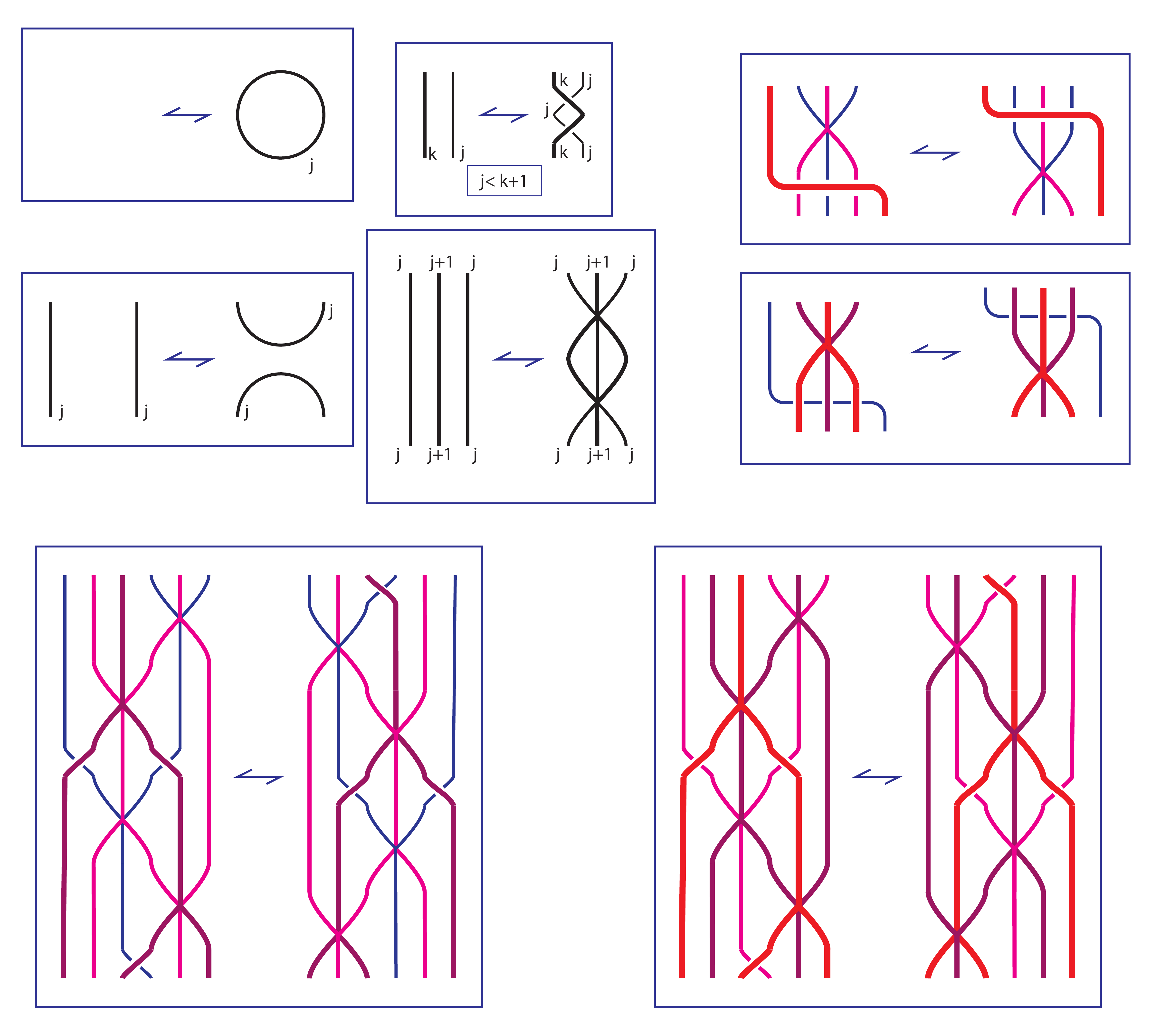}}
\]
\caption{The moves for $\Sigma_5$ charts}
\label{almostKmoves}
\end{figure}

With the exception of adding or removing a simple arc, there are  some justifications for the remaining changes.

\subsection{A category that is associated to a group presentation}
\label{coolcat}
In general, suppose that a group presentation $G=\langle X:R\rangle$ is given. There is a category that can be constructed as follows. We call this category ${\mathcal {C}}(G)$ although it depends upon the presentation $\langle X:R\rangle$. The objects in the category are unreduced words in the free group $F(X)$ on the generating set $X$. If two words represent different elements in $G$, then the hom-set between them is empty. If the words represent the same element in $G$, then the set of morphisms between them is generated by the relations in $R$ as well as the trivial relations $aa^{-1}=1,$ and $a^{-1}a=1$.  This category has a monoidal structure that is given by the juxtaposition of words. Thus a word, as a source object, is arranged, as usual, along the bottom horizontal edge of a rectangular disk. And the target object is arranged along the top horizontal edge. A morphism is a diagram in the interior of the rectangle.  In this way, the identity morphism from a word to itself is represented by a collection of vertical arcs between the corresponding letters in the identical words. 
The generating morphisms are invertible by declaration. That declaration is precisely the set of moves that are depicted in Fig.~\ref{dichartmoves} in the case of the dihedral group.

A word in the free group $F(X)$ is essentially the same as a loop that is confined to the $1$-skeleton of the classifying space $BG= K(G,1)$.  A morphism is a homotopy of such a loop. The $\cup$ and $\cap$ morphisms correspond to the relations $1=aa^{-1}, $ $1=a^{-1}a$, $a^{-1}a=1$, and $ aa^{-1}=1$. They are realized by the usual homotopy from ``wind and then unwind'' to ``remain constant." Higher order homotopies can be constructed that indicate the invertibility of these $\cup$s or $\cap$s operations. Whence one might write $\cup^{-1} =\cap$. Similarly, the other moves to dihedral charts can be thought of as higher homotopies between sliding back and forth over a relator disk in the classifying space.  The idea of mapping the knot exterior into $BG$ is suggested, in part, by Hilden's discussion \cite{HildenPac78}.

There are additional identities among relations for the dihedral group, but we won't need them here. These identities and those in Fig.~\ref{dichartmoves} will be seen in Lemmas~\ref{MainLemma1} and~\ref{MainLemma2} as consequences of the corresponding relations in the symmetric group. 

The structure of the classifying space, at least through the $3$-skeleton, of the symmetric group $\Sigma_5$ is known to us. Four loops  that correspond to the generators $t_1=(12)$ through $t_4=(45)$ are attached at the unique vertex. The $2$-cells  are planar duals to crossings and white vertices (the unoriented versions of those depicted in Fig.~\ref{BraidChart}).  The $3$-cells and the other identities among relations are depicted in Fig.~\ref{almostKmoves}. 

Lemma~\ref{MainLemma1}  and Lemma~\ref{MainLemma2}  will be stated and proved using the graphical notation. Lemma~\ref{MainLemma2}   and  Corollary~\ref{cor1} mean that the morphisms in the category ${\mathcal {C}}(D_5)$ are invertible. That the $\cup$ and $\cap$ morphisms in $D_5$ are invertible, and mutually inverse, is  a consequence of the corresponding fact in $\Sigma_5$. The invertibility of $\cup$s and $\cap$s in $\Sigma_5$ is true because the corresponding  disks in $\R^3$ are related by type II bubble and saddle moves. More specifically, we take the moves depicted in Fig.~\ref{almostKmoves} as given.

\section{Lemmas}
\label{thelemmas}

\begin{figure}[htb]
\[
\scalebox{.2}{\includegraphics{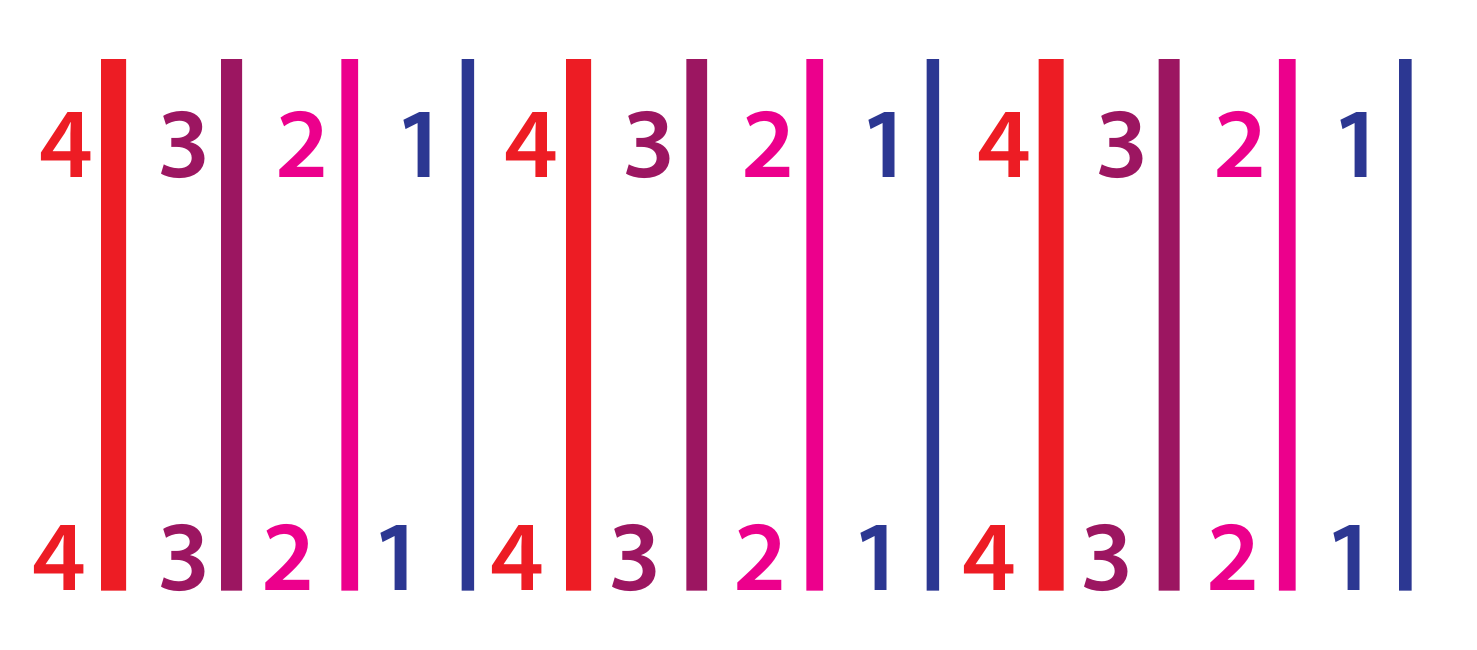}}
\]
\caption{The left-hand side of Lemma~\ref{MainLemma1}}
\label{Lemma1Left}
\end{figure}

\begin{lemma}
\label{MainLemma1}
\[
\scalebox{.3}{\includegraphics{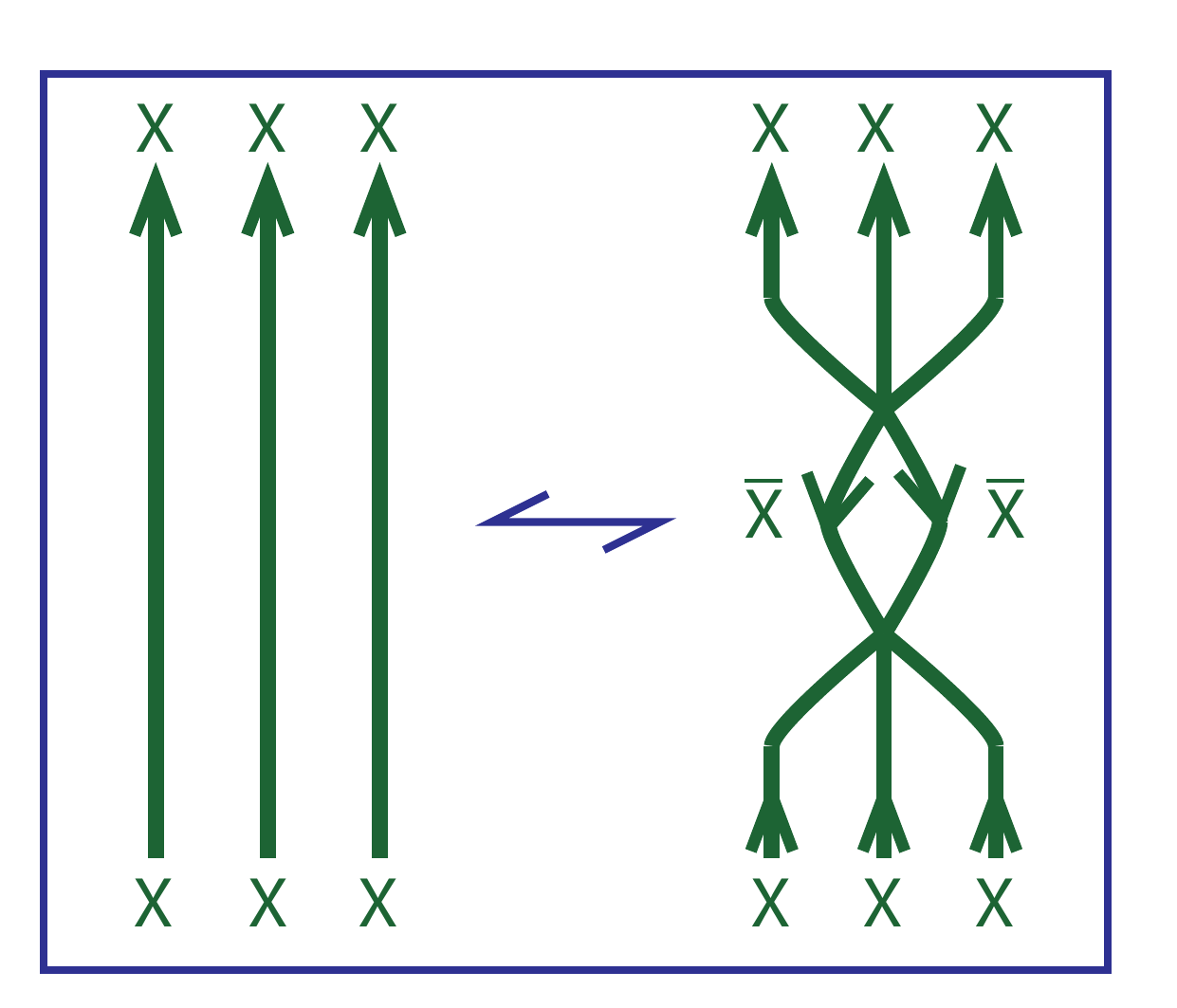}}
\]
\end{lemma}
The meaning of the statement of the lemma is as follows.  Either diagram  in the statement represents a family of  five rectangularly shaped disks that intersect amongst each other in the $3$-dimensional space that is $D \times [0, 6]$ where $D$ is an isometric copy of one of these rectangles. These disks are {\it intertwined}. That is, each is embedded, but along the bottom and top edges of the disks the permutation pattern $(15432)^3$ appears. That $5$-cycle $(15432)=(45)(34)(23)(12)=t_4t_3t_2t_1$ indicates that the first disk successively intersects the second through fifth disks, and that pattern continues thrice. Compare that description with the second illustration on the right-hand, top, side of Fig.~\ref{legend}.  On the left-hand side of the equivalence, that pattern of  $(15432)^3$  continues without interruption as the  identity morphism on the element $(45)(34)(23)(12)(45)(34)(23)(12)(45)(34)(23)(12)$ in the category ${\mathcal {C}}(\Sigma_5).$ See Fig.~\ref{Lemma1Left}.

At the level of permutation cycles, one can compute that  \[ (15432)(15432)(15432)=(13524) =(12345)(12345).\] That algebraic calculation has a visual or diagrammatic counterpart that involves the crossing diagrams that we have labeled $1$ through $4$ (appropriately colored). The diagrammatic description of this identity involves the white vertices and crossings of the permutation charts. 

\begin{figure}[htb]
\[
\scalebox{.08}{\includegraphics{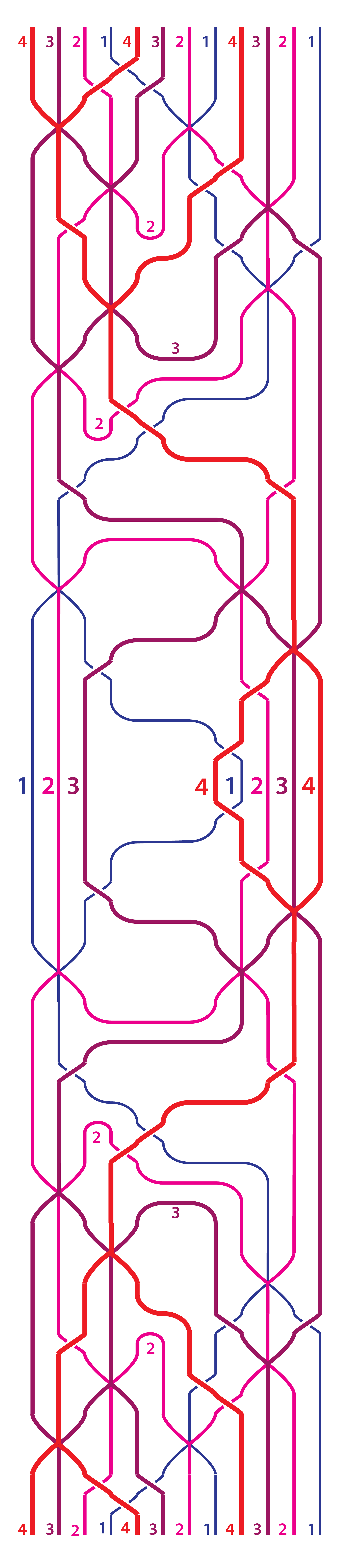}}
\]
\caption{The right-hand side of Lemma~\ref{MainLemma1}}
\label{Lemma1Right}
\end{figure}

Specifically, the five intertwined disks project to a planar diagram in which, at its center, the transposition sequence appears as\[(12)(23)(34)(45)(12)(23)(34)(45).\] At the top and bottom of the diagram, the pattern $(15432)^3$ appears. Between the vertical extremes and the center are white vertices and crossings. Thus the diagram on the right of the statement is a synopsis of the illustration in Fig.~\ref{Lemma1Right}.

The proof of the lemma will consist of a collection of diagrams that are related, one to the next, by applications of some of the moves that appear in Fig.~\ref{almostKmoves}. At the level of the classifying space $B\Sigma_5$ the isotopy of the  loop $(t_4t_3t_2t_1)^3(t_4t_3t_2t_1)^{-3}$  to the identity is isotopic to a loop that slides over a collection of hexagons (for the white vertices) and squares (for the crossings) until it reaches the path $(t_1t_2t_3t_4)^2$ and then slides back over the same sequence of hexagons and squares. The extreme isotopic disks are represented by the diagrams in Figs.~\ref{Lemma1Left} and ~\ref{Lemma1Right}. The isotopy of disks between them is the sequence of intertwined disks represented in the proof.


\clearpage

\noindent
{\sc Proof of of Lemma~\ref{MainLemma1}.}

\[
\raisebox{1.2cm}{
\scalebox{.23}{\includegraphics{Lemma1Left.pdf}}}
\raisebox{2cm}{
{\Huge{$\leftrightharpoons$}}}
\scalebox{.18}{\includegraphics{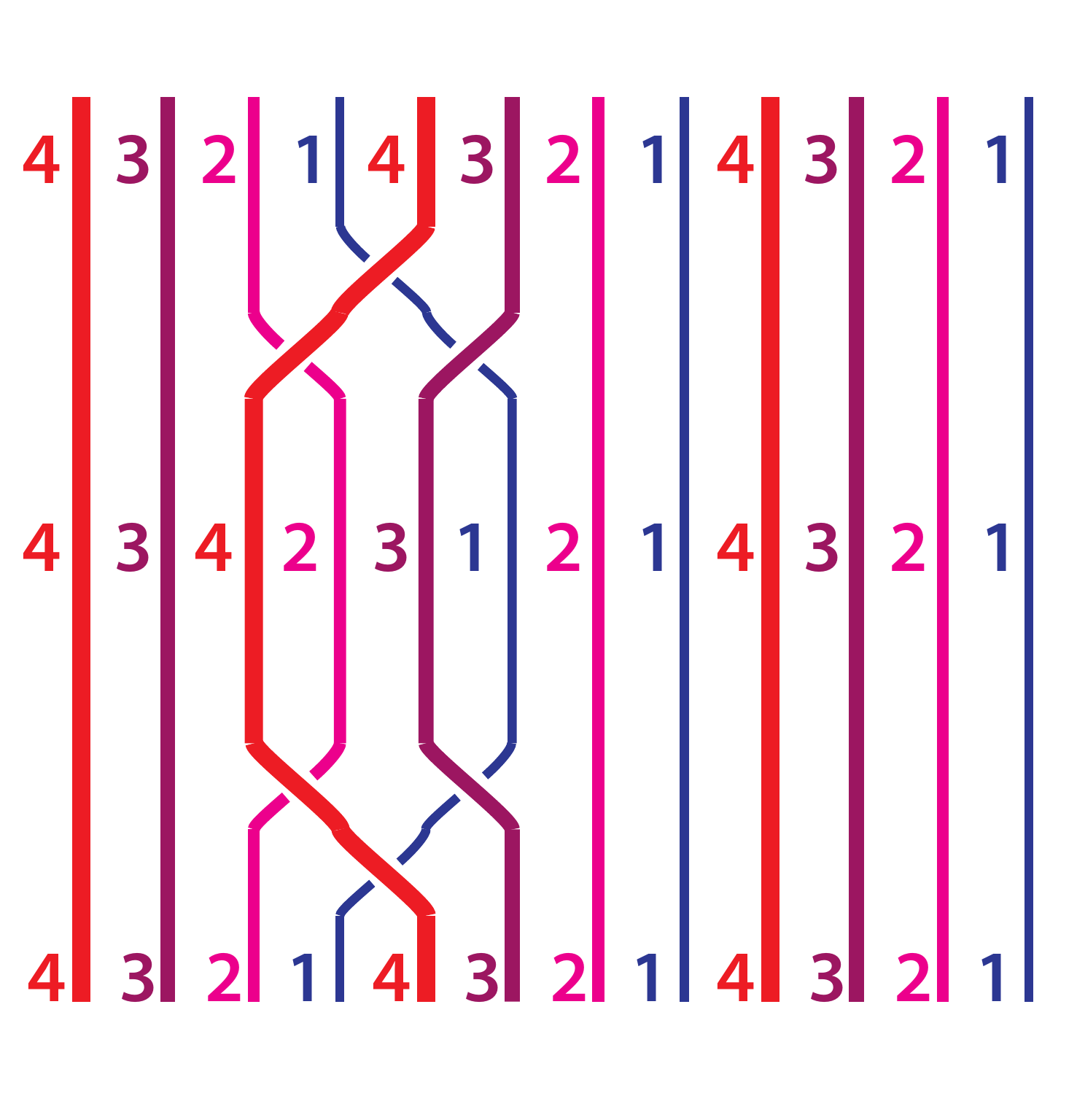}}
\raisebox{2cm}{
{\Huge{$\leftrightharpoons$}}}
\]

\[
\raisebox{1.5cm}
{\scalebox{.18}{\includegraphics{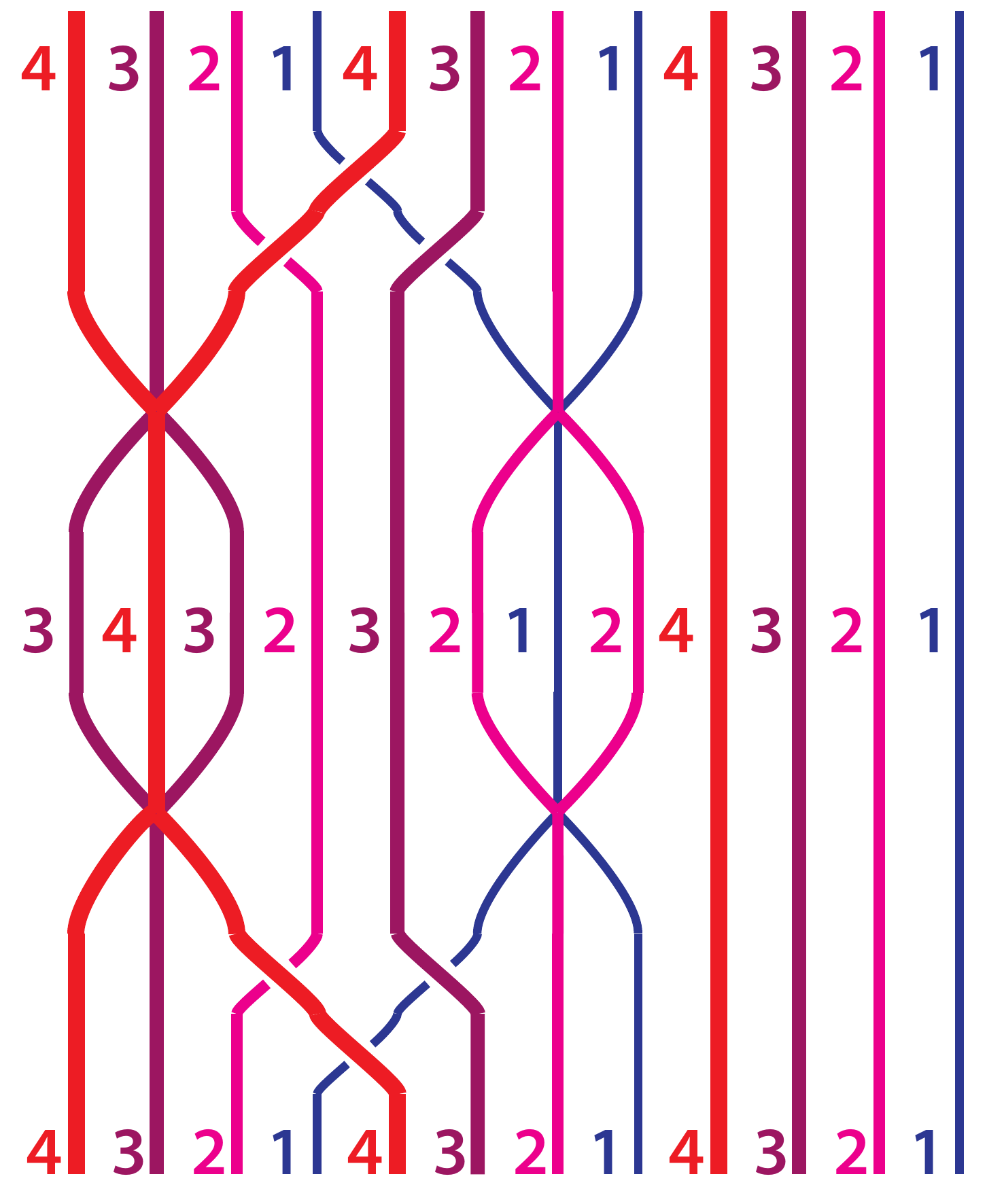}}}
\raisebox{3.75cm}{
{\Huge{$\leftrightharpoons$}}}
\scalebox{.18}{\includegraphics{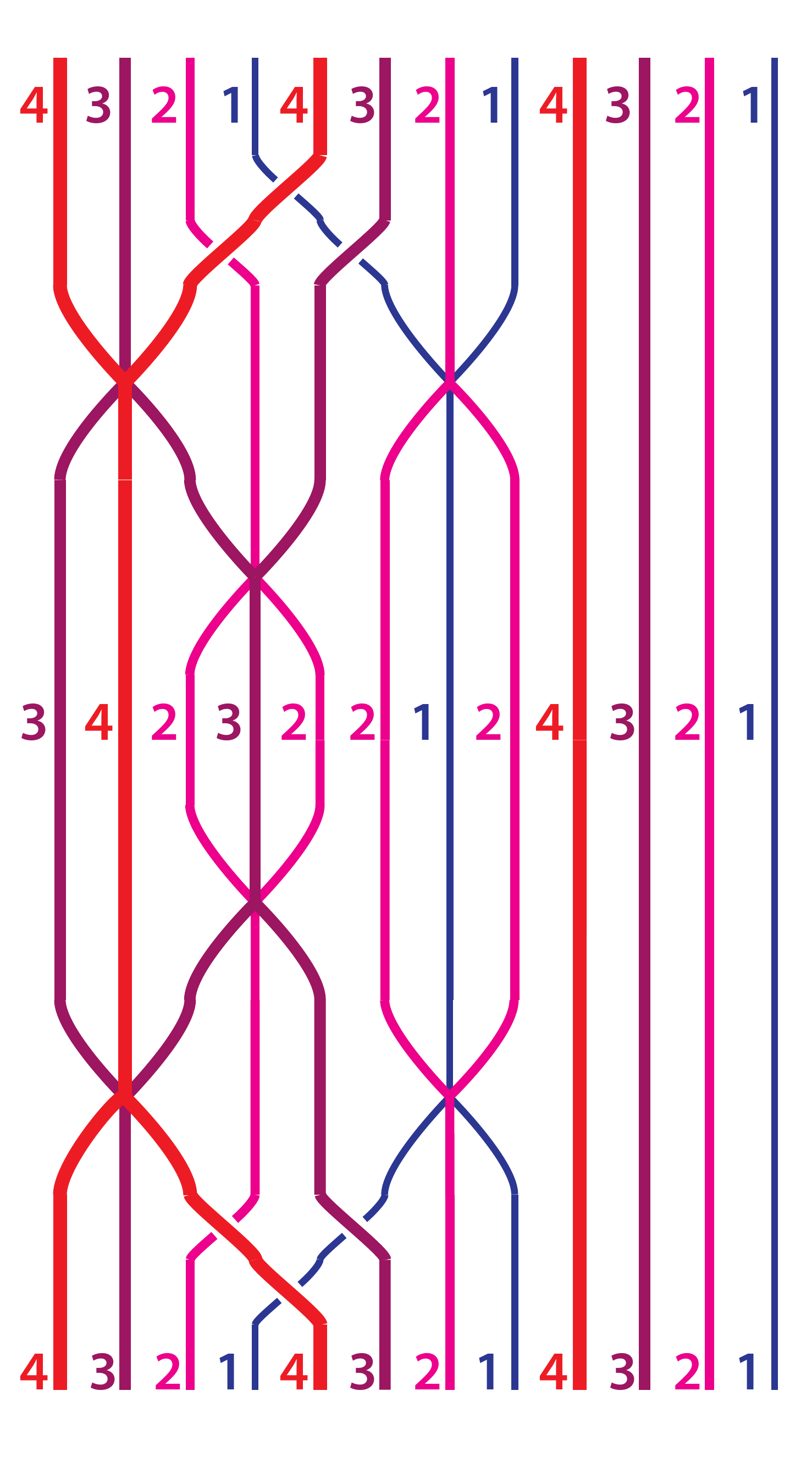}}
\raisebox{3.75cm}{
{\Huge{$\leftrightharpoons$}}}
\]

\[
\raisebox{.75cm}
{\scalebox{.18}{\includegraphics{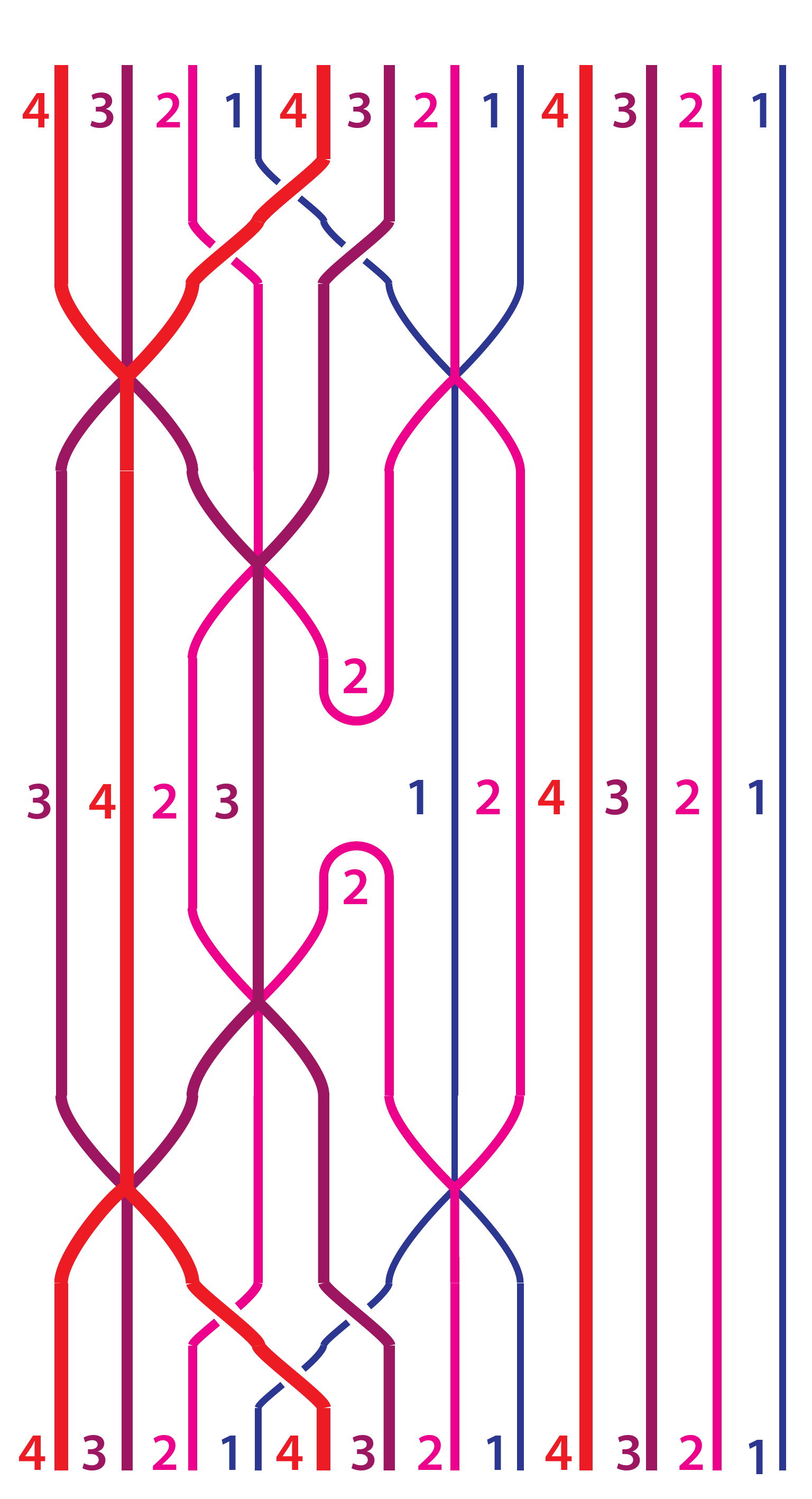}}}
\raisebox{4.6cm}{
{\huge{$\leftrightharpoons$}}}
\scalebox{.18}{\includegraphics{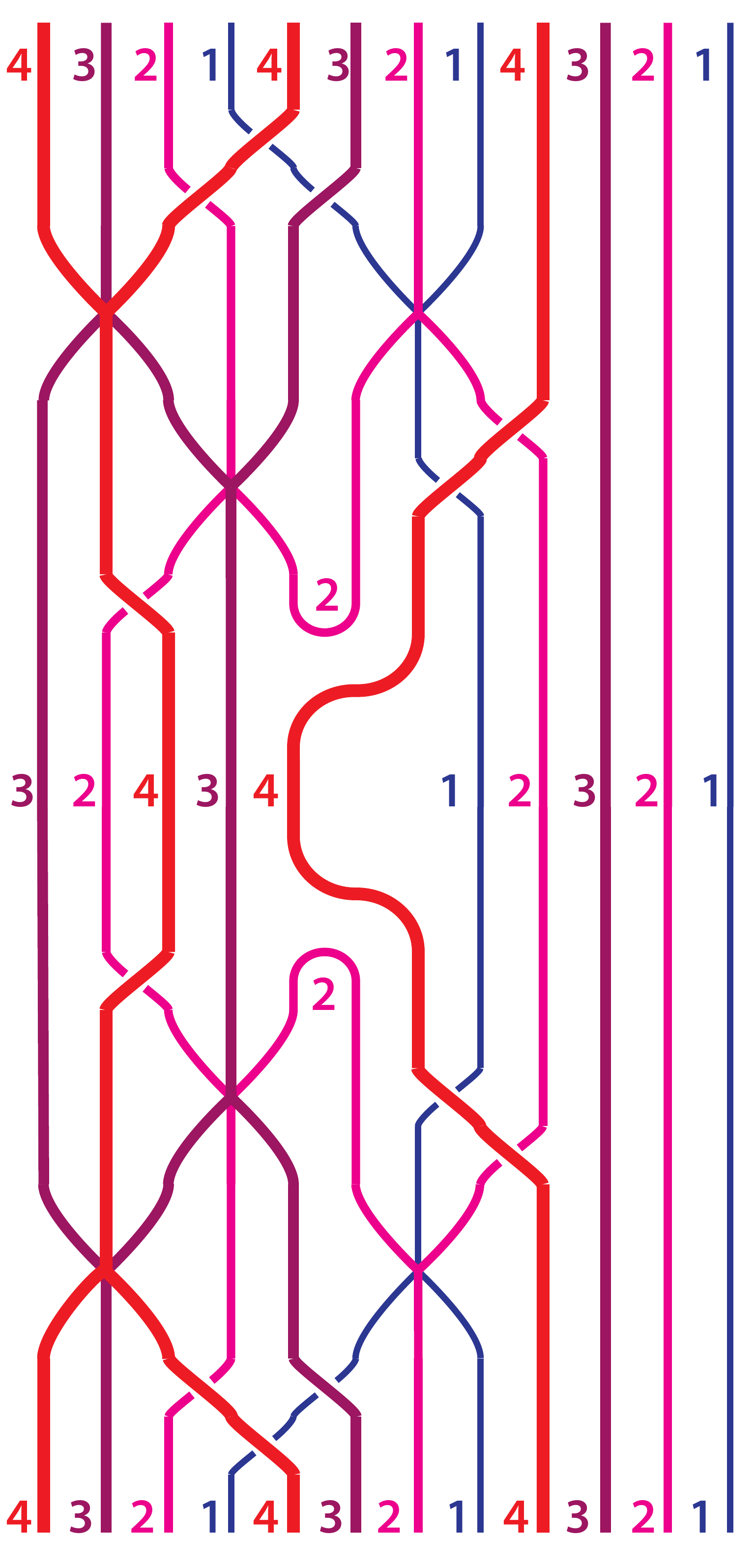}}
\raisebox{4.6cm}{
{\huge{$\leftrightharpoons$}}}
\]

\[
\raisebox{6.25cm}{
{\huge{$\leftrightharpoons$}}}
\scalebox{.18}{\includegraphics{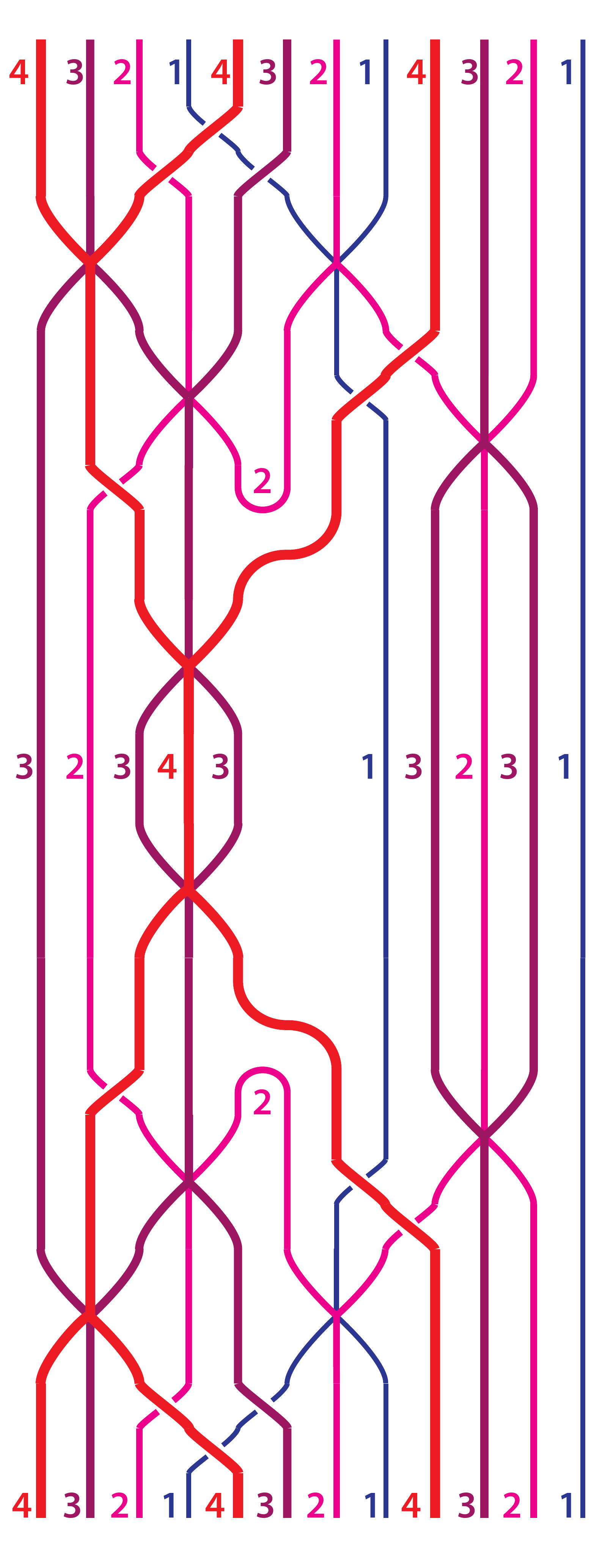}}
\raisebox{6.25cm}{
{\huge{$\leftrightharpoons$}}}
\raisebox{-.15cm}
{\scalebox{.18}{\includegraphics{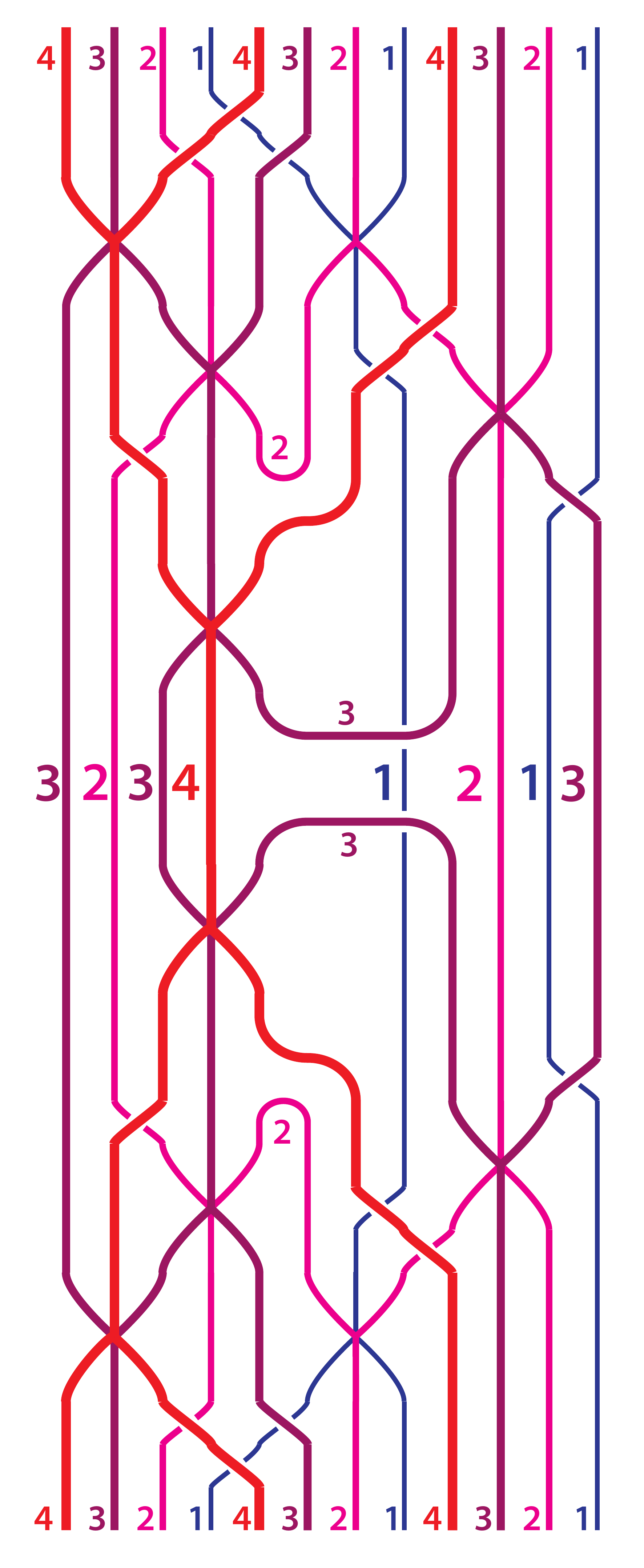}}}
\]

\[
\raisebox{7cm}{
{\huge{$\leftrightharpoons$}}}
\scalebox{.18}{\includegraphics{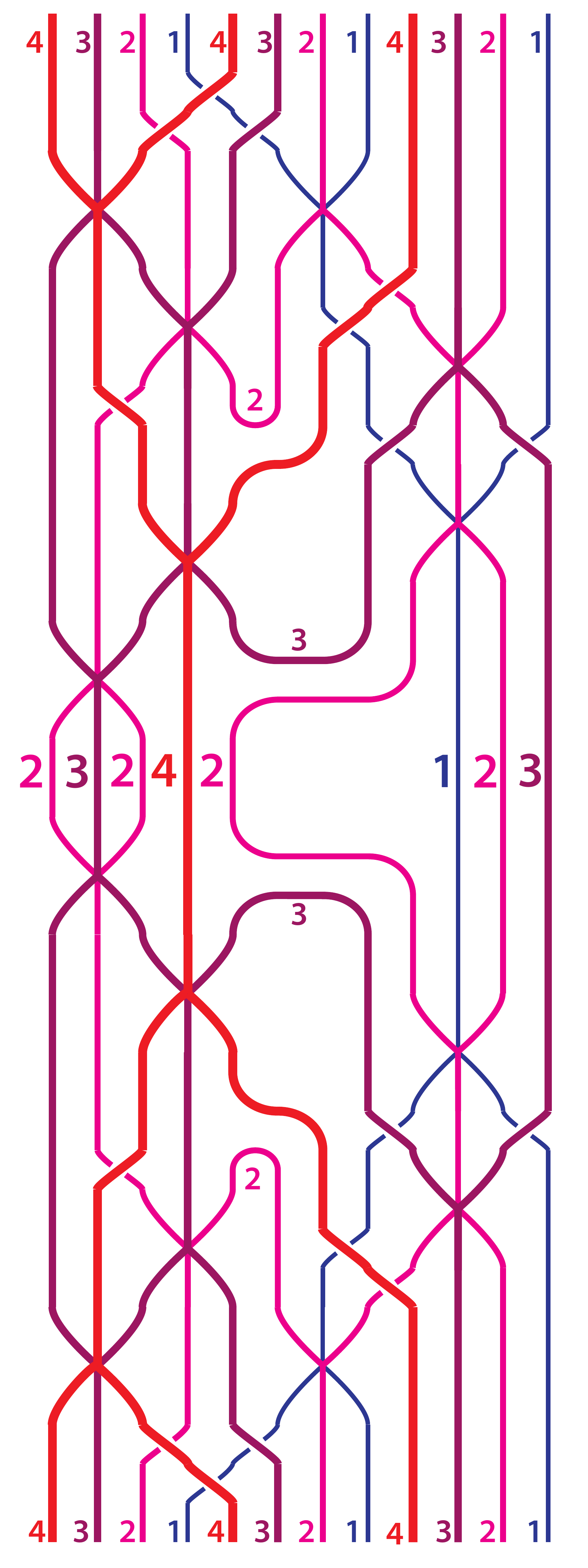}}
\raisebox{7cm}{
{\huge{$\leftrightharpoons$}}}
\raisebox{-.9cm}
{\scalebox{.18}{\includegraphics{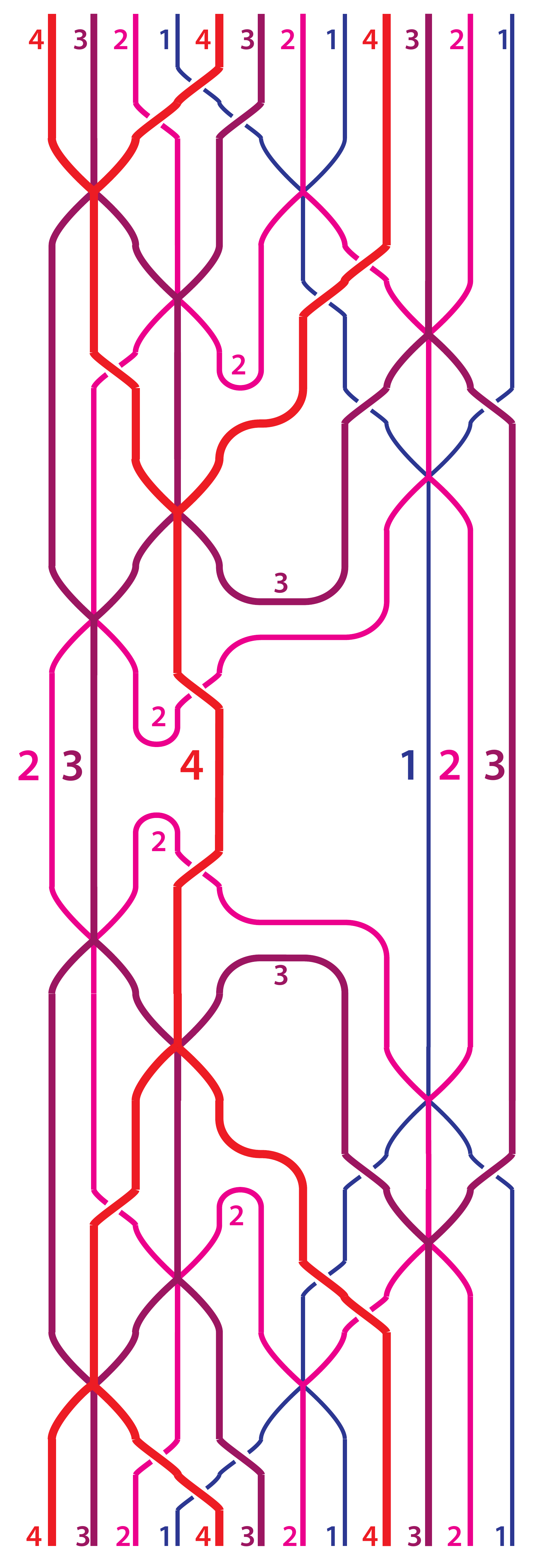}}}
\]

\[
\raisebox{9.25cm}{
{\huge{$\leftrightharpoons$}}}
\scalebox{.18}{\includegraphics{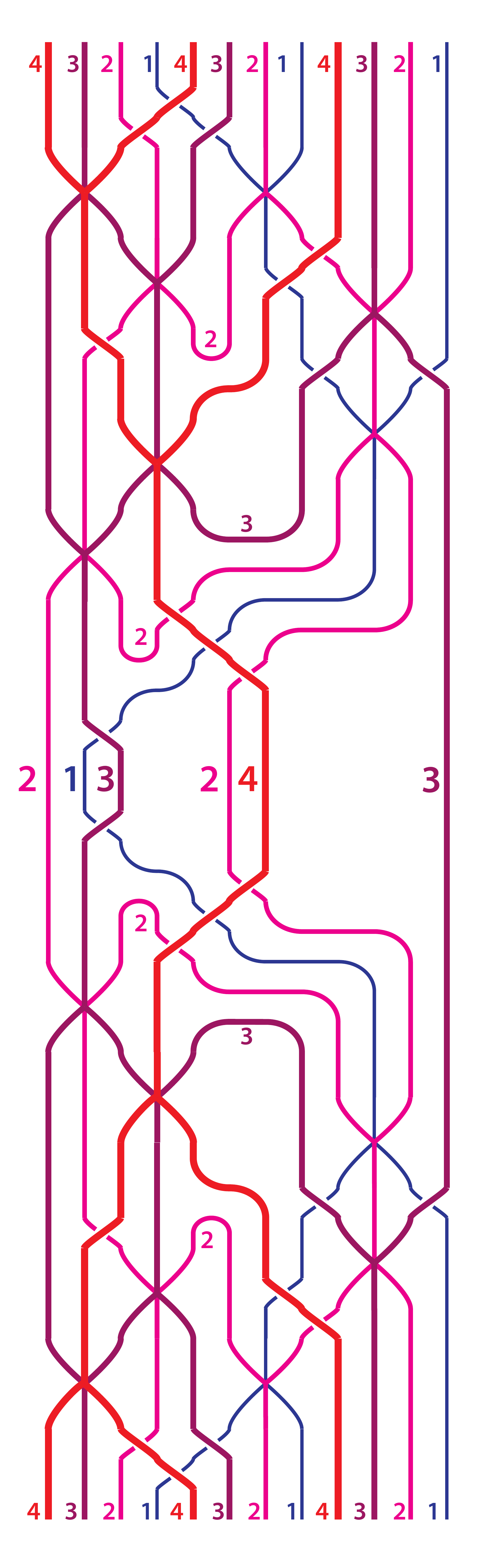}}
\raisebox{9.25cm}{
{\huge{$\leftrightharpoons$}}}
\raisebox{1.cm}
{\scalebox{.15}{\includegraphics{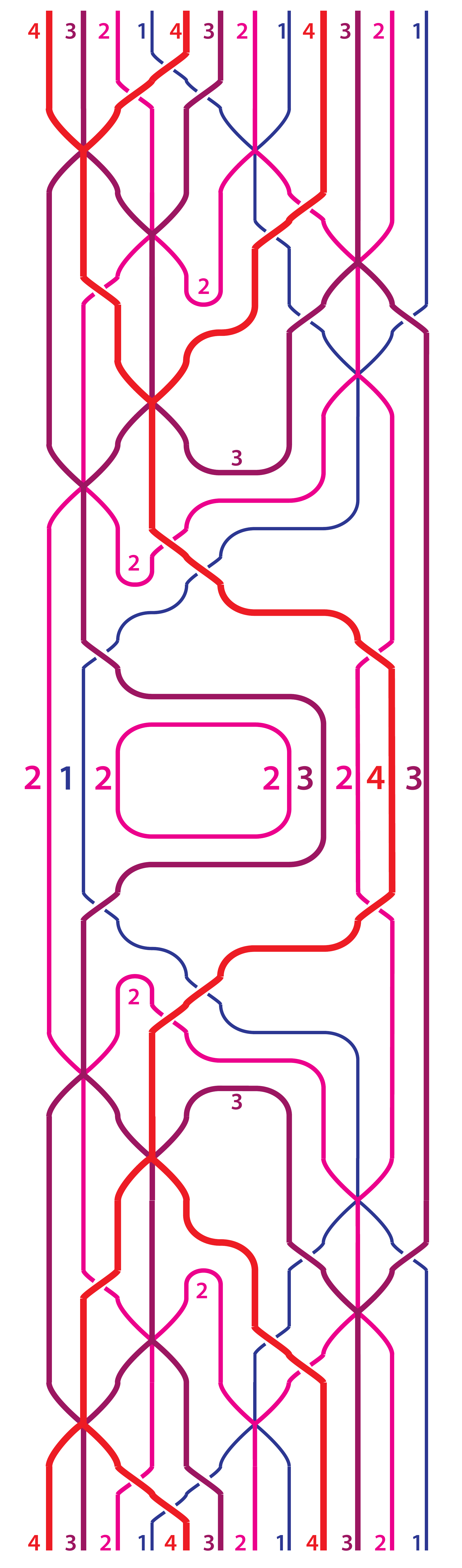}}}
\]

\[
\raisebox{11.25cm}{
{\huge{$\leftrightharpoons$}}}
\raisebox{2.25cm}{
\scalebox{.15}{\includegraphics{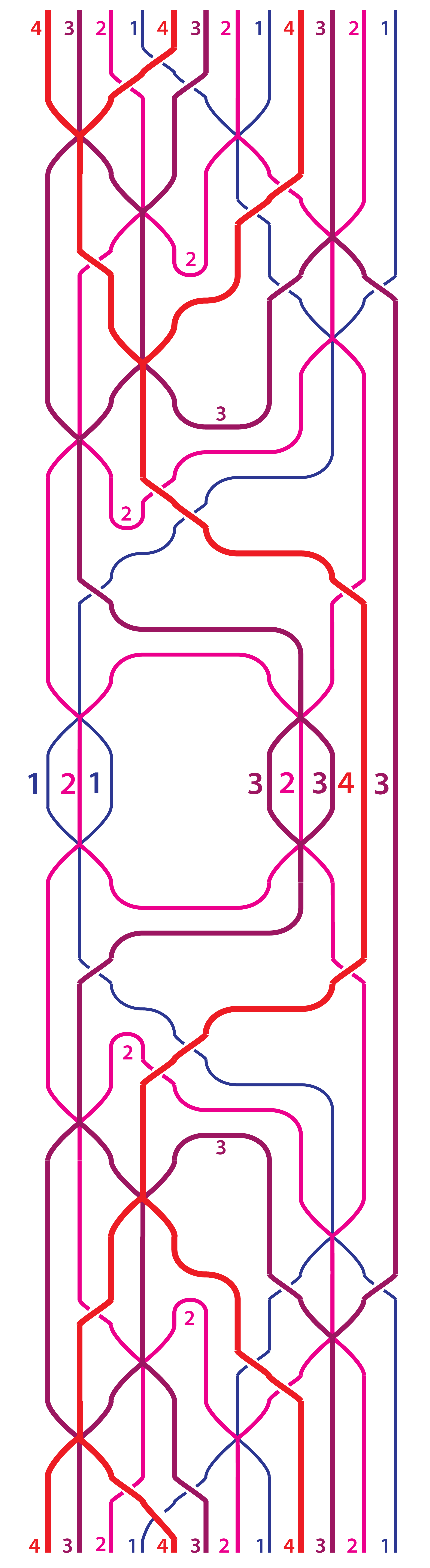}}}
\raisebox{11.25cm}{
{\huge{$\leftrightharpoons$}}}
\raisebox{2.5cm}
{\scalebox{.13}{\includegraphics{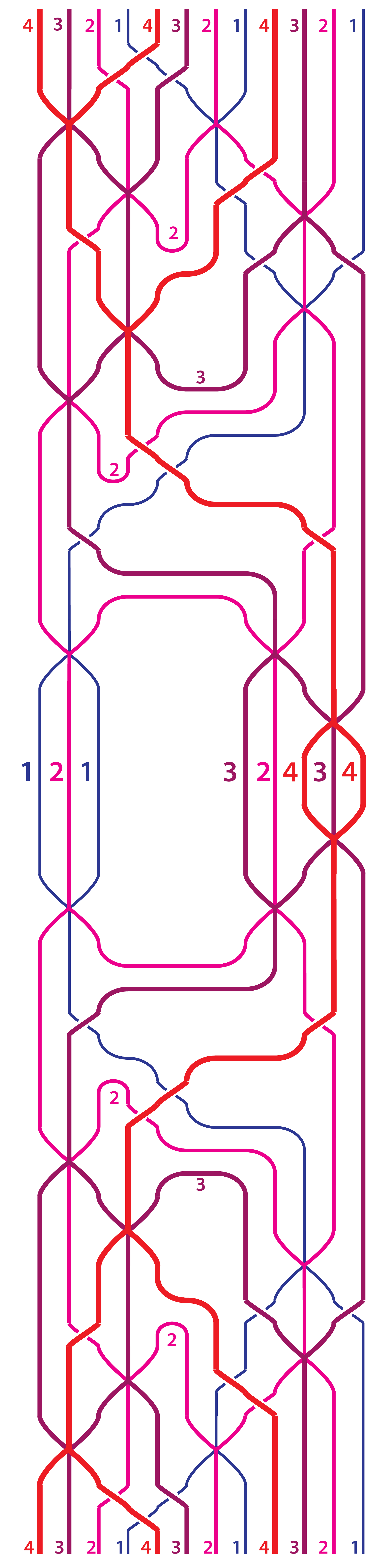}}}
\]

\[
\raisebox{12cm}{
{\huge{$\leftrightharpoons$}}}
\raisebox{1.9cm}{
\scalebox{.13}{\includegraphics{PfL1S14.pdf}}}
\raisebox{12cm}{
{{$\bullet$}}}
\]


\begin{corollary}
\label{cor1}
\[\scalebox{.3}{\includegraphics{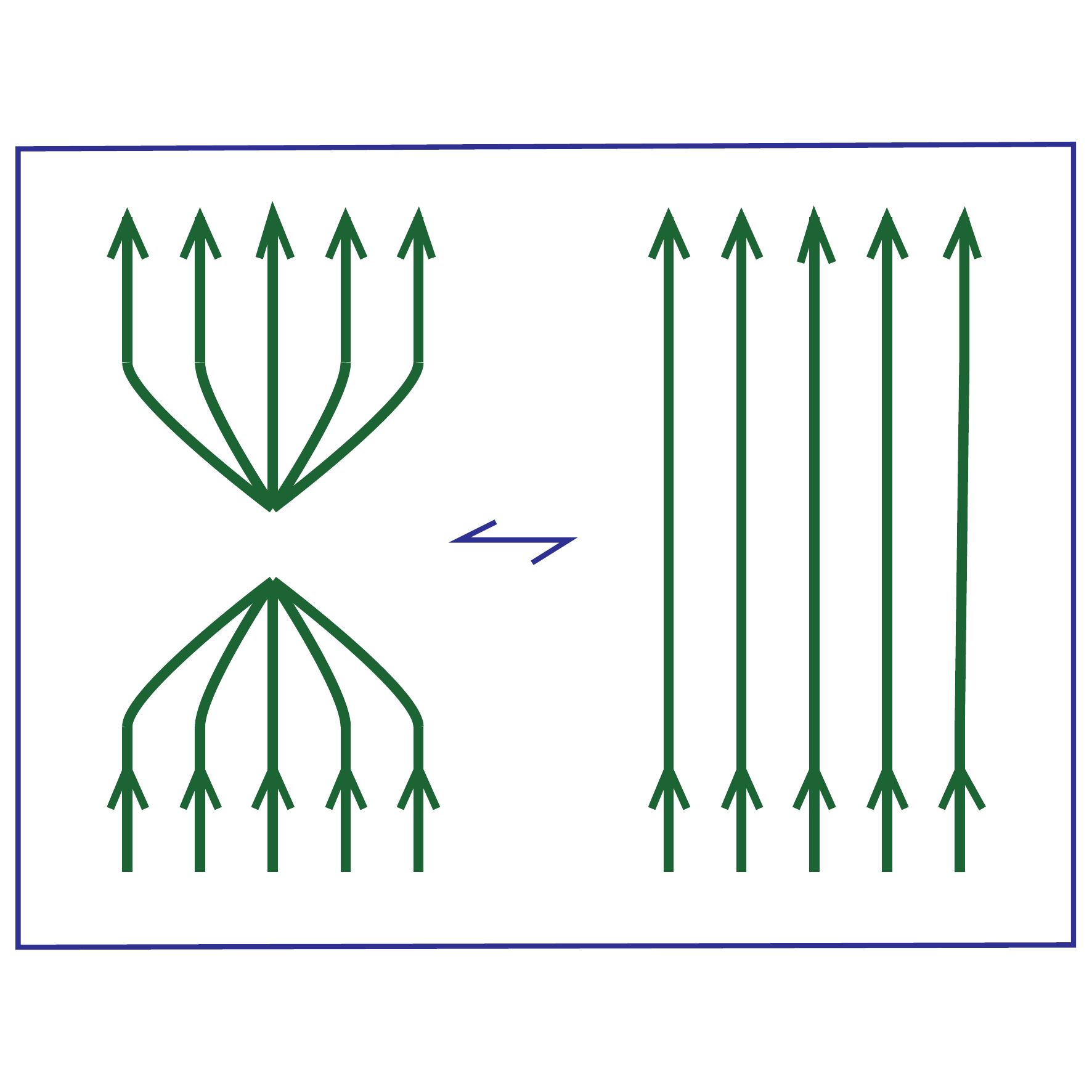}}\]
\end{corollary}
{\sc Proof.} 
\[\scalebox{.1}{\includegraphics{PrCor1.pdf}}\]

\begin{lemma}
\label{MainLemma2}
\[\scalebox{.25}{\includegraphics{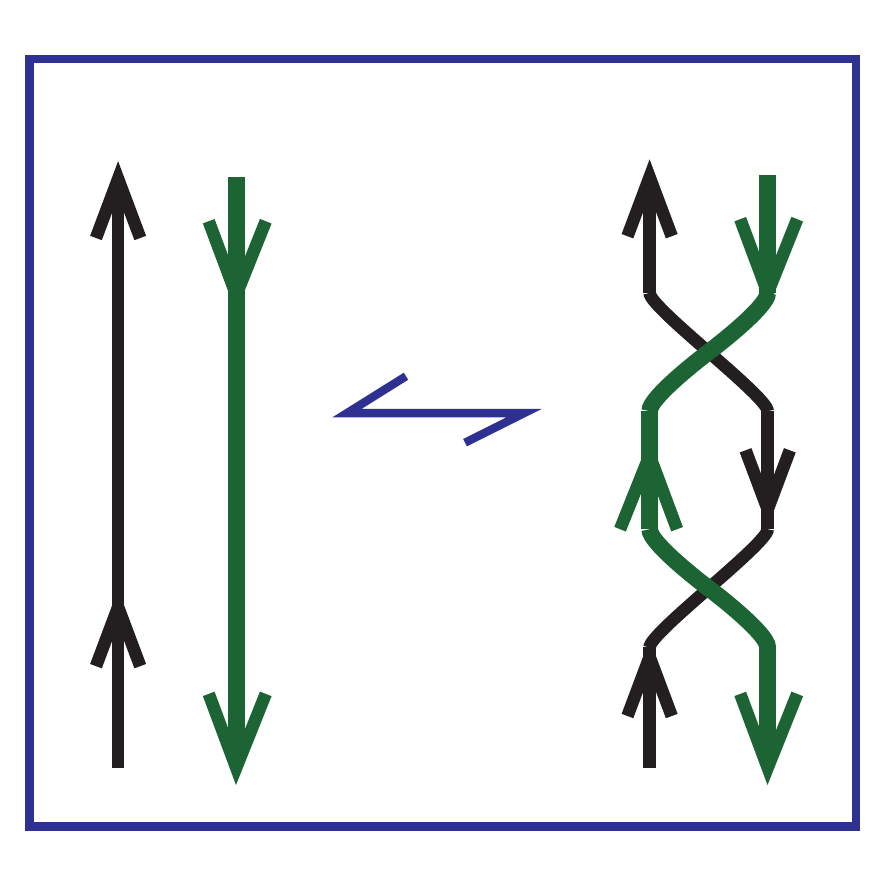}}\]
\end{lemma}
{\sc Proof.} 
\[\scalebox{.15}{\includegraphics{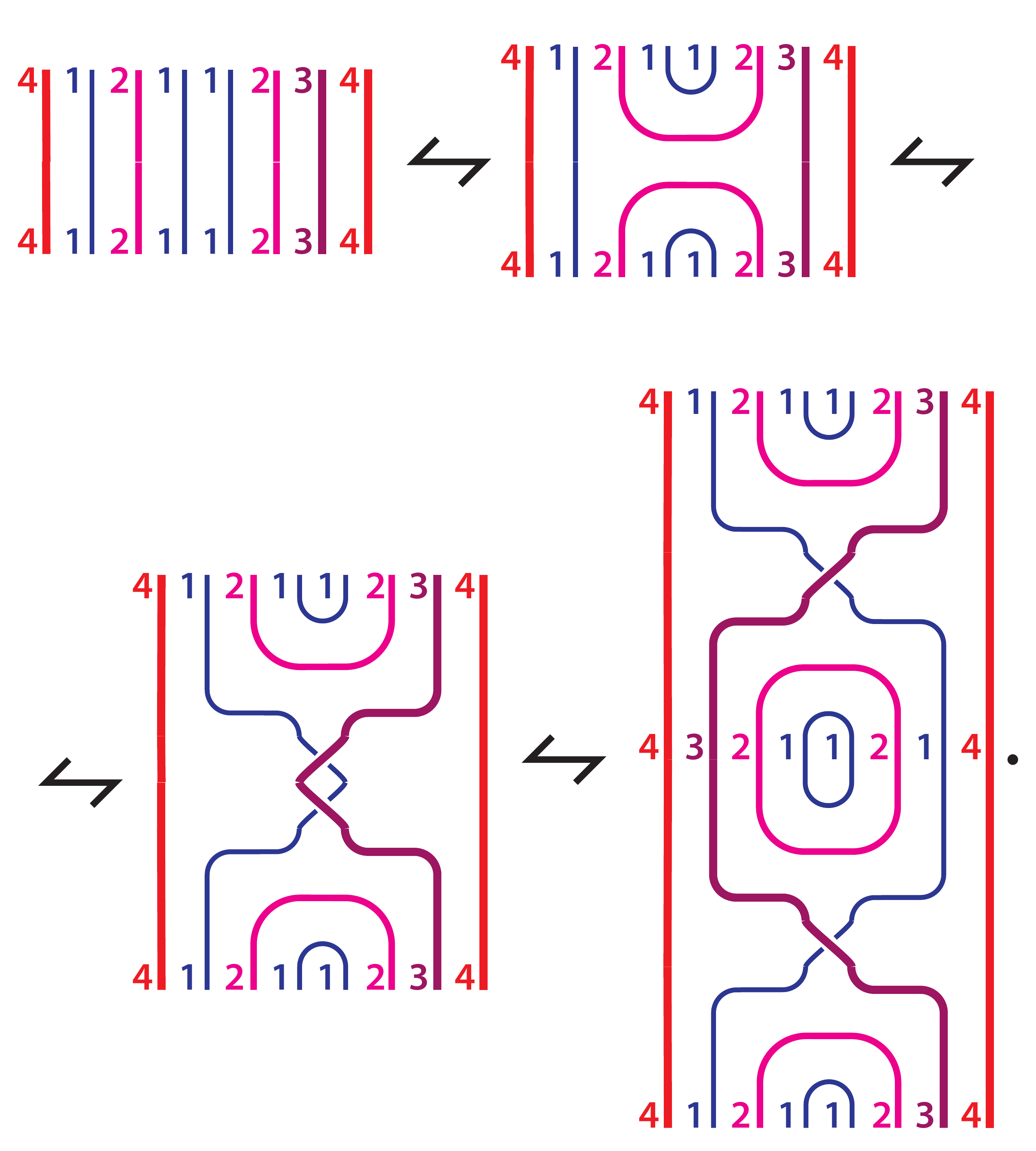}}\]

\noindent
{\bf Remark.}
Let a finite group $G = \langle X: R \rangle$ be given with a permutation representation $\Sigma_n \stackrel{\psi}{\longleftarrow} G$. Consider the category ${\mathcal C}(G)$ whose definition was sketched in Section~\ref{coolcat}. The  morphisms in ${\mathcal C}(G)$ are generated by the relators in $R$ and the trivial relations such as $aa^{-1}=1$. In such presentations, we write the relators as words in $F(X)$. For example, if $G=D_5$, then the relators are $r^2$, $x^5$, and $rxrx$. Thus we write $w\in R$ to indicate that a word $w$ is a relator. 

The relators are graphically depicted as vertices, dual to the $2$-dimensional disks in the classifying space $BG=K(G,1)$. In addition, there are $\cup$ and $\cap$ morphisms that correspond to the trivial relations. The techniques that we demonstrated above in the proofs of Lemmas~\ref {MainLemma1}  and~\ref{MainLemma2} indicate that, in general, the morphisms in ${\mathcal C}(G)$ are invertible.

Here is why. Consider the braid chart moves of type CI in Kamada's book \cite{KamBook}, page 143, Figure 18.17. By removing the orientation arrows upon these, they become permutation chart moves  as in Fig.~\ref{almostKmoves}. These are the moves that generate homotopies of immersed surfaces in $3$-dimensional space. In particular, white vertices, crossings, $\cup$s, and $\cap$s are invertible morphisms in the corresponding category ${\mathcal C}(\Sigma_n)$. The permutation chart moves generate the $3$-cells in the classifying space $B\Sigma_n$. A word $w\in R$ maps, via the group homomorphism $\psi$, to a word $\psi(w)$ in the standard generators $t_1$ through $t_{n-1}$ in $\Sigma_n$. The word $\psi(w)$ is a consequence of the relations that correspond to $\cup$s, $\cap$s, crossings, and white vertices. Such a word, $\psi(w)$ is invertible in ${\mathcal C}(\Sigma_n)$ since $\cup$s, $\cap$s, crossings, and white vertices are invertible in the category associated to $\Sigma_n$.  Thus the invertibility of relators in $R$ is a consequence. 

However, the general proof of invertibility of the morphisms in ${\mathcal C}(G)$  sketched above is only existential. That is to say,  even though we know an algorithm to compute, for example, that $(15432)(15432)(15432)=(13524)$, we explicitly used a string diagram for the corresponding words in the generators $t_1$ through $t_4$, and simplified it using the permutation relations. Then having explicitly written down a simplification, we reflected it graphically, and step-by-step canceled (or added canceling pairs of) the vertices.

{\bf Those step-by-step simplifications of  intertwined disks give an explicit folding of the branched covers.} Before proceeding to the proof of Theorem~\ref{DihedralThm}, we indicate some of the substitutions between the dihedral charts and their associated permutation charts in the steps indicated in Section~\ref{thesteps}. 

One should view the illustrations on the next pages as an all too brief summary of some calculations that are occurring in a diagrammatic realm. In an analogous situation, a string of published algebraic simplifications may skip several steps that were filled in by the authors in their private notes. 
We have carefully done the homotopy moves upon the sphere in $3$-space that is folded as a dihedral cover of $S^2$, but we cannot fit the corresponding illustrations into the margin. For example, even though the line segments that correspond to crossings are spaced two centimeters apart in the illustrations, several drawings are upon canvases that measure roughly $2$m $\times$ $1.5$ meters. And our calculation spans about forty canvases even when multiple steps occur between successive canvases. 

We notice that the simple branched dihedral cover of $S^2$ that has  eight simple branch points (four red and four pink) is also a $2$-sphere.  Consequently, the $D_5$ cover of $S^3$ that is branched along $T(2,5)$ is a $3$-sphere.

\newpage

The  illustration upon the current page corresponds to Step 2 above. It indicates a branched  dihedral cover of $S^2$ with four branch points as it is folded into $S^2 \times [0,6]$. On the top left,  a $1$-handle is attached between the fourth and fifth  sphere (red or $t_4$), and a $1$-handle is attached between the $1$st and $3$rd sphere (pink  or $t_2$). On the bottom the $2$nd and $4$th sphere are connected, as are the $1$st and $5$th.

\[\scalebox{.06}{\includegraphics{Step2inPerm.pdf}}\]

\newpage

The  illustration  upon this page corresponds to Step 12 above.  The four branch points from the center of the previous drawing have wrapped around each other following the swirling pattern in the steps 2 through 12 presented above. In this step and many of the subsequent drawings, the size of the illustration is around $3$ square meters or larger. 
\[\scalebox{.06}{\includegraphics{Step12inPerm.pdf}}\]

\newpage

The  illustration  upon this page roughly corresponds to Step 17 above.  Lemma~\ref{MainLemma2} and its analogues have been applied so that the parallel vertical arcs on the left and right side of the figure correspond to powers of the rotation $x$ in the dihedral group. The folded $S^2$ in $S^2\times [0,6]$ is still homotopic to the first one illustrated.

\[\scalebox{.06}{\includegraphics{Step17inPerm.pdf}}\]

\newpage

The  illustration  upon this page roughly corresponds to Step 19 above. The proof of Corollary~\ref{cor1} has been applied. 

\[\scalebox{.06}{\includegraphics{Step21inPerm.pdf}}\]

\newpage

In the next two illustrations the order of cancelation between steps 19 and 20 and that between steps 20 and 21 have been interchanged. These steps correspond to passing the two minimal  points of the knot $T(2,5)$.

\[\scalebox{.06}{\includegraphics{Step23inPerm.pdf}}\]

\newpage

In the illustration here, the diagram indicates the self-intersection set of five $2$-dimensional spheres that are immersed into $S^2\times [0,6]$. Each is embedded, but they intersect generically. There are twenty triple intersections among these spheres. 

\[\scalebox{.06}{\includegraphics{Step24inPerm.pdf}}\]

\newpage

The $L$ shaped diagram of intersections has been flattened. Look among the triple points and crossing changes and see which can cancel immediately. 

\[\scalebox{.06}{\includegraphics{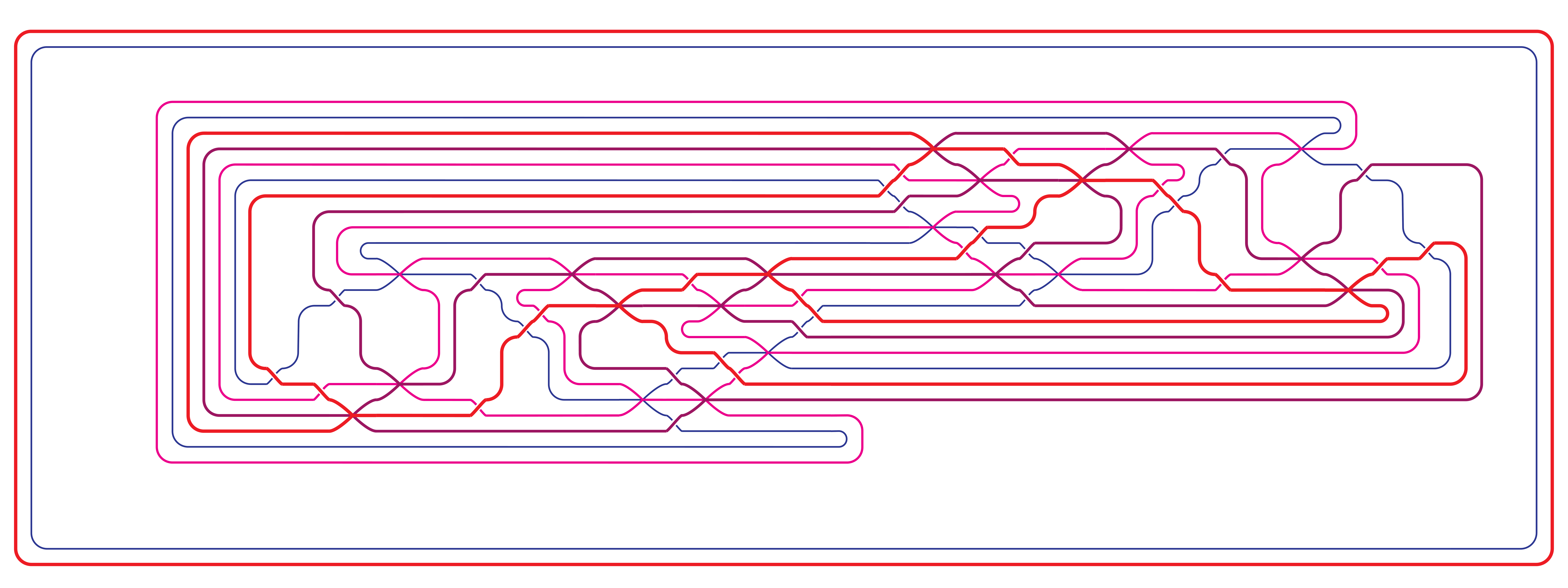}}\]

\newpage

After about five pages of modifications, we were able to move the spheres into the configuration that appears here. The highlighted region indicates an area in which the quadruple point (tetrahedral or Zamolodchikov) move will be used. 

\[\scalebox{.06}{\includegraphics{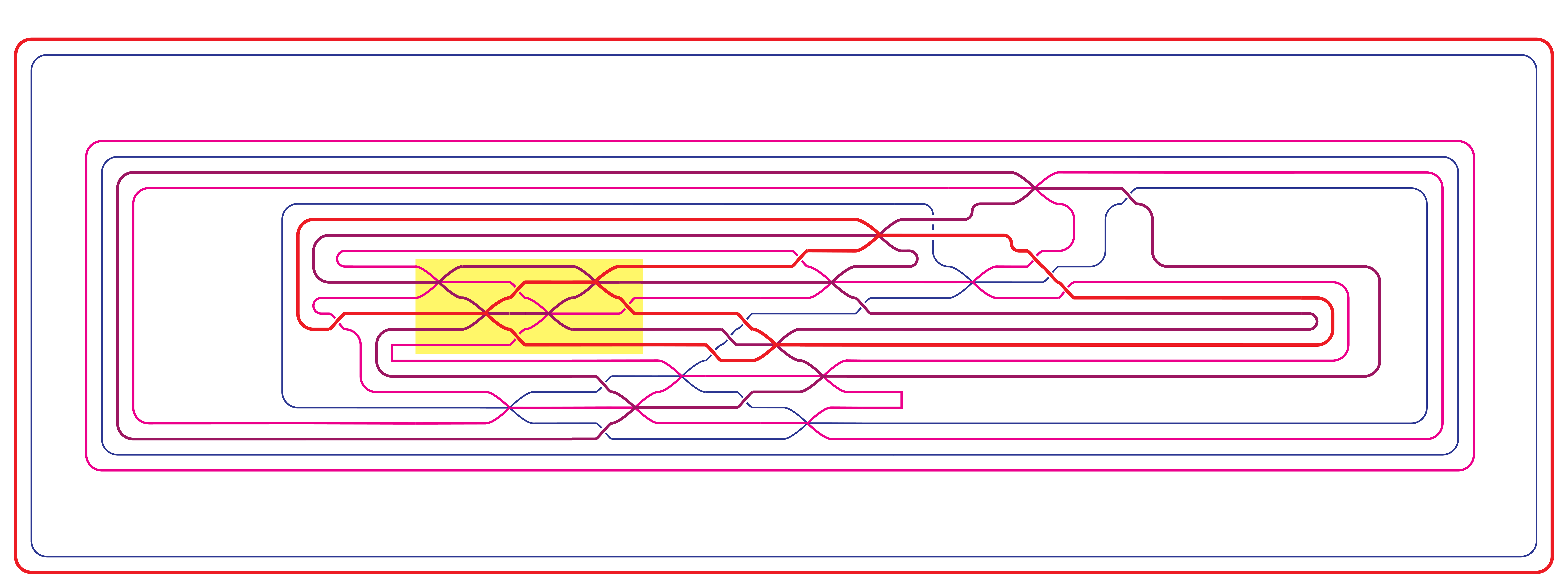}}\]

The highlighted region  here indicates the region in which the quadruple point  move occurred.

\[\scalebox{.06}{\includegraphics{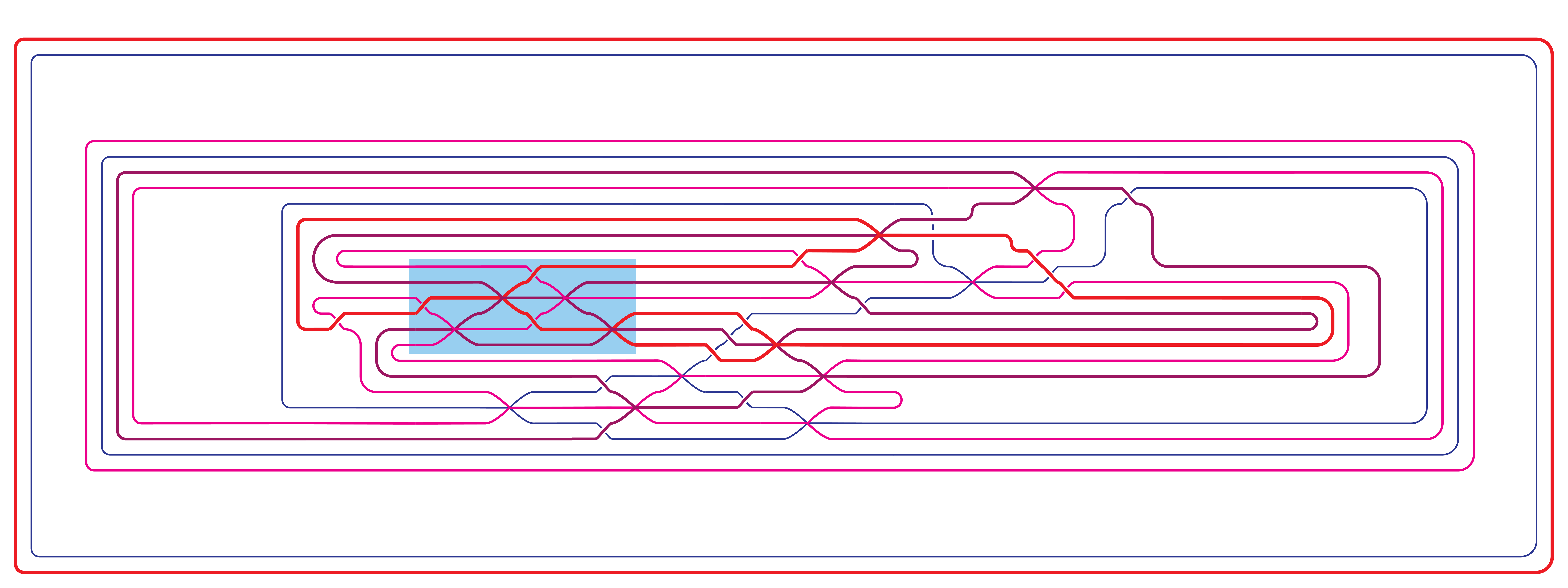}}\]

\newpage

Several more modifications occur, and another quadruple point move is highlighted. 

\[\scalebox{.06}{\includegraphics{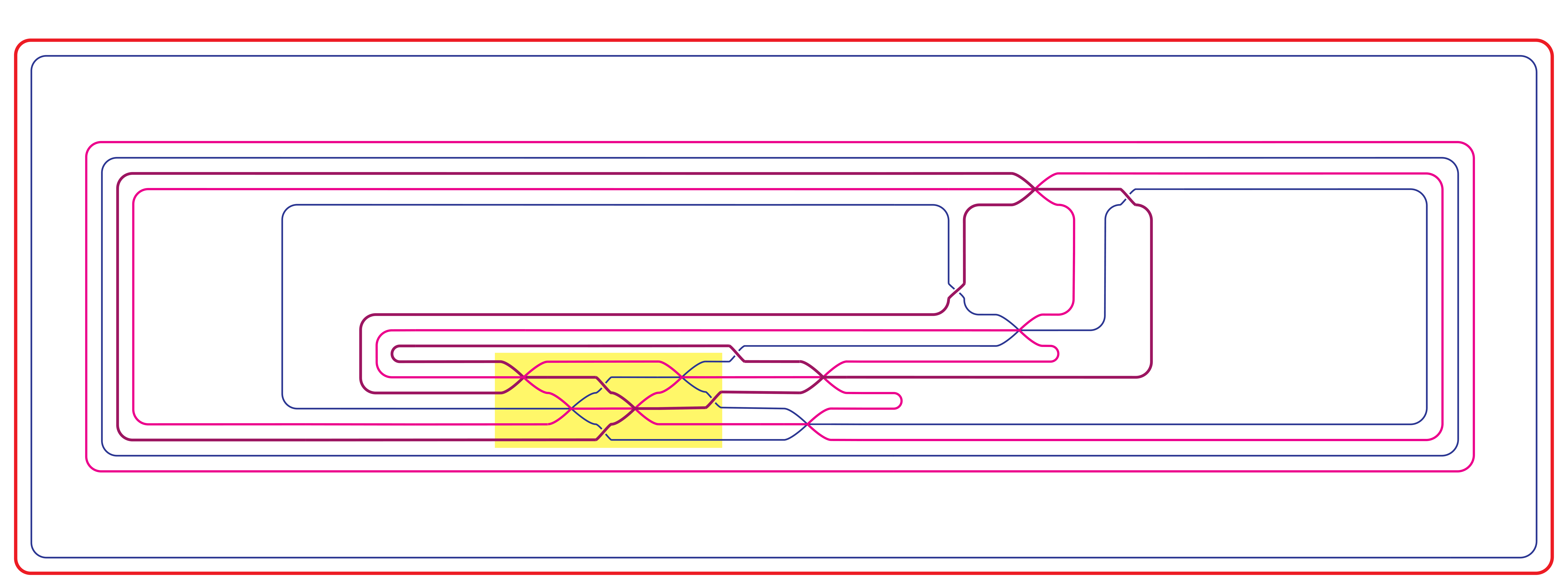}}\]

In the illustration below, it is possible to determine that some pairs of triple points can be canceled, and in a few more moves, a collection of nested loops of intersections between pairs of spheres results. 

\[\scalebox{.06}{\includegraphics{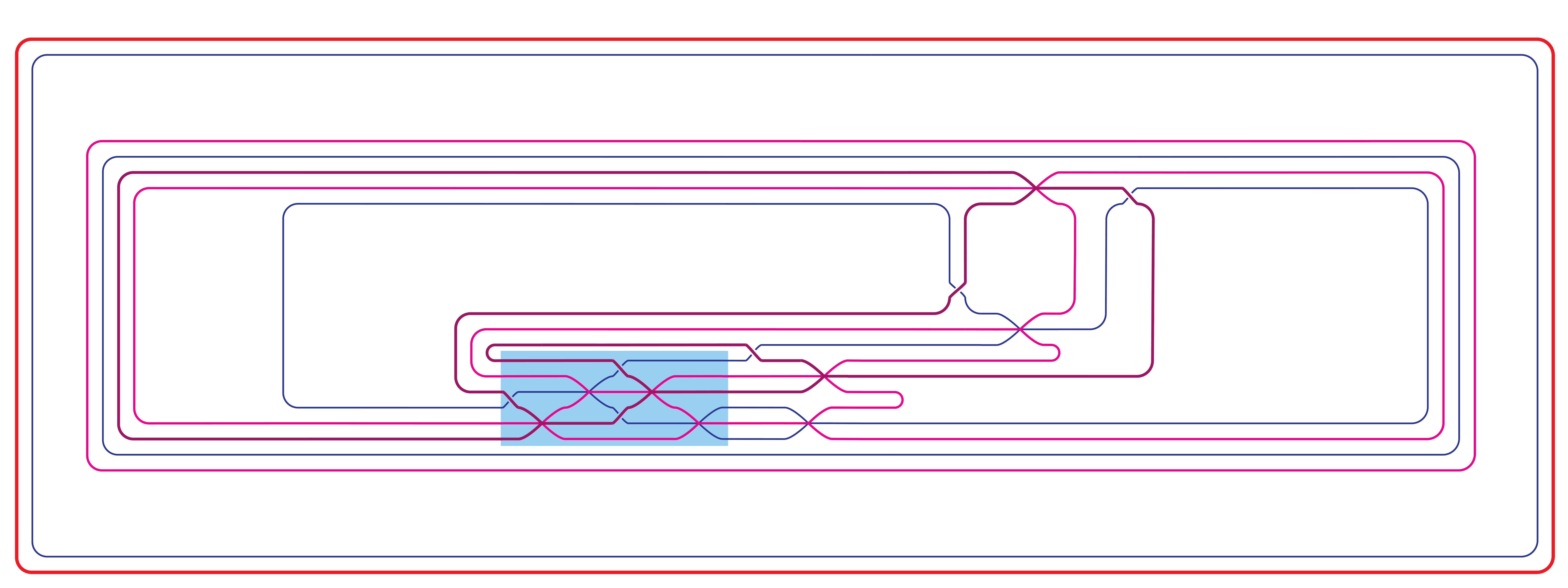}}\]

\[\scalebox{.06}{\includegraphics{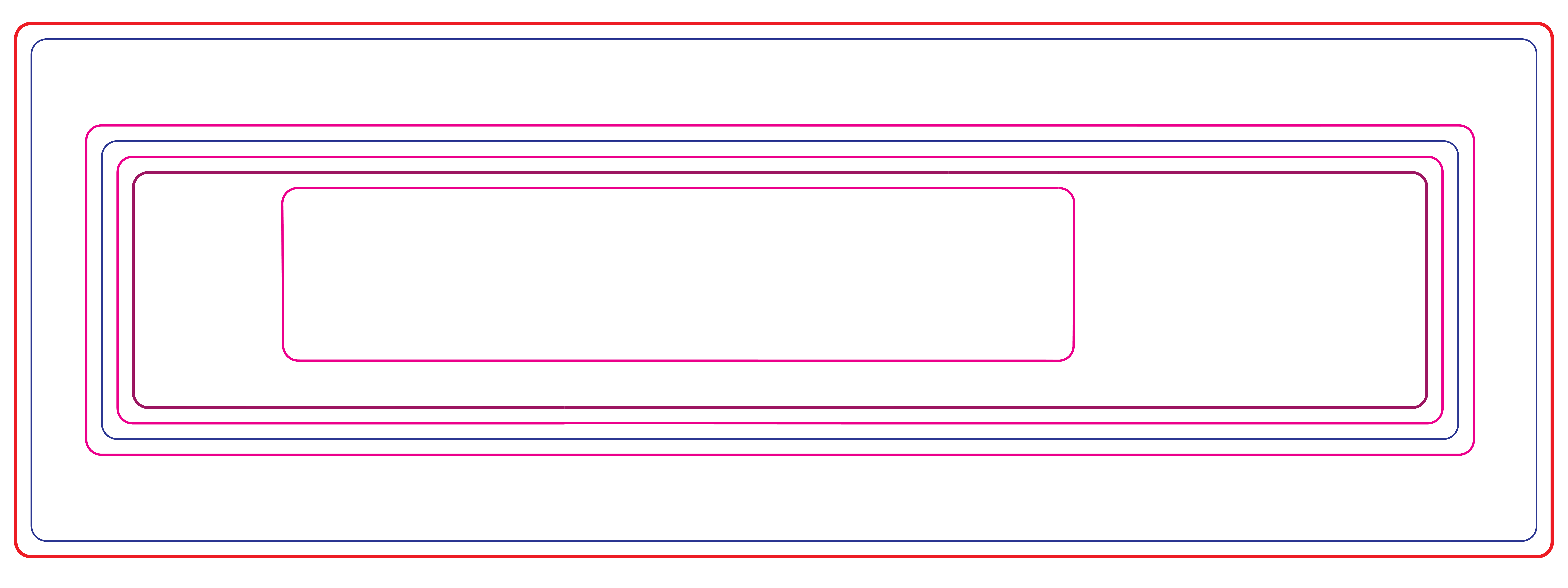}}\]
These nested loops of double points may be removed sequentially from smallest to largest by using the type II bubble move.

\section{Choosing reflections}

Let $D_n$ denote the dihedral group when $1<n=2k+1$ is  positive and odd. 
Then $D_n= \langle r, x: r^2=x^n=(rx)^2=1 \rangle.$ All of the reflections ($\{r x^j \}$) are conjugate, by means of  a power of $x$,  to a fixed reflection. To make sure that the permutation representation that we choose does not get too messy diagrammatically, we choose the reflection
\[r_k= (k-1,k+1)(k-2,k+2)\cdot \cdots \cdot (1, n-2)(n-1, n)\]
as the generating reflection in the group $D_{2k+1}$. So when $n=5$, $r=r_2=(1,3)(4,5)$; when $n=7$, $r=r_3=(2,4)(1,5)(6,7)$; and when $n=9$, $r=r_4=(3,5)(2,6)(1,7)(8,9).$
The transposition $(n-1,n)=t_{n-1}$ is a standard generator of the symmetric group $\Sigma_n$ into which the dihedral group is represented. And the remaining transpositions that are factors of $r$ are all conjugates of $t_{k}=(k,k+1)$.

For example, consider $n=5= 2\cdot 2 +1$. The reflection $r_2=(1,3)(4,5)$ is given. The transposition $(1,3)$ was written above as $(1,3)=t_1 t_2 t_1$. 
When $n=7$, the reflection $r_3=(2,4)(1,5)(6,7)$, and $(2,4)=t_2t_3t_2$ while $(1,5)=t_4t_1(t_2t_3t_2)t_1t_4.$ 

In the general case $n=2k+1$,  the reflection $r_k$ can be written  as a product of the last generator $t_{n-1}$ and a product of conjugates of $(k-1,k+1)= t_{k-1}t_kt_{k-1}.$

For consistency with the $n=5$ case, we choose $x=(1,n,n-1, \ldots, 2)$ to be the generating rotation. 

Suppose that a  $D_n$-cover of $S^2$, called $M^2$,  that has two branch points is given. The exterior $E(K)= S^2 \setminus {\mbox{\rm int}} \left(D^2_{(-)} \cup D^2_{(+)}\right)$ is homeomorphic to an annulus. The ``knot" $K$ is a pair of points, that can be identified with the north and south poles of $S^2$.  A homomorphism $D_n \stackrel{\phi}{\longleftarrow} \pi_1(E(K))=\Z$, then, is determined by sending the generator $1 \in \Z$ to either a reflection $x^j r_k x^{-j}$ or a rotation,   $x^j$. 

In case $1 \mapsto x^j$, there is a folding of  $M$ that is given by means of a permutation chart that is analogous to that in the top of  Fig.~\ref{fromnone2five}. In such a chart, the arcs are labeled from top to bottom \[ [t_{n-1}, t_{n-2},\ldots, t_1], [t_{n-1}, t_{n-2}, \ldots, t_1], \ldots , [t_{n-1}, t_{n-2} \ldots, t_1] \] where the sequence of decreasing permutation generators is of length $j$. 
Indeed, if such a chart has its two non-simple branch points on the left and the right of the figure, then the permutation chart can be oriented by arrows pointing left-to-right, and an embedded braiding
$S^2 \times [0,n+1] \times [-1,1] \stackrel{f}{\longleftarrow} M$ is achieved by lifting the permutation generators $t_\ell$ to braid generators $\sigma_\ell$.

It is worth remarking that the long word $[t_{n-1}, t_{n-2},\ldots, t_1]^j$ can be reduced into something more manageable. At least, the exponent can be reduced to lie within the interval $[-k,k]$. After that reduction, one looks for reduced expressions that represent  the cycles $(1,n, n-1, \ldots, 2)^{\pm j}.$

\begin{figure}[htb]
\[
\scalebox{.25}{\includegraphics{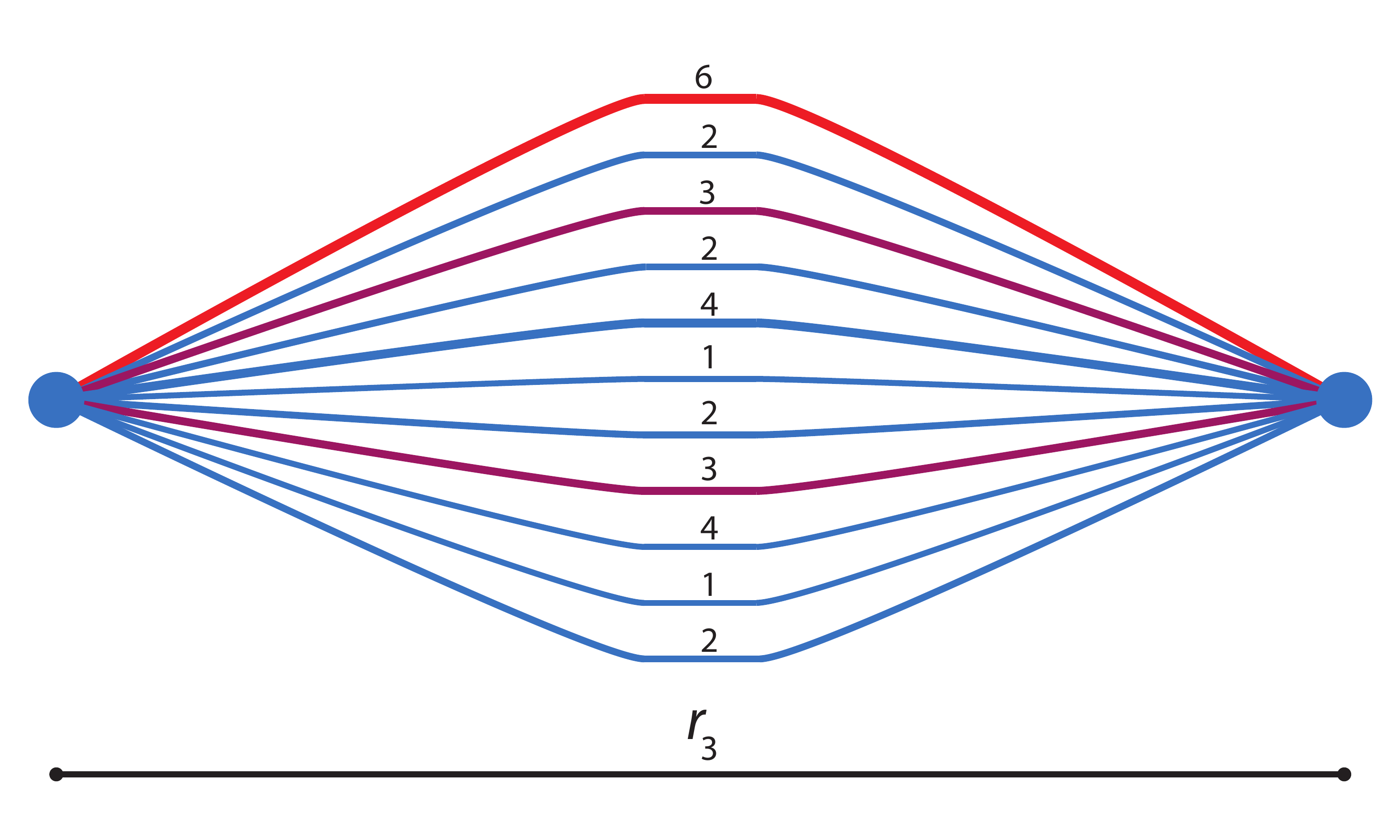}}
\]
\caption{The branched cover of $S^2$ that is associated to the homomorphism $D_7 \stackrel{\phi}{\longleftarrow} \pi_1(S^1 \times [0,1])$ which sends the generator to $r_3$}
\label{squashedovalnest7}
\end{figure}

In case $1 \mapsto x^j r_k x^{-j}$, the permutation chart depicted in Fig.~\ref{squashedovalnest7} is surrounded by $j \cdot (n-1)$ ovals that are labeled (again) sequentially by  
\[ [t_{n-1}, t_{n-2},\ldots, t_1], [t_{n-1}, t_{n-2}, \ldots, t_1], \ldots , [t_{n-1}, t_{n-2} \ldots, t_1]. \] Or more economically, reduced expressions as outlined in the previous paragraph. Orientations can be given for the ovals, and left-to-right orientations are imposed upon the reflection $r_k$, so that again there is an embedded braiding
$S^2 \times [0,n+1] \times [-1,1] \stackrel{f}{\longleftarrow} M.$

 \begin{figure}[htb]
\[
\scalebox{.3}{\includegraphics{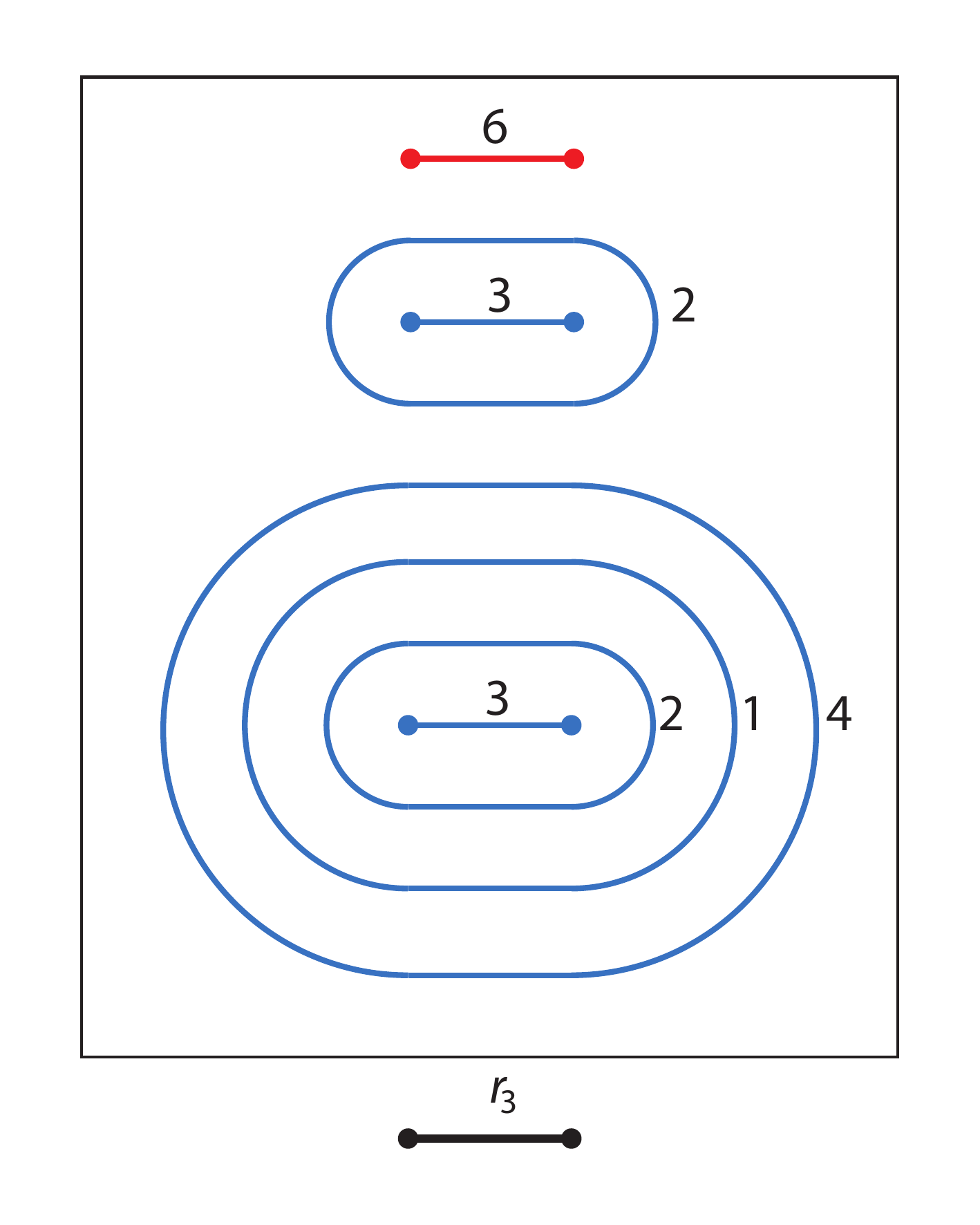}}
\]
\caption{An oval nest depiction of a simple permutation chart associated to the reflection $r_3 \in D_7$}
\label{ovalnest7}
\end{figure}

 The permutation chart that is associated to 
the  homomorphism of $\pi_1(S^1 \times [0,1])$ in which the generator maps to the reflection $r_3 \in D_7$ is depicted in Fig.~\ref{squashedovalnest7}.  Each arc in this figure represents a double curve of a singular surface that is mapped into $S^2 \times [0,8]$. The index upon each of these curves demonstrates the proximity of the intersection curve with the sphere $S^2 \times \{8\}$ upon which the branch cover projects. Larger integers indicate double curves that are closer. 
In Fig.~\ref{ovalnest7}, a  perturbation of the chart  from Fig.~\ref{squashedovalnest7} is illustrated. 
 The branch points have been perturbed into being simple branch points.

In this way, we have embedded braidings for any dihedral branched cover of $S^2$ in which the cover has two branch points. Also the process of going from Fig.~\ref{squashedovalnest7} to Fig.~\ref{ovalnest7} represents a de-singularization. So the embedded braiding of such $M$s can be perturbed into simple surface braids.

\section{Proof of Theorem~\ref{Dihedral2}}
\label{mequals2}
\begin{figure}[htb]
\[
\scalebox{.25}{\includegraphics{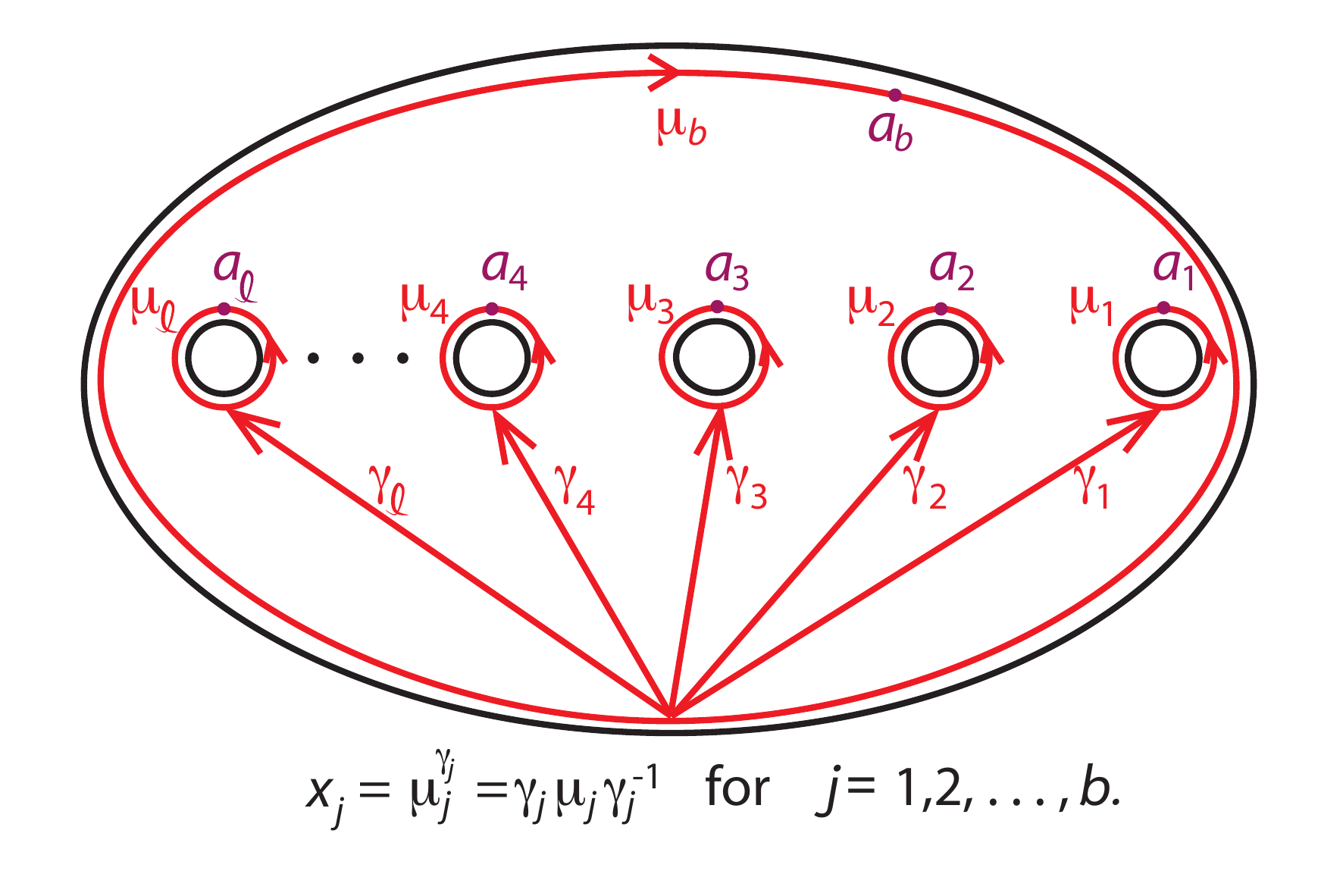}}
\]
\caption{A choice of generators $x_j$ for the free group $F_\ell$ and their images $\phi(x_j)=a_j \in D_n$}
\label{freegroup}
\end{figure}

Recall the hypotheses of the theorem.
A branched cover  of $S^2$ is given that has a finite number of branch points $K=\{ p_1, p_2, \ldots, p_b \}$. By assumption, $1<b$. Let $M$ denote the covering that is associated to the group homomorphism $D_n \stackrel{\phi}{\longleftarrow} \pi_1(E(K)).$ The compact planar surface $E(K)$ is depicted in Fig.~\ref{freegroup}. It has $b$ boundary components. So its fundamental group is the free group $F_\ell$ of rank $\ell=b-1$.

We begin by describing the unbranched cover over $E(K)$. 
In Fig. ~\ref{freegroup}, generators $x_j = \mu_j^{\gamma_j}= \gamma_j \mu_j \gamma_j^{-1}$ ($j= 1, \ldots, \ell$) for $\pi_1(E(K))$ are chosen to be
 the  based paths formed from the counterclockwise meridional loops corresponding to the ``interior" boundary components of the planar surface $E(K)$. The homomorphism $\phi$ is specified by $\phi(x_j)=a_j$. In addition, the element $x_b=\mu_b$ satisfies the relation $\prod_{j=1}^b x_j=1$ since the loop $\mu_b^{-1}$ is homotopic (via the surface depicted) to $\mu_1^{\gamma_1} \mu_2^{\gamma_2} \cdot \cdots  \cdot \mu_\ell^{\gamma_\ell}.$  Consequently, the element $a_1 a_2 \cdots a_b$ is trivial in $D_n$.

 \begin{figure}[htb]
\[
\scalebox{.25}{\includegraphics{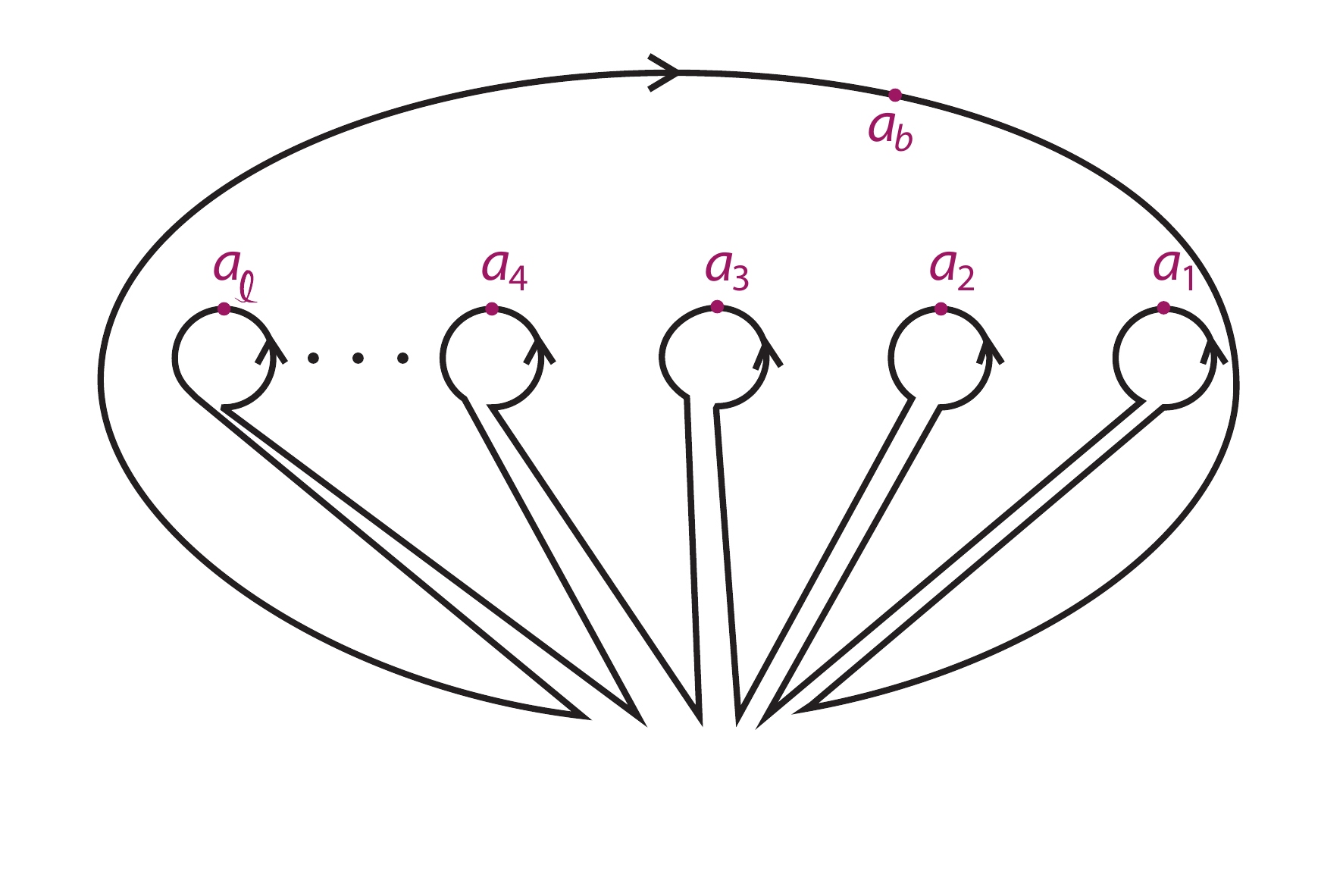}}
\]
\caption{A null homotopy for the product $a_1 a_2 \cdots a_b=1$}
\label{sawtooth}
\end{figure}
 
 The saw-toothed disk that is depicted in Fig.~\ref{sawtooth} indicates the trivialization of $a_1 a_2 \cdots a_b$. The boundary, $B$,  of this disk can be thought of as a  loop that maps into   $BD_n$. The interior meridians map to the corresponding loops in the classifying space.  For definiteness pretend that the edges of the teeth map to the base point of $BD_n$. The loop $B$  extends to the interior of the saw-toothed disk, and that disk maps to the $2$-skeleton of the classifying space. By pulling back the disks that correspond to the relators $x^n$, $r^2$, and $xrxr$,  a dihedral chart without black vertices  is obtained in $E(K)$.
 From such a chart, an immersion of the unbranched cover of $E(K)$ that is associated with $\phi$ is to be constructed.

\begin{figure}[htb]
\[
\scalebox{.1}{\includegraphics{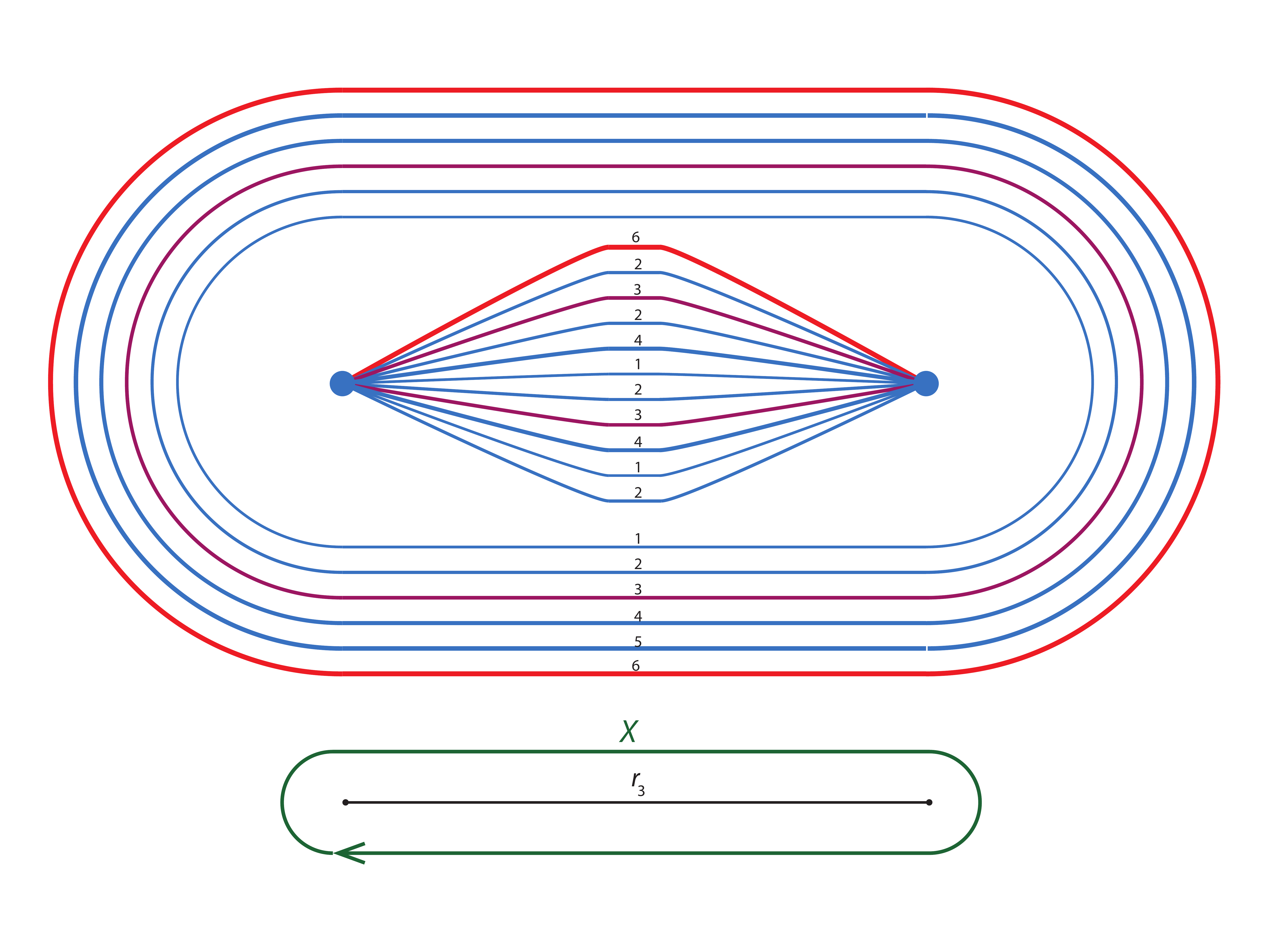}}
\]
\caption{Conjugating the reflection $r_3$ by the rotation $x$}
\label{conjrefl}
\end{figure}

The element $a_j$ is either a reflection or a rotation in $D_n$.  The reflections are all conjugate to $r_k$ where $n=2k+1$. The rotations are all powers of $x$. Thus the reflections are of the form $x^j r_k x^{-j}$, and the rotations are of the form $x^j$. Represent $D_n$ into the permutation group $\Sigma_n$ via $x\mapsto (1,n,n-1, \ldots, 2)$ and 
$r_k \mapsto (k-1,k+1)(k-2,k+2) \cdots  (1, n-2)(n-1, n)$. Now decorate the boundary of $E(K)$ nearby the meridian $\mu_j$ as follows. If $a_j$ is a reflection, then write, along the meridian $\mu_j$, the word in the generators  $t_i$ that represents  the permutation $(k-1,k+1)(k-2,k+2) \cdots  (1, n-2)(n-1, n)$ or its inverse. See for example Fig.~\ref{squashedovalnest7}. To achieve the conjugation by $x^{i_j}$ write that conjugate in terms of the generators $t_i$ along $\mu_j$ and surround the syllable $r_k$  with  concentric bowls that express the cancelations  of the conjugations by $(t_1 t_2 \cdots t_{n-1})$ as indicated in Fig.~\ref{conjrefl}. 

If $a_j$ is a rotation, then decorate the boundary of $E(K)$ nearby $\mu_j$ by the reduced word that represents a power of the product $(t_1 t_2 \cdots t_{n-1})$. The decoration is an explicit  word in the generators, $t_i$,  of $\Sigma_n$ along segments of the boundary that represents this cycle. The word corresponds to the  permutation representation of the cycle that is a power of $x$. 

In the manner described by the prior two paragraphs, a word  built from the  $t_i$,  that represents the identity element in $\Sigma_n$, appears along the boundary of the saw-toothed disk.  Therefore, a permutation chart without black vertices can be found in the interior of the disk as a trivialization chart for the word that represents a trivial element in $\Sigma_n$. Specifically, analogues of  Corollary~\ref{cor1} and Lemma~\ref{MainLemma2}, hold for the permutation charts. So the relator vertices in the dihedral charts can be replaced, as in the proofs of the lemmas, with crossings, white vertices, and, when  needed, by inserting type II bubble and saddle moves. Such a  permutation chart  without black vertices describes a folding of the cover of $E(K)$ that is immersed into $E(K) \times [0,n+1]$.

In order to branch along $K$, we add  black vertices. The black vertices are either sources or sinks with labels $x^j$ or $r_k$. The conjugates of $r_k$ by powers of $x$ are achieved by surrounding the corresponding black vertex by  semi-circular bowls as imagined in Fig.~\ref{conjrefl}. The last step in the construction of the folding is to 
cone the  permutation words along the meridians to  being permutation branch points. These will either be of the form as in Fig.~\ref{ex5} or the left of Fig.~\ref{resolve}; see also Fig.~\ref{squashedovalnest7}. In this way, the branched covering of $S^2$ that is associated to the given homomorphism $\phi$ is folded into $S^2\times [0,n+1]$. See also Fig.~\ref{catproduct} that is illustrated below. 

Finally, the folding can be lifted to an immersed braiding in $S^2\times [0,n+1]\times [-1,1]$, by attempting to orient all the edges in the permutation chart. Orientations have to pass consistently through white vertices (triple points), so that the permutation chart lifts to a braid chart (that has non-simple black vertices). In case no such lifting exists, nodes can be inserted at which orientations of the edges of the chart change. Such nodes explicitly represent double points of the branched surface as it lifts into the $4$-dimensional space $S^2\times [0,n+1]\times [-1,1]$.  This completes the proof of Theorem~\ref{Dihedral2}.

 \begin{figure}[htb]
\[
\scalebox{.12}{\includegraphics{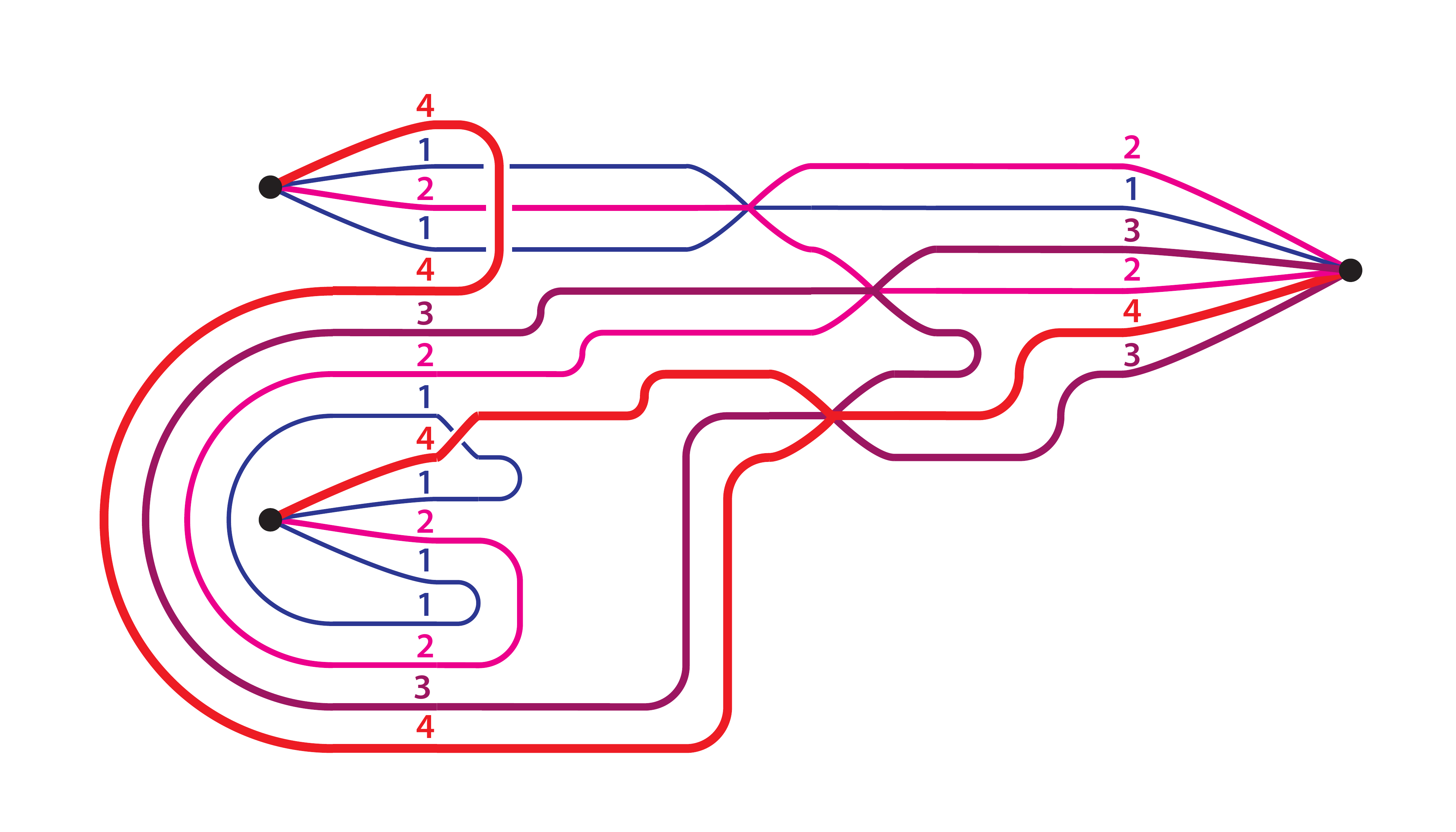}}
\]
\caption{A branched cover of $S^2$ that has three non-simple branch points}
\label{catproduct}
\end{figure}

\noindent
{\bf Example.} In the dihedral group $D_5$, consider the product 
$r(xrx^{-1})x^2= x^{-1}r^2x^{-1}x^2=1.$ Using the representation into $\Sigma_5$ in which $r=(13)(45)$ and $x=(15432)$, a folding of the branched cover of $S^2$ that has three non-simple branch points is constructed, as in the proof of Theorem~\ref{Dihedral2}. The branch point on the top left of Fig.~\ref{catproduct} has a branch point that is colored by $r=r_2$. The bottom left branch point is also colored by $r$, but it has been conjugated by $x$. 

The product of generating transpositions that emerges from the left edge of the illustration is 
\[(45)(12)(23)(12)(45)(34)(23)(12)(45)(12)(23)(12)(12)(23)(34)(45)\]
\[=(12)(23)(12)(45)(45)(34)(23)(45)(12)(12)(34)(45)\]
\[=(12)(23)(12)(34)(23)(45)(34)(45)\]
\[=(23)(12)(23)(34)(23)(34)(45)(34)\]
\[=(23)(12)(34)(23)(34)(34)(45)(34)\]
\[=(23)(12)(34)(23)(45)(34).\]
And notice that $x^{-2}=(12345)^2=(13524)$ which factors as above. Moreover, the gist of the algebraic calculation is illustrated via  Fig.~\ref{catproduct}. The branched cover that is illustrated corresponds to the homomorphism from the free group $F_2$ into $D_5$ in which one generator is mapped to $r$ and the other is mapped to $xrx^{-1}$.  The folding depicted in Fig.~\ref{catproduct} illustrates the method of the proof.

\section{Proof of Theorem~\ref{DihedralThm}}
\label{ProofofDihedral}

\begin{figure}[htb]
\[
\scalebox{.5}{\includegraphics{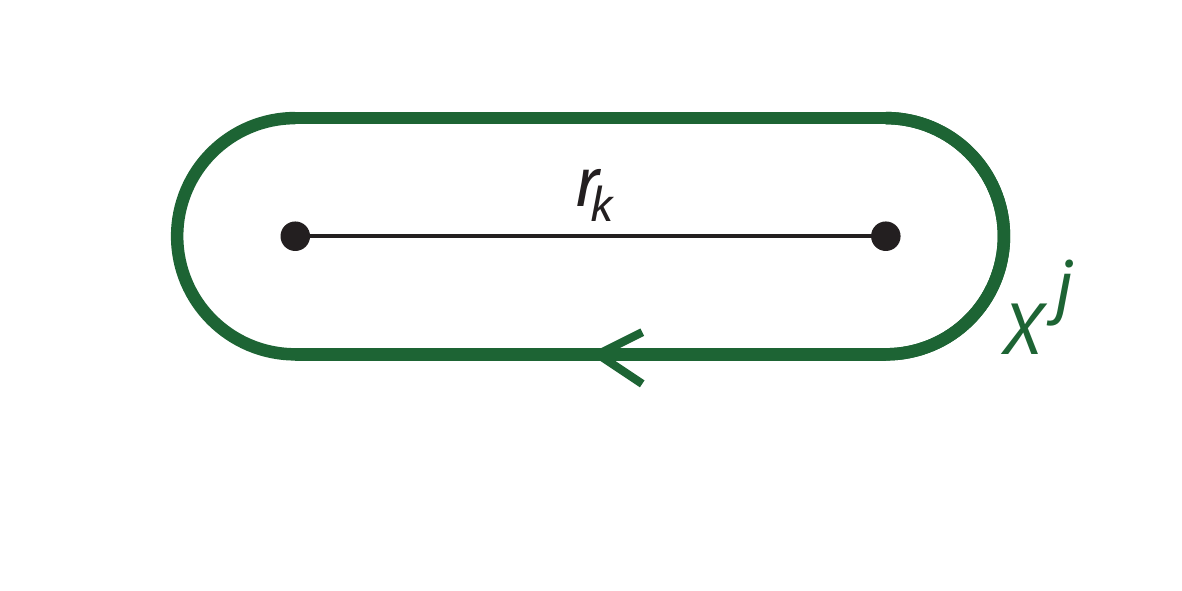}}
\]
\caption{A dihedral cover of $S^2$ that corresponds to the assignment $\mu \mapsto x^{-j}r_k x^{j}$}
\label{singleoval}
\end{figure}

\noindent
{\sc Proof.} Suppose that a knot (or link) $K$ is given. Consider a diagram, $D(K)$, that consists of a collection of arcs in the plane. Color the diagram $D(K)$ 
 in such a way that a reflection in $D_n$ is assigned to each arc.  Such a Fox coloring corresponds to a homomorphism 
$D_n \stackrel{\phi}{\longleftarrow} \pi_1 ( E(K))$ of the fundamental group of the knot's exterior to the dihedral group. A given arc in the diagram of $K$ is colored as depicted in Fig.~\ref{singleoval}. Differing arcs will usually have differing degrees of conjugation, $j$.

The case in which all arcs are assigned the same rotation is another type of homomorphism that is handled in the proof of  Theorem~\ref{CyclicThm} below. More generally, it can be demonstrated that any homomorphism $\phi$ either assigns reflections to all the arcs of a diagram of $K$, or assigns the same rotation to all the arcs. 

\begin{figure}[htb]
\[
\scalebox{.15}{\includegraphics{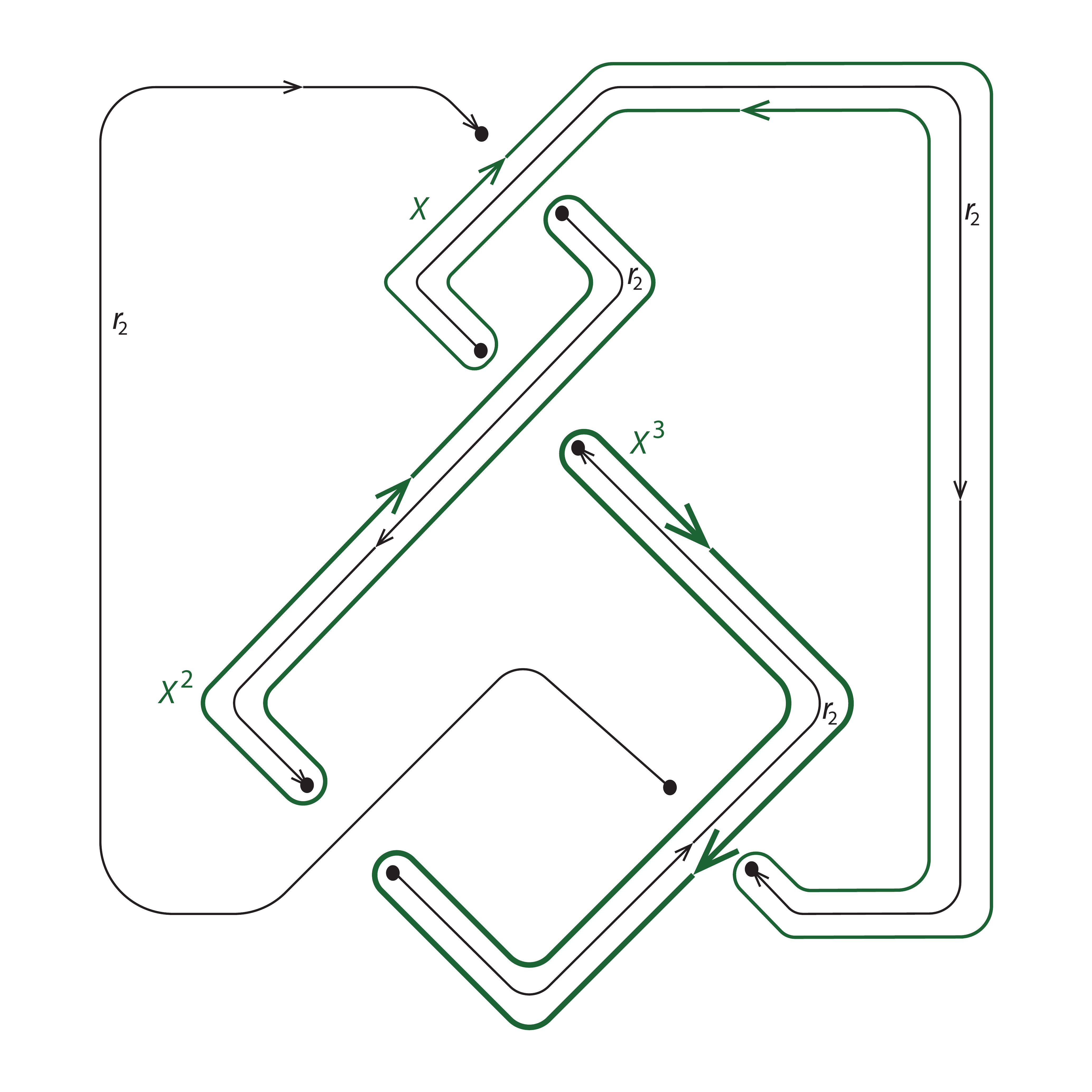}}
\]
\caption{The dihedral chart $C(4_1)_1$ that corresponds to a dihedral coloring of the Figure $8$ knot $4_1$ }
\label{D5coloredFig8}
\end{figure}

The diagram $D(K)$ corresponds to a height function  in which the over-arcs (the arcs that are drawn) represent maxima and the under-arcs that are deleted at crossings represent minima. According to Fig.~\ref{singleoval}, the reflection $r_k$ is assigned to each over-arc, and  the arcs are encircled by differing degrees of rotation, $j$. The diagram $D(K)$ then is replaced with a dihedral chart, $C(K)_1$. The dihedral chart represents a branched cover of $S^2$ that has twice as many branch points as the crossing number of the diagram. Such a chart $C(4_1)_1$ is depicted for the figure-$8$ knot, $4_1$, in Fig.~\ref{D5coloredFig8}.

The essence of the proof is to create a sequence of charts $C(K)_i$ that are (1) related to each other by chart moves, (2) have the properties that 
the branch points remain in the same position throughout the sequence, and (3) result in charts for which the under-arcs are colored by $r_k$. The index $i$ that parameterizes this sequence is undetermined at the outset. 

Perhaps it is a good time to talk about the previous chart $C(K)_0$ and its predecessor $C(K)_{-1}$. The chart $C(K)_0$ represents a trivial covering of $S^2$. It is obtained from $C(K)_1$ by removing all the arcs that are labeled by $r_k$. When we conceive of $S^3$ in a polar/tropical decomposition,  $B^3_N \cup (S^2 \times I) \cup B^3_S$,  the chart $C(K)_1$ appears at the top  of the tropical region that contains the knot. The chart $C(K)_0$ consists of a collection of simple closed curves that are labeled by powers of the reflection $x$ and that are like auras of some of the arcs in the diagram $D(K)$ of the knot $K$. 

The chart $C(K)_{-1}$ is empty, and represents $n$ nested spheres that trivially cover the sphere that contains the empty chart. As one travels from the north pole in $S^3$, these $2$-spheres are successively born. In moving from  the chart $C(K)_0$ to  the chart $C(K)_1$ arcs that are colored by the reflection $r_k$ are born. These correspond to the critical points of $K$ as it passes to the plane of the diagram.

The reason the narrative followed that digression is to explain that after the branch points are reconnected so that the under-arcs are colored by $r_k$, these arcs are no longer surrounded by powers of the rotation, and they can die happily in the southern hemisphere of $S^3$. The chart that remains represents a trivial covering of $S^2$. As in the case of the dihedral cover of the torus knot $T(2,5)$, this  chart without black vertices will be represented into a permutation chart, and  we will argue that it can subsequently be simplified by means of the analogues of the Roseman moves that do not involve branch points.


 \begin{figure}[htb]
\[
\scalebox{.25}{\includegraphics{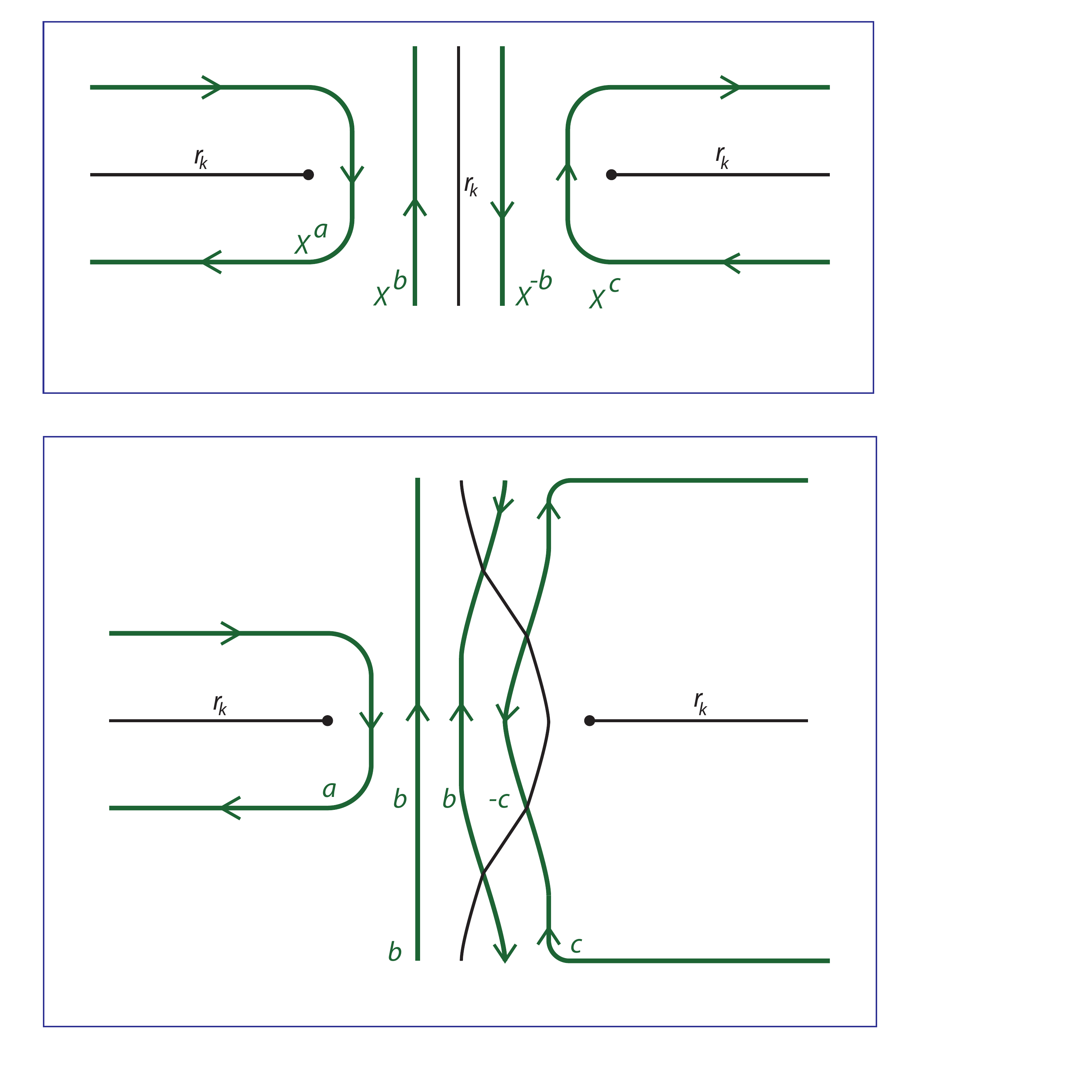}}
\]
\caption{The local view of a dihedral colored crossing}
\label{dihcrossing}
\end{figure}

 \begin{figure}[htb]
\[
\scalebox{.25}{\includegraphics{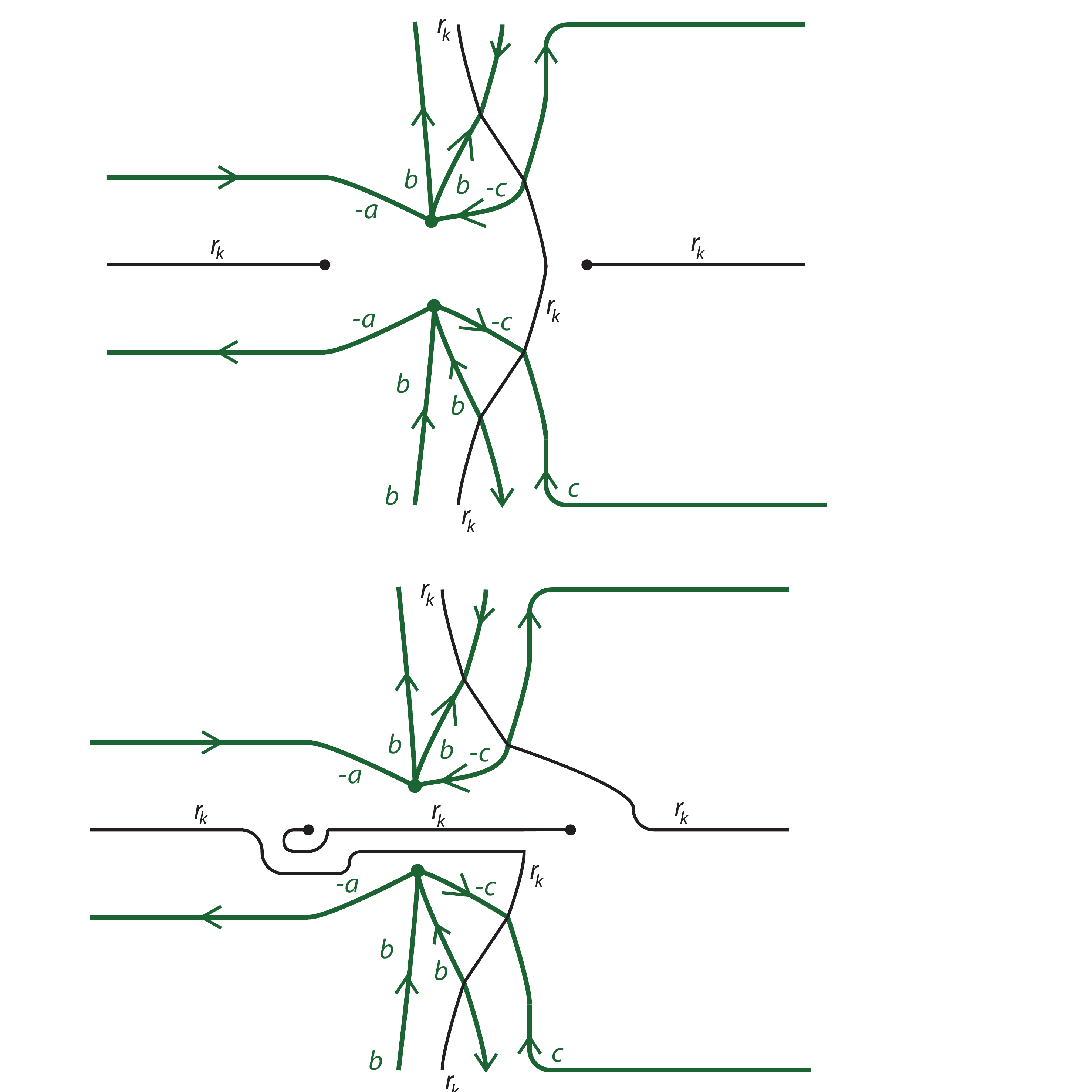}}
\]
\caption{The essential step in proving Theorem~\ref{DihedralThm}}
\label{proofofthm2}
\end{figure}

\clearpage

Recall from Section~\ref{coolcat}, the category ${\mathcal C}{(D_n)}$, in which the objects are words in the generators $r$ and $x$, and the morphisms are generated by the relators $r^2$, $x^n$,   $rxrx$, and the trivial relators $r^{\pm 1}r^{\mp 1}$ and $x^{\pm 1}x^{\mp 1}$. In this category, the generating morphisms are all invertible. Thus the analogues of the chart moves presented in Fig.~\ref{dichartmoves} hold.

Consider the local description of a coloring of a crossing in the chart $C(K)_1$ as depicted at the top of Fig.~\ref{dihcrossing}. 
The under crossing arc on the left is colored by $x^{-a} r_k x^{a}$. The central over-crossing arc is colored  $x^{b} r_k x^{-b}$.  The right-hand-side of the figure indicates the under crossing arc that is colored $x^{-c} r_k x^{c}$. The relation $c=2b-a$ holds among the powers of the surrounding rotations. At the bottom of the Figure, the dihedral relations have been applied. For notational brevity, the symbol $x^j$ has been replaced by $j$, for $j=a,b,c.$  
The vertices in the illustration also are a short-hand depiction of multiple applications of the relation $xr= r^{-1}x^{-1}$.

The Fig.~\ref{proofofthm2} contains the next level of detail for the proof. Since $2b-a=c$, the rotations that separate the end points of the arcs colored $r_k$ can be moved out of the way by multiple applications of the invertibility of the morphism $(rx)^2$, and then broken by applying invertibility of $x^n$. Now the black vertices are separated by a black arc, and some saddle moves can interconnect the two black vertices as in the bottom of the illustration. 

The  process outlined in Figs.~\ref{dihcrossing} and~\ref{proofofthm2} can occur at each crossing in the diagram, and the chart can be reduced to one with no black vertices.

We hope that it is not too difficult to observe that the analogous  diagrammatic relation to Corollary~\ref{cor1} holds for general dihedral charts. In this way the dihedral black vertices can be joined via the moves given in Fig.~\ref{proofofthm2}. The short black arc in the bottom of this illustration corresponds to the under-crossing arc at the corresponding crossing. 

In this way, a dihedral chart that corresponds to the  Fox coloring of the knot $K$ can be replaced by a chart that has no black vertices incident to arcs labeled by $r_k$. 

Furthermore, a permutation chart can be constructed for the original dihedral chart   $C(K)_1$ (and its predessors) that corresponds to the Fox coloring of the knot $K$. As the degree $n$ of the dihedral group $D_n$ increases, the proofs of the analogues of  Lemmas~\ref{MainLemma1} and ~\ref{MainLemma2} become more and more baroque. 

In the same fashion, an arc that is labeled by $r_k$ represents a permutation chart such as that given in Fig.~\ref{squashedovalnest7}. The perturbation of these radiant branch points to simple branch points is depicted in Fig.~\ref{ovalnest7}. The length of the rotation $x \in D_n$, when written in terms of the $t_i$\/s is, of course, $n-1$. So the ovals that surround $r_k$ in the permutation representation of the initial chart $C(k)_1$ might involve many loops.  Thus it could take a while to remove the loops of intersection among the spheres from $C(K)_{0}$ to an empty chart. For this reason, it is best to reduce the powers of the rotation to representative words in $\Sigma_n$ that are reduced.

On the other hand, the process of going from $C(K)_0$ to $C(K)_1$ just corresponds to passing the critical points that are the births of the arcs in the diagram.

When all the black vertices of the dihedral chart have been removed, the chart represents a family of $n$ mutually intersecting $2$-spheres. These can be brought to a position of being concentric by means of the analogues to the Roseman moves that do not involve branch points.

Thus the proof of Theorem~\ref{DihedralThm} is complete.

 \begin{figure}[htb]
\[
\scalebox{.1}{\includegraphics{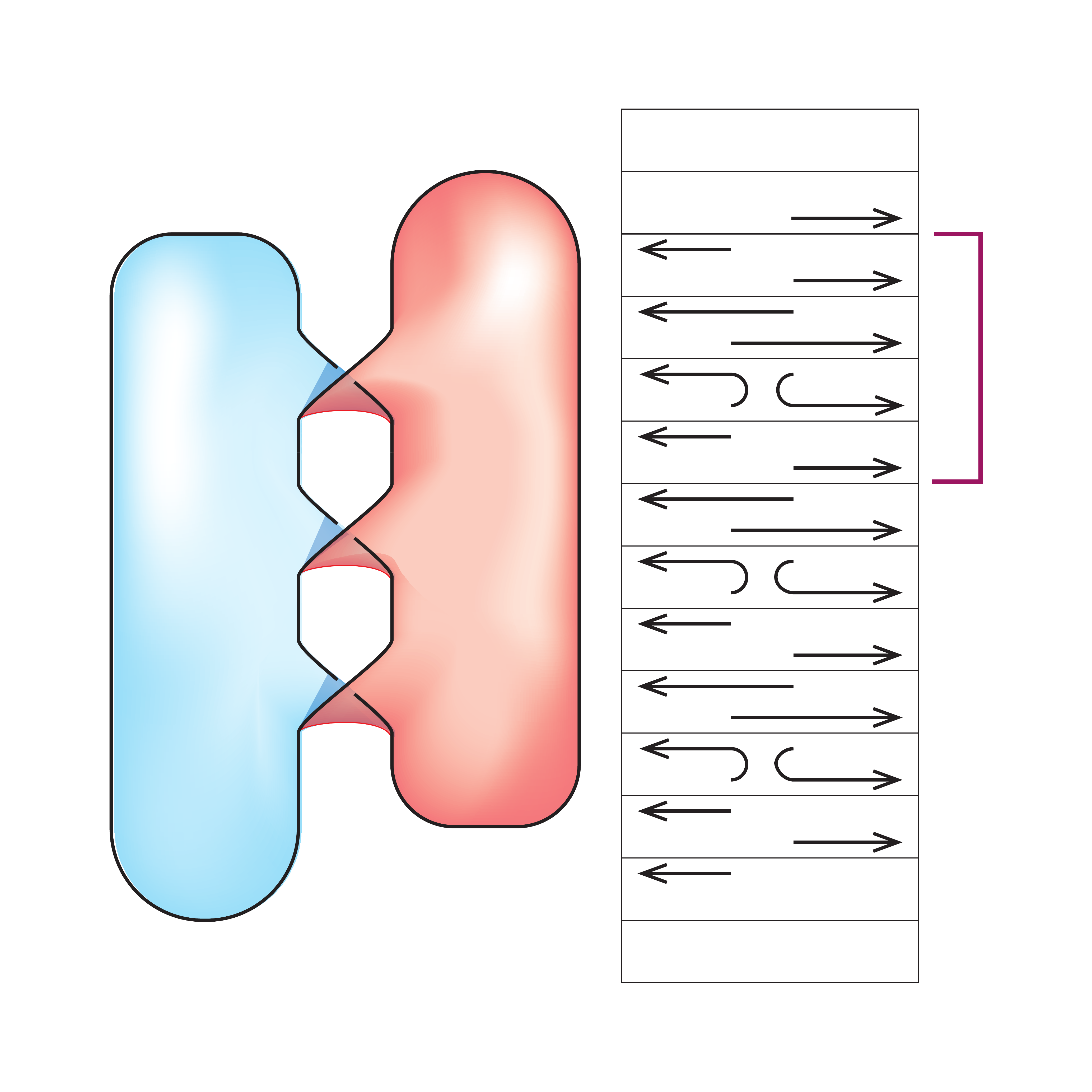}}
\]
\caption{A Seifert surface for the trefoil, its generic cross sections, and an indication of the first crossing}
\label{Seifert4Trefoil}
\end{figure}

\clearpage

\section{Proof of Theorem~\ref{CyclicThm}}
\label{cyclicproof}

Theorem~\ref{CyclicThm}   states that the cyclic branched cover of $S^3$ when branched over a knot or link $K\subset S^3$ can be folded and embedded into  the $5$-dimensional space $S^3\times [0,n+1]\times [-1,1]$. 
The construction upon which the proof relies is easy to extrapolate from a favored example. 

The $5$-fold cyclic branched cover of $S^3$ that is branched over the trefoil knot $K=3_1$ is known \cite{KirSch,Rolfsen}  to be the Poincar\'{e} homology sphere.  Such a cyclic branched cover of $S^3$ that is branched over a knot or link is constructed by taking multiple copies of $S^3$, cutting each along a Seifert surface of the knot, and gluing one to the other back to front. To construct the  braiding we cut the  $3$-sphere by a family of latitudinal $2$-spheres. Each sphere is at a fixed height of the knot diagram. The generic intersection of such a sphere with the Seifert surface will be used to construct a braid chart with non-simple, yet cyclic, branch points. The changes between successive generic intersections will indicate how to interpolate between successive charts.

\begin{figure}[htb]
\[
\scalebox{.1}{\includegraphics{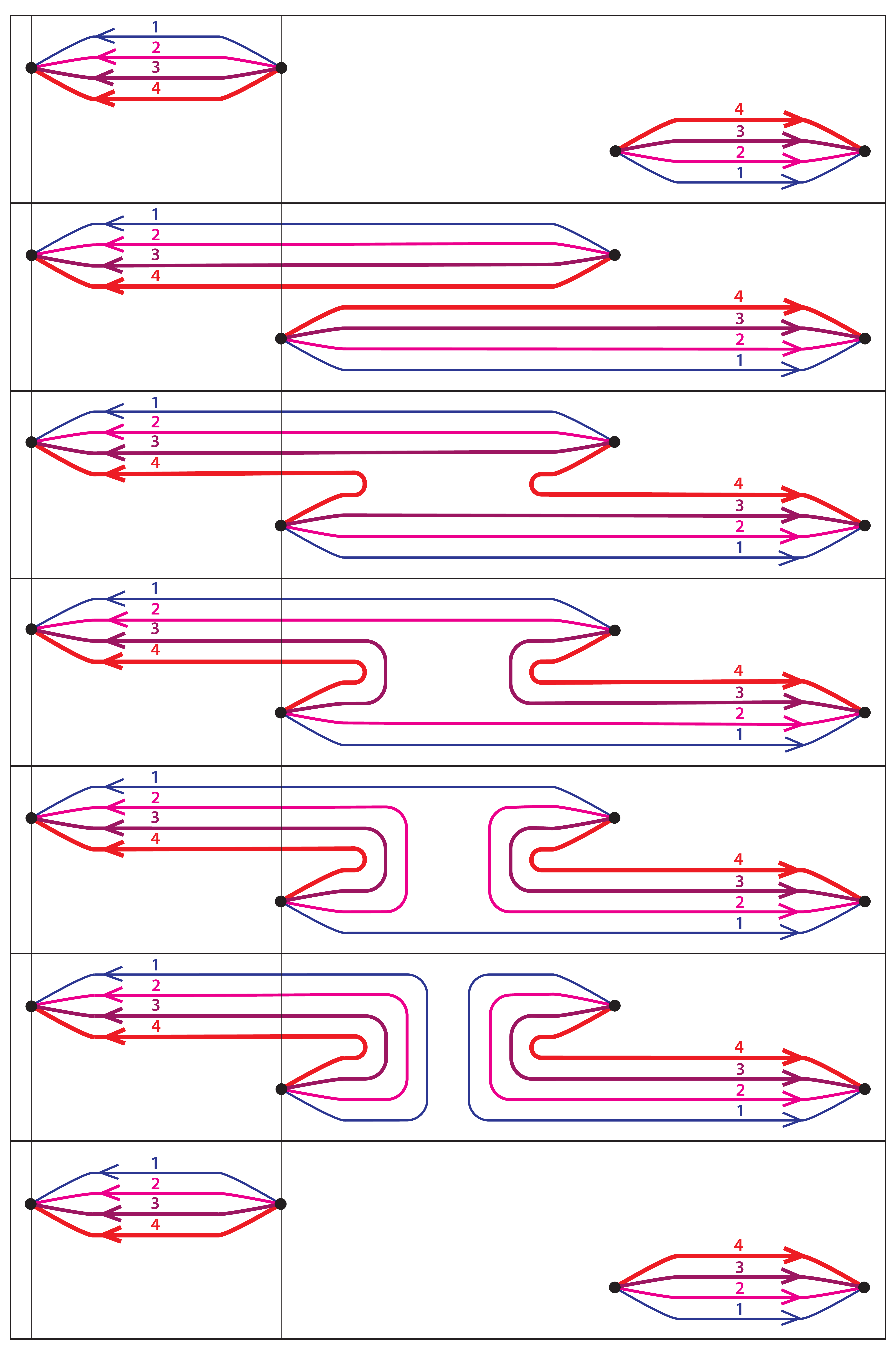}}
\]
\caption{The transitions of the $5$-fold cyclic cover at a crossing of the trefoil knot}
\label{fiveoldatacrossing}
\end{figure}

In Fig.~\ref{Seifert4Trefoil}, an oriented Seifert surface for the trefoil is indicated. On the right of the illustration, generic cross-sections are depicted, and the bracket on the far right demonstrates the cross-sectional twisted band that occurs along the crossing. The main part of the construction of the braiding occurs in this area. The slices of the $5$-fold cyclic branched cover of the trefoil that are associated such a twisted band are illustrated in Fig.~\ref{fiveoldatacrossing}. Replicate these slices three times. In this way, the most interesting part of the branched covering will be folded. To complete the braiding, allow two  pairs of cyclic branch points to be born as in the transition from the bottom drawing in Fig.~\ref{fromnone2five} to the top drawing therein. Similarly, allow them to disappear as if  going from top-to-bottom. 

The example here generalizes in two ways. First an arbitrary knot or link can be considered. The Seifert surface for such a knot is constructed similarly to the one illustrated in Fig.~\ref{Seifert4Trefoil} via the Seifert algorithm. Indeed, twisted bands interconnect the Seifert circles. At such twisted bands oriented versions of moves from Fig.~\ref{fiveoldatacrossing} can be implemented. The saddle moves between levels are found in Fig.~\ref{almostKmoves} and these are the immersed versions of Roseman type-II saddle move. 

The next generalization occurs when the $n$-fold cyclic branched cover involves $n\ne 5$. We imagine it is easy to abstract from the example given.

Finally, we point out that the edges in the permutation charts for the cyclic coverings can be oriented by using the orientations of the knots or links that are the branched sets: If an arc in the diagram points upward, then the edges in the chart point towards that branch point. The saddle moves will all be orientation preserving. 

This completes the proof.

\section{Conclusion and future work}
\label{Concl}

In this paper, we have outlined and provided much detail on how to construct foldings and braidings of branched covers of the $2$-dimensional and $3$-dimensional spheres. For specificity, we concentrated upon dihedral covers since the dihedral groups $D_n$ appear nicely as subsets of the permutation groups $\Sigma_n$. In particular, permutation charts are obtained from braid charts by removing orientations. Since braid charts are used to represent surface braids which are embedded surfaces in $S^2\times [0, n+1] \times [-1,1]$, permutation charts can be used to represent folded surfaces that are generically mapped to $S^2\times [0, n+1]$ away from the branch set. 

The ideas of braid, permutation, and dihedral charts can be abstracted to any group, $G$,  that has a given finite presentation. If the group $G$ is finite, it certainly has a finite permutation representation. Of course, the group acts upon itself, but as in the case of the dihedral group, we are interested in presentations that arise from a geometric action of the group. Theorems~\ref{Dihedral2} and~\ref{DihedralThm} appear to generalize easily.

In particular, the simplification of the permutation charts after the branch points have been removed depends upon the triviality of the higher homotopy groups in $B\Sigma_n$. Yet the higher dimensional cells in this classifying space are characterized as products of permutahedra. Simplifications of the immersed spheres can be achieved by the projections of the Roseman moves.

 Let $K=3_1$ denote the trefoil knot. There are  specific homomorphisms $A_3 \stackrel{\phi_3}{\longleftarrow} \pi_1(E(3_1))$ and 
$A_5 \stackrel{\phi_5}{\longleftarrow} \pi_1(E(3_1))$ in which the all the arcs in the standard diagram of the trefoil are assigned $3$-cycles or $5$-cycles, respectively. 
The later case was given by Fox in \cite{QuickTrip}. The former case corresponds to a quandle coloring of the trefoil by the four element tetrahedral quandle $Q(4,1)$. In both these cases,  our techniques can yield foldings of the associated branched coverings.  To convince yourself that the generalizations of 
Theorems~\ref{Dihedral2} and~\ref{DihedralThm} to arbitrary finite groups are possible, we encourage you to construct foldings that are associated to these colorings. 

Given a folding of a covering into $S^3 \times [0,n+1]$ it is not too difficult to find an immersed braiding into $S^3 \times [0,n+1] \times [-1,1].$ One attempts to consistently orient the edges in the sequence of permutation charts that describe the folding. Some choices of orientations may not pass easily through white vertices (triple points) of a chart, or the orientation may not be able to pass consistently between successive charts. The keys to consistent orientations are found within the possible proofs of the corresponding lemmas above. Thus foldings in $S^3 \times [0,n+1]$ provide possible examples of $3$-manifolds in $4$-space that do not lift to embeddings in $5$-space. Our examples provide a proving ground for liftability criteria.

We are hopeful that these folded maps of branched covers might be useful for computing invariants of knots and links that are derived from the coverings.

\subsection*{Acknowledgements.} JSC gratefully acknowledges a conversation with John Etnyre many years ago which started his interest in the subject. He would like to thank Patricia Cahn, Cameron Gordon, Lou Kauffman, Allison Miller, and Maggie Miller for helpful conversations and hints. Some of the work was discussed at the 2024 trisectors workshop in Lincoln, Nebraska. The work of Seonmi Choi was supported by the Basic Science Research Program through the National Research Foundation of Korea (NRF) funded by the Ministry of Education (No. 2021R1I1A1A01049100).
The work of Byoerhi Kim was supported by National Research Foundation of Korea (NRF) Grant No. 2019R1A3B2067839 and No. 2022R1A6A3A01086872.

\bigskip
\footnotesize
\bibliographystyle{abbrv}
\bibliography{biblio}

\end{document}